\documentclass[a4paper,12pt]{amsart}
\usepackage{amssymb}
\input xypic
\xyoption{all}
\xyoption{arc}

\newtheorem{theo}{Theorem}[section]

\newtheorem{prop}[theo]{Proposition}

\newtheorem{lem}[theo]{Lemma}

\newtheorem{coro}[theo]{Corollary}

\def\remark#1{{\refstepcounter{theo}\label{#1}\noindent\sc Remark  
\arabic{section}.\arabic{theo} - }}
\def\example#1{{\refstepcounter{theo}\label{#1}\noindent\sc Example 
\arabic{section}.\arabic{theo} - }}

\def\equat{\refstepcounter{theo}$$~}
\def\endequat{\leqno{{\boldsymbol{(}}
{\boldsymbol{\arabic{section}.\arabic{theo})}}}~$$}

\newcounter{numero}[section]

\def\FM{{\mathbb{F}}}

\def\NM{{\mathbb{N}}}

\def\QM{{\mathbb{Q}}}

\def\ZM{{\mathbb{Z}}}

\def\SG{{\mathfrak S}}

\def\bGB{{\boldsymbol{\mathfrak b}}}

\def\gGB{{\boldsymbol{\mathfrak g}}}
\def\hGB{{\boldsymbol{\mathfrak h}}}

\def\lGB{{\boldsymbol{\mathfrak l}}}

\def\pGB{{\boldsymbol{\mathfrak p}}}

\def\tGB{{\boldsymbol{\mathfrak t}}}
\def\uGB{{\boldsymbol{\mathfrak u}}}
\def\vGB{{\boldsymbol{\mathfrak v}}}

\def\a{\alpha}
\def\b{\beta}

\def\g{\gamma}

\def\d{\delta}
\def\D{\Delta}
\def\e{\varepsilon}
\def\ph{\varphi}
\def\ch{\chi}
\def\l{\lambda}

\def\m{\mu}
\def\n{\nu}
\def\o{\omega}
\def\O{\Omega}
\def\r{\rho}
\def\s{\sigma}

\def\th{\theta}
\def\Th{\Theta}
\def\t{\tau}
\def\x{\xi}

\def\z{\zeta}

\def\alpb{{\boldsymbol{\alpha}}}

\def\mub{{\boldsymbol{\mu}}}

\def\Sigb{{\boldsymbol{\Sig}}}

\def\AC{{\mathcal{A}}}
\def\BC{{\mathcal{B}}}

\def\DC{{\mathcal{D}}}
\def\EC{{\mathcal{E}}}
\def\FC{{\mathcal{F}}}

\def\HC{{\mathcal{H}}}

\def\KC{{\mathcal{K}}}
\def\LC{{\mathcal{L}}}
\def\MC{{\mathcal{M}}}

\def\UC{{\mathcal{U}}}

\def\WC{{\mathcal{W}}}
\def\XC{{\mathcal{X}}}

\def\ZC{{\mathcal{Z}}}

\def\FCt{{\tilde{\mathcal{F}}}}

\def\KCt{{\tilde{\mathcal{K}}}}

\def\FCh{{\hat{\mathcal{F}}}}

\def\NCB{{\boldsymbol{\mathcal{N}}}}
\def\OCB{{\boldsymbol{\mathcal{O}}}}

\def\SCB{{\boldsymbol{\mathcal{S}}}}

\def\Ab{{\mathbf A}}
\def\Bb{{\mathbf B}}
\def\Cb{{\mathbf C}}

\def\Fb{{\mathbf F}}
\def\Gb{{\mathbf G}}
\def\Hb{{\mathbf H}}

\def\Lb{{\mathbf L}}
\def\Mb{{\mathbf M}}

\def\Ob{{\mathbf O}}
\def\Pb{{\mathbf P}}
\def\Qb{{\mathbf Q}}

\def\Sb{{\mathbf S}}
\def\Tb{{\mathbf T}}
\def\Ub{{\mathbf U}}
\def\Vb{{\mathbf V}}

\def\Xb{{\mathbf X}}
\def\Yb{{\mathbf Y}}
\def\Zb{{\mathbf Z}}

\def\ab{{\mathbf a}}

\def\db{{\mathbf d}}

\def\fb{{\mathbf f}}

\def\hb{{\mathbf h}}
\def\ib{{\mathbf i}}

\def\nb{{\mathbf n}}

\def\pb{{\mathbf p}}

\def\sb{{\mathbf s}}

\def\vb{{\mathbf v}}

\def\xb{{\mathbf x}}

\def\fti{{\tilde{f}}}
\def\gti{{\tilde{g}}}

\def\uti{{\tilde{u}}}

\def\Kti{{\tilde{K}}}

\def\aha{{\hat{a}}}

\def\fha{{\hat{f}}}

\def\cdo{{\dot{c}}}

\def\sdo{{\dot{s}}}

\def\wdo{{\dot{w}}}

\def\zdo{{\dot{z}}}

\def\fba{{\bar{f}}}

\def\sba{{\bar{s}}}

\def\wba{{\bar{w}}}

\def\xbt{{\tilde{\xb}}}

\def\Fbt{{\tilde{\Fb}}}

\def\Obt{{\tilde{\Ob}}}

\def\Xbt{{\tilde{\Xb}}}
\def\Ybt{{\tilde{\Yb}}}

\def\Gbh{{\hat{\Gb}}}

\def\Xbh{{\hat{\Xb}}}
\def\Ybh{{\hat{\Yb}}}

\def\alpb{{\boldsymbol{\alpha}}}

\def\mub{{\boldsymbol{\mu}}}

\def\Sigb{{\boldsymbol{\Sigma}}}

\def\chit{{\tilde{\chi}}}

\def\sigt{{\tilde{\s}}}

\def\thet{{\tilde{\theta}}}

\def\zett{{\tilde{\zeta}}}

\def\alpba{{\bar{\alpha}}}
\def\betba{{\bar{\beta}}}

\def\piba{{\bar{\pi}}}

\def\rhoba{{\bar{\rho}}}

\def\gamh{{\hat{\gamma}}}

\def\Gambh{{\hat{\boldsymbol{\Gamma}}}}

\def\ad{\mathop{\mathrm{ad}}\nolimits}
\def\Ad{\mathop{\mathrm{Ad}}\nolimits}

\def\deg{\mathop{\mathrm{deg}}\nolimits}
\def\diag{\mathop{\mathrm{diag}}\nolimits}

\def\End{\mathop{\mathrm{End}}\nolimits}

\def\Id{\mathop{\mathrm{Id}}\nolimits}

\def\Im{\mathop{\mathrm{Im}}\nolimits}

\def\Ind{\mathop{\mathrm{Ind}}\nolimits}

\def\INT{\mathop{\mathrm{int}}\nolimits}
\def\Irr{\mathop{\mathrm{Irr}}\nolimits}

\def\Ker{\mathop{\mathrm{Ker}}\nolimits}

\def\reg{{\mathrm{reg}}}
\def\Reg{{\mathrm{Reg}}}
\def\res{\mathop{\mathrm{res}}\nolimits}
\def\Res{\mathop{\mathrm{Res}}\nolimits}

\def\Stab{\mathop{\mathrm{Stab}}\nolimits}

\def\supp{\mathop{\mathrm{supp}}\nolimits}

\def\Tr{\mathop{\mathrm{Tr}}\nolimits}

\def\uni{{\mathrm{uni}}}

\def\tete#1{\par\leavevmode\makebox[0.7cm]{$(\mathrm{#1})$}}

\def\to{\rightarrow}
\def\longto{\longrightarrow}
\def\injto{\hookrightarrow}

\def\longmapright#1{\hspace{0.3em}\smash{
     \mathop{\longrightarrow}\limits^{#1}}\hspace{0.3em}}

\def\fonction#1#2#3#4#5{\begin{array}{rccc}
{#1} : & {#2} & \longto & {#3} \\
& {#4} & \longmapsto & {#5} 
\end{array}}

\def\fonctio#1#2#3#4{\begin{array}{ccc}
{#1} & \longto & {#2} \\
{#3} & \longmapsto & {#4} 
\end{array}}

\def\ci{\circ}
\def\pr{\prime}
\def\ve{\vee}
\def\we{\wedge}

\def\incl{\hspace{0.05cm}{\subset}\hspace{0.05cm}}
\def\notincl{\hspace{0.05cm}{\not\subset}\hspace{0.05cm}}

\def\vide{\varnothing}

\def\fq{\FM_q}

\def\ql{{\QM_\el}}
\def\qlb{{\overline{\QM}_\el}}

\def\DS{\displaystyle}
\def\SS{\scriptstyle}
\def\SSS{\scriptscriptstyle}

\def\fin{~$\SS \blacksquare$}
\def\finl{~$\SS \square$}
\def\el{\ell}

\def\matrice#1{\left(\begin{array}{ccccccccccccccccccc}#1\end{array}\right)}
\def\matricecrochet#1{\left[\begin{array}{ccccccccccccccccccc}#1\end{array}\right]}

\def\lexp#1#2{\kern\scriptspace\vphantom{#2}^{#1}\kern-\scriptspace#2}
\def\le{\hspace{0.1em}\mathop{\leqslant}\nolimits\hspace{0.1em}}
\def\ge{\hspace{0.1em}\mathop{\geqslant}\nolimits\hspace{0.1em}}

\mathchardef\lllllll="3278
\def\SEC{$\lllllll$}
\mathchardef\inferieur="321E
\mathchardef\superieur="321F
\def\arobas{\char'100}

\def\sec{\section}
\def\sub{\subsection}

\def\bi{\bigskip}
\def\med{\medskip}
\def\sma{\smallskip}
\def\eqna{\begin{eqnarray*}}
\def\endeqna{\end{eqnarray*}}

\def\mor{morphism }

\def\endo{endomorphism }

\def\iso{isomorphism }
\def\isos{isomorphisms }
\def\auto{automorphism }

\def\proof{\noindent{\sc{Proof}~-} }
\def\rem{\noindent{\sc{Remark}~-} }

\def\tor{maximal torus }

\def\borel{Borel subgroup }

\def\para{parabolic subgroup }
\def\paras{parabolic subgroups }
\def\levi{Levi subgroup }
\def\levis{Levi subgroups }

\def\car{character }

\def\irr{irreducible }

\def\resp{respectively }

\def\cf{{\it cf.}~}

\def\na{\nabla}
\def\unb{{\mathbf{1}}}
\def\ins{{\mathrm{ins}}}

\def\mini{{\mathrm{min}}}

\def\et{{\mathrm{\acute{e}t}}}
\def\nablat{{\tilde{\nabla}}}
\def\phan{\phantom{\bullet}}
\def\nil{{\mathrm{nil}}}

\begin{document}

\begin{centerline}{\Large \bf Actions of relative Weyl groups I}\end{centerline}

\bi

\begin{centerline}{\sc C\'edric Bonnaf\'e
\footnote{{\sc CNRS - UMR 6623, Universit\'e de Franche-Comt\'e, 
D\'epartement de Math\'ematiques, 16 Route de Gray, 25000 BESAN\c{C}ON - FRANCE}, 
{\tt bonnafe\arobas math.univ-fcomte.fr}}}\end{centerline}

\bi

\begin{centerline}{\today}\end{centerline}

\bi

\bi

\begin{quotation}{\small \noindent {\sc Abstract - } 
We construct a new isomorphism between the endomorphism algebra 
of an induced cuspidal character sheaf and the group algebra 
of the relative Weyl group involved. We show it differs from 
the isomorphism of Lusztig by a linear character, and we relate 
this linear character to some 
stabilizers. Some consequences for characteristic functions of 
character sheaves are obtained. In the forthcoming second part, we 
will compute explicitly this linear character whenever the cuspidal 
local system is supported by the regular unipotent class and, 
as an application of these methods, we obtain a refinement of 
Digne, Lehrer 
and Michel's theorem on Lusztig restriction of Gel'fand-Graev 
characters.}\end{quotation}

\bigskip

\begin{quotation}
\end{quotation}

\bigskip

\begin{centerline}{\sc General introduction}\end{centerline}

\bi

A theorem of Digne, Lehrer and Michel says that the Lusztig restriction of 
a Gel'fand-Graev \car of a finite reductive group $\Gb^F$ is still a Gel'fand-Graev 
character \cite[Theorem 3.7]{DLM2}. However, 
an ambiguity remains on the character obtained (whenever the center of $\Gb$ 
is not connected, there are several Gel'fand-Graev characters). 
The original aim of this series of two papers was to resolve this ambiguity. 
For this, we needed to study more deeply the structure of the 
endomorphism algebra of an induced cuspidal character sheaf~: 
for instance, we wanted to follow the action of a Frobenius 
endomorphism through this algebra.

This led us to this first part, in which we develop another approach 
for computing explicitly this endomorphism algebra. One of the main goals 
is to construct another isomorphism between this endomorphism algebra and 
the group algebra of the 
relative Weyl group involved (one was already been constructed by 
Lusztig \cite[Theorem 9.2]{luicc}). In the case where $\Gb$ is symplectic 
or special orthogonal, this new isomorphism was constructed 
and computed explicitly by J.L. Waldspurger \cite[\SEC VIII.8]{waldspurger}. 

By comparing the isomorphisms, 
we get some immediate consequences for finite reductive groups. 
Note that the results of this part are valid for any 
cuspidal local system supported by a unipotent class and have 
a chance to be useful for computing values of characters 
at unipotent elements. 

In the second forthcoming part, we will restrict our attention to the case where 
the cuspidal local system is supported by the regular unipotent 
class. We are then able to compute explicitly the generalized 
Springer correspondence through this new isomorphism. This result 
is valid only for $p$ good. As an application 
of these (sometimes fastidious) computations, we get the desired more 
precise version of Digne, Lehrer and Michel's theorem. 
It must be said that this final result is valid only when the cardinality 
of the finite field is large enough.

\bi

\begin{centerline}{\sc Introduction to the first part}\end{centerline}

\bi

Let $\Gb$ be a connected reductive group defined over an algebraically 
closed field $\FM$, let $\Lb$ be a \levi of a \para $\Pb=\Lb.\Vb$ 
of $\Gb$, let $\Cb$ 
be a unipotent class of $\Lb$, and let $v \in \Cb$. We first explain how the action of 
the finite group $W_\Gb(\Lb,\Cb)=N_\Gb(\Lb,\Cb)/\Lb$ on some varieties 
introduced by Lusztig \cite[\SEC 3 and 4]{luicc} can be extended 
``by density'' to some slightly bigger varieties (see \SEC \ref{sub relatif}). 
We then generalize this construction to extend the action of 
$W_\Gb(\Lb,v) =N_\Gb(\Lb,v)/C_\Lb^\ci(v)$ to other varieties 
covering the previous ones (see \SEC \ref{sub nor}). 
One of our goals is to determine the stabilizer $H$ of an element 
lying over a representative $u$ of the induced unipotent class of $\Cb$ 
(we choose $u$ in $v\Vb$). 

Whenever $C_\Lb(v)/C_\Lb^\ci(v) \simeq C_\Gb(v)/C_\Gb^\ci(v)$, then 
$W_\Gb(\Lb,v)=W_\Gb^\ci(\Lb,v) \times A_\Lb(v)$ (where $A_\Lb(v)=C_\Lb(v)/C_\Lb^\ci(v)$ 
and $W_\Gb^\ci(\Lb,v)=N_\Gb(\Lb) \cap C_\Gb^\ci(v)/C_\Lb^\ci(v)$). 
If moreover $C_\Gb(u) \incl \Pb$, then we show that there exists a morphism 
$\ph_{\Lb,v}^\Gb : W_\Gb^\circ(\Lb,v) \to A_\Lb(v)=C_\Lb(v)/C_\Lb^\ci(v)$ such that
$$H=\{(w,a) \in W_\Gb^\ci(\Lb,v) \times A_\Lb(v)~|~a=\ph_{\Lb,v}^\Gb(w)\}$$
(see \SEC \ref{02} and \ref{01}). We also provide some reduction 
arguments to compute explicitly the morphisms $\ph_{\Lb,v}^\Gb$ 
(see \SEC \ref{ele sec}).

From \SEC \ref{sec endo} to the end, we assume that $\Cb$ supports 
a cuspidal local system $\EC$ (we denote by $\z$ the 
\car of the finite group $A_\Lb(v)$ associated to $\EC$). Let 
$K$ denote the perverse sheaf obtained from the datum $(\Cb,\EC)$ 
by parabolic induction \cite[4.1.1]{luicc}, and let $\AC$ denote 
its endomorphism algebra. In this case, $W_\Gb(\Lb,\Cb)$ is 
equal to $W_\Gb(\Lb)$ (by \cite[Theorem 9.2]{luicc}) and is 
isomorphic to $W_\Gb^\ci(\Lb,v)$. 
Lusztig \cite[Theorem 9.2]{luicc} 
constructed a canonical isomorphism 
$\Th : \qlb W_\Gb^\ci(\Lb,v) \to \AC$. The aim of \SEC \ref{sec endo} 
is to construct another explicit isomorphism $\Th' : \qlb W_\Gb^\ci(\Lb,v) \to \AC$ 
using the varieties previously introduced in \SEC \ref{mor phi}. 
It turns out that $\Th$ and $\Th'$ differ by a linear character 
$\g_{\Lb,v,\z}^\Gb$ (see Corollary \ref{epsilon}). 
For computing explicitly this linear character, one could use the 
characterization of Lusztig in terms of the action on the perverse sheaf 
$K$. We skip this difficulty by using work done in the previous 
sections~: we show that the linear character 
$\g_{\Lb,v,\z}^\Gb$ 
is known whenever the morphism $\ph_{\Lb,v}^\Gb$ is known. (Indeed, 
the fact that $\Cb$ supports a cuspidal local system implies that 
$A_\Lb(v)=A_\Gb(v)$ and $C_\Gb(u) \incl \Pb$, so that 
the morphism $\ph_{\Lb,v}^\Gb$ is defined~: see Theorem \ref{alv}.) They are related 
by the formula $\g_{\Lb,v,\z}^\Gb={1 \over \z(1)}(\z \ci \ph_{\Lb,v}^\Gb)$ 
(which means that, in this case, $\ph_{\Lb,v}^\Gb$ has values in the center 
of the character $\z$). 

In \SEC \ref{part finite}, we assume further that $\FM$ is an algebraic 
closure of a finite field and that $\Gb$ is endowed with a Frobenius 
endomorphism. We then explain what kind of refinements may be obtained 
by using the previous results about characteristic functions of 
character sheaves.

Whenever $\Gb$ is special orthogonal or symplectic, then the linear character 
$\g_{\Lb,v,\z}^\Gb$ has been computed directly (without using explicitly 
the morphism $\ph_{\Lb,v}^\Gb$) by J.L. Waldspurger \cite[Lemma VIII.9]{waldspurger}. 
In the future part II, we will assume throughout that $v$ is regular. 
Under this hypothesis, 
we will compute explicitly the morphisms $\ph_{\Lb,v}^\Gb$ 
even when $\Cb$ does not support a cuspidal local system. 
We then follow Digne, Lehrer and 
Michel's method~: knowing $\g_{\Lb,v,\z}^\Gb$ and from 
the explicit nature of the isomorphism $\Th'$, we get a slightly 
more precise version of their theorem on Lusztig restriction 
of Gel'fand-Graev characters.

\bi

\def\isol{{\mathrm{iso}}}

\section*{Notation}~

\med

\noindent{\bf Fields, varieties, sheaves.} We fix an algebraically 
closed field $\FM$ and we denote by $p$ its 
characteristic. All algebraic varieties and all algebraic groups will be 
considered over $\FM$. We also fix a prime number $\el$ different from $p$. 
Let $\qlb$ denote an algebraic 
closure of the $\el$-adic field $\ql$. 

If $\Xb$ is an algebraic variety (over $\FM$), we also denote 
by $\qlb$ the constant $\el$-adic sheaf associated to $\qlb$ (if necessary, 
we denote it by $(\qlb)_\Xb$). By a constructible 
sheaf (respectively a local system) on $\Xb$ we mean 
a constructible $\qlb$-sheaf (respectively a $\qlb$-local system).
Let $\DC\Xb$ denote the 
bounded derived category of constructible sheaves on $\Xb$. If $K \in \DC\Xb$ 
and $i \in \ZM$, we denote by $\HC^iK$ the $i$-th cohomology sheaf of $K$ 
and if $x \in \Xb$, then $\HC^i_x K$ denotes the stalk at $x$ of the 
constructible sheaf $\HC^iK$. If $K \in \DC \Xb$, 
we denote by $DK$ its Verdier dual. If $\LC$ is a constructible sheaf on $\Xb$, 
we identify it with its image in $\DC\Xb$, that is 
the complex, concentrated in degree $0$, whose $0$th term is $\LC$. 

Let $K \in \DC \Xb$. We say that $K$ is a {\it perverse sheaf} if the following two 
conditions hold~:
$${\mathrm{(a)}}\hspace{2.3cm}\forall i \in \ZM,~\dim \supp \HC^i K \le -i,$$
$${\mathrm{(b)}}\hspace{2cm}\forall i \in \ZM,~\dim \supp \HC^i D K \le -i.$$
We denote by $\MC\Xb$ the full subcategory of $\DC\Xb$ whose objects are perverse 
sheaves~: this is an abelian category \cite[2.14, 1.3.6]{BBD}.

Let $\Yb$ be a locally closed, smooth, irreducible subvariety of $\Xb$ and let 
$\LC$ be a local system on $\Yb$. We denote by $IC(\overline{\Yb},\LC)$ the 
Deligne-Goresky-MacPherson intersection cohomology complex of $\overline{\Yb}$ 
with coefficients in $\LC$. We often identify $IC(\overline{\Yb},\LC)$ with 
its extension by zero to $\Xb$~; $IC(\overline{\Yb},\LC)[\dim \Yb]$ is a perverse sheaf 
on $\Xb$. 

\bi

\noindent{\bf Algebraic groups.} If $\Hb$ is a linear algebraic group, we will denote 
by $\Hb^\ci$ the neutral component of $\Hb$, by $\Hb_\uni$ the closed 
subvariety of $\Hb$ consisting of unipotent elements of $\Hb$, and 
by $\Zb(\Hb)$ the center of $\Hb$. If $h \in \Hb$, then 
$A_\Hb(h)$ denotes the finite group $C_\Hb(h)/C_\Hb^\ci(h)$, $(h)_\Hb$ denotes the 
conjugacy class of $h$ in $\Hb$, and $h_s$ (\resp $h_u$) denotes 
the semisimple (\resp unipotent) part of $h$. If $\hGB$ is 
the Lie algebra of $\Hb$, we denote by $\Ad h : \hGB \to \hGB$ the 
differential at $1$ of the \auto $\Hb \to \Hb$, $x \mapsto \lexp{h}{x}=hxh^{-1}$.

If $\Xb$ and $\Yb$ are varieties, and if $\Xb$ (\resp $\Yb$) is 
endowed with an action of $\Hb$ on the right (\resp left), 
then we denote, when it exists, $\Xb \times_\Hb \Yb$ the quotient 
of $\Xb \times \Yb$ by the diagonal left action of $\Hb$ given by 
$h.(x,y)=(xh^{-1},hy)$ for any $h \in \Hb$ and $(x,y) \in \Xb \times \Yb$. 
If $(x,y) \in \Xb \times \Yb$, and if $\Xb \times_\Hb \Yb$ exists, 
we denote by $x *_\Hb y$ the image of $(x,y)$ in $\Xb \times_\Hb \Yb$ 
by the canonical morphism.

Finally, if $X_1$,\dots, $X_n$ are subsets or elements of $\Hb$, we denote 
by $N_\Hb(X_1,\dots,X_n)$ the intersection of the normalizers $N_\Hb(X_i)$ 
of $X_i$ in $\Hb$ ($1 \le i \le n$). 

\bi

\noindent{\bf Reductive group.} We fix once and for all a connected reductive algebraic 
group $\Gb$. We fix a \borel $\Bb$ of $\Gb$ and a maximal torus $\Tb$ of $\Bb$. 
We denote by $X(\Tb)$ (\resp $Y(\Tb)$) the lattice of rational characters 
(\resp of one-parameter subgroups) of $\Tb$. 
Let $W=N_\Gb(\Tb)/\Tb$. 
Let $\Phi$ denote the root system of $\Gb$ relative to $\Tb$ and let 
$\Phi^+$ (\resp $\D$) denote the set of positive roots 
(\resp the basis) of $\Phi$ associated to $\Bb$. 
For each root $\a \in \Phi$, we denote by 
$\Ub_\a$ the one-parameter unipotent subgroup of $\Gb$ 
normalized by $\Tb$ associated to $\a$. 

We also fix in this paper 
a \para $\Pb$ of $\Gb$ and a \levi $\Lb$ of $\Pb$. We denote by $\pi_\Lb : \Pb \to \Lb$ 
the canonical projection with kernel $\Vb$, the unipotent radical of $\Pb$. 
We denote by $\Phi_\Lb$ the root system 
of $\Lb$ relative to $\Tb$~; we have $\Phi_\Lb \incl \Phi$. Finally, $W_\Lb$ 
denotes the Weyl group of $\Lb$ relative to $\Tb$. 

\bi

\sec{Preliminaries\label{sec steinberg}}~

\med

\sub{Centralizers} We start this subsection by recalling two well-known 
results on centralizers of elements in reductive groups. The first one is 
due to Lusztig \cite[Proposition 1.2]{luicc}, while the second one is due 
to Spaltenstein \cite[Proposition 3]{holt}.

\bi

\begin{lem}[{\bf Lusztig}]\label{dim p}
$(1)$ Let $l \in \Lb$ and $g \in \Gb$. Then 
$$\dim \{x \Pb~|~x^{-1}gx \in (l)_\Lb .\Vb\} \le 
{1 \over 2}(\dim C_\Gb(g)-\dim C_\Lb(l)).$$

\tete{2} If $g \in \Pb$, then $\dim C_\Pb(g) \ge \dim C_\Lb(\pi_\Lb(g))$.
\end{lem}

\bi

\begin{lem}[{\bf Spaltenstein}]\label{spa connexe}
If $l \in \Lb$, then $C_\Vb(l)$ is connected.
\end{lem}

\bi

We will now give several applications of the two previous lemmas. 
We first need the following technical result~:

\bi

\begin{lem}\label{equiv class}
Let $l \in \Lb$. Then the following are equivalent~:

\tete{a} $C_\Gb^\ci(l) \incl \Lb$~;

\tete{b} $C_\Gb^\ci(l_s) \incl \Lb$~;

\tete{c} $C_\Vb(l)=\{1\}$~;

\tete{d} $\dim C_\Vb(l)=0$.
\end{lem}

\bi

\proof It is clear that (b) implies (a), and that (a) implies (d). 
Moreover, by Lemma \ref{spa connexe}, (c) is equivalent to (d). 
It remains to prove that (c) implies (b). 

For this, let $s$ (\resp $u$) denote the semisimple (\resp unipotent) 
part of $l$, and assume that $C_\Gb^\ci(s) \notincl \Lb$. We want 
to prove that $C_\Vb(l) \not= \{1\}$. Without loss of generality, 
we may and do assume that $s \in \Tb$ and $u \in \Ub \cap \Lb$. 

Let $\Gb'=C_\Gb^\ci(s)$, $\Ub'=C_\Ub(s)$, $\Vb'=C_\Vb(s)$ and $\Bb'=C_\Bb(s)$. 
Then, by \cite[Corollary 11.12]{borel}, $u \in \Gb'$. 
Moreover, by Lemma \ref{spa connexe}, $\Ub'$ and $\Vb'$ are connected. 
So $\Bb'=\Tb.\Ub'$ is connected. Let $\Phi_s$ denote the root system 
of $\Gb'$ relative to $\Tb$, and let $\Phi_s^+$ 
be the positive root system associated to 
the Borel subgroup $\Bb'$ of $\Gb'$. 
Since $\Gb' \notincl \Lb$, there exists an irreducible component $\Psi$ of 
the root system $\Phi_s$ such that $\Psi \notincl \Phi_\Lb$. 
Let $\a_s$ denote the highest root of $\Psi$ with respect to $\Psi \cap \Phi_s^+$~: 
then $\a_s \not\in \Phi_\Lb$. But $\Ub_{\a_s}$ is central in $\Ub'$, so 
$\Ub_{\a_s} \incl C_{\Vb'}(u)=C_\Vb(l)$. 
Therefore, $C_\Vb(l) \not= \{1\}$. The proof of Lemma 
\ref{equiv class} is now complete.\fin

\bi

Now, let $\OCB$ be the set of elements $l \in \Lb$ such that 
$C_\Vb(l)=\{1\}$. This variety has been used by Lusztig in 
\cite[\SEC 2]{luspringer} where it was denoted by $\UC$.

\bi

\begin{lem}\label{U conjugue}
The set $\OCB$ is a dense open subset of $\Lb$ and the map
$$\fonctio{\OCB \times \Vb}{\OCB.\Vb}{(l,x)}{xlx^{-1}}$$
is an isomorphism of varieties.
\end{lem}

\bi

\proof The group $\Vb$ acts on $\Gb$ by conjugation. So, by 
\cite[Proposition 1.4]{humphreys}, the set 
$$\NCB=\{g \in \Gb~|~\dim C_\Vb(g) =0 \}$$
is an open subset of $\Gb$. Therefore, $\NCB \cap \Lb$ is an open 
subset of $\Lb$. Moreover, $\NCB \cap \Lb$ is not empty 
since any $\Gb$-regular element of $\Tb$ belongs 
to $\NCB \cap \Lb$. But, by Lemma \ref{spa connexe}, $\NCB \cap \Lb=\OCB$. 
This proves the first assertion of the lemma.

Now, let $f$ denote the \mor defined in Lemma \ref{U conjugue}. 
Let $l \in \OCB$. Then the map $f_l : \Vb \to l.\Vb$, $x \mapsto xlx^{-1}$ 
is injective by definition of $\OCB$, and its image is closed because it 
is an orbit under a unipotent group \cite[Proposition 4.10]{borel}. 
By comparing dimensions, we get that $f_l$ is bijective. 
As this holds for every $l \in \OCB$, $f$ is bijective.

Moreover, the variety $\OCB.\Vb \simeq \OCB \times \Vb$ is smooth. 
Hence, to prove that $f$ is an isomorphism, it is enough to prove   
that the differential $(df)_{(l,1)}$ is surjective for 
some $l \in \OCB$ (see \cite[Theorems AG.17.3 and AG.18.2]{borel}). 

Now, let $t \in \Tb$ be a $\Gb$-regular element (so that $t \in \OCB$). 
The tangent space to $\OCB$ at $t$ may be identified with the Lie algebra 
$\lGB$ of $\Lb$ via the translation by $t$. By writing 
$f(l,x)=l.(l^{-1}xlx^{-1})$ for every $(l,x) \in \OCB \times \Vb$, 
the differential $(df)_{(t,1)}$ may be identified with the map 
$$\fonctio{\lGB \oplus \vGB}{\lGB \oplus \vGB}{l \oplus x}{l \oplus 
(\Ad t^{-1} - \Id_\vGB)(x),}$$
where $\vGB$ denotes the Lie algebra of $\Vb$. 
The bijectivity of $(df)_{(t,1)}$ follows immediately from the fact 
that the eigenvalues of $\Ad t^{-1}$ are equal to $\a(t)^{-1}$ for 
$\a \in \Phi^+-\Phi_\Lb$, so they are different from $1$ by 
the regularity of $t$.\fin

\bi

Lemma \ref{U conjugue} implies immediately the following~:

\bi

\begin{coro}\label{U conjugue bis}
Let $\SCB$ be a locally closed subvariety of $\OCB$. Then the map
$$\fonctio{\SCB \times \Vb}{\SCB.\Vb}{(l,x)}{xlx^{-1}}$$
is an \iso of varieties.
\end{coro}

\bi

\noindent{\sc Notation - } If $\SCB$ is a locally closed subvariety 
of $\Lb$, we denote by $\SCB_\reg$ (or $\SCB_{\reg,\Gb}$ if there is 
some ambiguity) the open subset $\SCB \cap \OCB$ of $\SCB$. It might 
be empty.\finl

\bi

\sub{Steinberg map} Let $\na : \Gb \to \Tb/W$ be the 
Steinberg map. Recall that for $g \in \Gb$, $\na(g)$ is defined to be 
the intersection of $\Tb$ with the conjugacy class 
of the semisimple part of $g$. Then $\na$ is a \mor of varieties 
\cite[\SEC{6}]{steinreg}. To compute the Steinberg map, we need to determine 
semisimple parts of elements of $\Gb$. In our situation, the following well-known 
lemma will be useful \cite[5.1]{luicc}~:

\bi

\begin{lem}\label{V conjugue}
If $g \in \Pb$, then the semisimple part of $g$ is 
$\Vb$-conjugate to the semisimple part of $\pi_\Lb(g)$. In particular, 
$\na(g)=\na(\pi_\Lb(g))$.
\end{lem}

\bi

\proof Let $s$ be the semisimple part of $g$. Then the semisimple part of 
$\pi_\Lb(g)$ is $\pi_\Lb(s)$. But, $s$ belongs to some \levi $\Lb_0$ of 
$\Pb$. Let $x \in \Vb$ be such that $\Lb=\lexp{x}{\Lb_0}$. 
Then $xsx^{-1} \in \Lb$ and $xsx^{-1}$ is the semisimple part of 
$xgx^{-1}$. Therefore, the semisimple part of 
$\pi_\Lb(xgx^{-1})=\pi_\Lb(g)$ is $\pi_\Lb(xsx^{-1})=xsx^{-1}$.\fin

\bi

\begin{lem}\label{nabla -2}
$\na(\Zb(\Lb)^\ci)$ is a closed subset of $\Tb/W$ and $\na(\Zb(\Lb)^\ci_\reg)$ 
is an open subset of $\na(\Zb(\Lb)^\ci)$. 
\end{lem}

\bi

\proof The restriction of $\na$ to $\Tb$ is a finite quotient morphism. 
In particular, it is open and closed. Since it is closed, $\na(\Zb(\Lb)^\ci)$ 
is a closed subset of $\Tb/W$. Since it is open, 
$\na(\Tb_\reg)$ is an open subset of $\Tb/W$. 
But $\na(\Zb(\Lb)^\ci_\reg)=\na(\Tb_\reg) \cap \na(\Zb(\Lb)^\ci)$. 
So the Lemma \ref{nabla -2} is proved.\fin

\bi 

Let $\na_\Lb : \Lb \to \Tb/W_\Lb$ denote the Steinberg map for the 
group $\Lb$. By Lemma \ref{equiv class}, we have
\equat\label{O T}
\OCB=\na_\Lb^{-1}(\Tb_\reg/W_\Lb).
\endequat

\bi

\sub{A family of morphisms\label{sub fam}} If $\SCB$ is a locally closed subvariety 
of $\Lb$ stable under conjugation by $\Lb$, then $\SCB.\Vb$ is a 
locally closed subvariety of $\Pb$ stable under conjugation by $\Pb$. 
We can therefore consider the quotients $\Gb \times_\Lb \SCB$ and 
$\Gb \times_\Pb \SCB.\Vb$. In this subsection, we will focus on the 
maps $\Gb \times_\Lb \SCB \to \Gb$, $g *_\Lb l \mapsto glg^{-1}$ 
and $\Gb \times_\Pb \SCB.\Vb \to \Gb$, $g *_\Pb x \mapsto gxg^{-1}$ 
which are well-defined morphisms of varieties. 

\bi

\remark{identification} If $\SCB$ is contained in $\OCB$, then 
the map $\Gb \times_\Lb \SCB \to \Gb \times_\Pb \SCB.\Vb$, 
$g *_\Lb l \mapsto g *_\Pb l$ is an isomorphism of varieties 
(by Corollary \ref{U conjugue bis}).\finl

\bi

The next result is well-known~:

\bi

\begin{lem}\label{proj mor}
The map $\Gb \times_\Pb \Pb \to \Gb$, $g *_\Pb x \mapsto gxg^{-1}$ is a projective 
surjective morphism of varieties. In particular, if $\Fb$ is a closed subvariety 
of $\Pb$ stable under conjugation by $\Pb$, then the map 
$\Gb *_\Pb \Fb \to \Gb$, $g *_\Pb x \mapsto gxg^{-1}$ is a projective 
morphism.
\end{lem}

\bi

\proof Let $\Xbt=\{(x,g\Pb) \in \Gb \times \Gb/\Pb~|~g^{-1}xg \in \Pb\}$. 
Then $\Xbt$ is a closed subvariety of $\Gb \times \Gb/\Pb$. Moreover, 
the variety $\Gb/\Pb$ is projective. Therefore, the projection 
$\pi : \Xbt \to \Gb$, $(x,g\Pb) \mapsto x$ is a projective morphism. 
Since every element of $\Gb$ is conjugate to an element of $\Bb$, 
$\pi$ is surjective.

But the maps $\Gb \times_\Pb \Pb \to \Xbt$, $g *_\Pb x \mapsto (gxg^{-1},g\Pb)$ 
and $\Xbt \to \Gb \times_\Pb \Pb$, $(x,g\Pb) \mapsto g *_\Pb g^{-1}xg$ 
are morphisms of varieties which are inverse of each other. 
Moreover, through these isomorphisms, the map constructed in Lemma \ref{proj mor} 
may be identified with $\pi$. The proof is now complete.\fin

\bi

The next result might be known but we have not come across it in the 
literature. 

\bi

\begin{lem}\label{lissable}
The morphisms of varieties 
$$\fonctio{\Gb \times_\Lb \OCB}{\Gb}{g*_\Lb l}{glg^{-1}}$$
$$\fonctio{\Gb \times_\Pb \OCB.\Vb}{\Gb}{g*_\Pb x}{gxg^{-1}}\leqno{\mathit{and}}$$
are \'etale. 
\end{lem}

\bi

\proof By Remark \ref{identification}, 
it is sufficient to prove that the morphism 
$$\fonction{f}{\Gb \times_\Pb \OCB.\Vb}{\Gb}{g *_\Pb x}{gxg^{-1}}$$
is \'etale. 

Since $\Gb \times_\Pb \OCB.\Vb$ and $\Gb$ are smooth varieties, 
$f$ is \'etale if and only if the 
differential of $f$ at any point of $\Gb \times_\Pb \OCB.\Vb$ 
is an isomorphism \cite[Proposition III.10.4]{hartshorne}. 
By $\Gb$-equivariance of the morphism $f$ ($\Gb$ acts on 
$\Gb \times_\Pb \OCB.\Vb$ by left translation on the first factor, 
and acts on $\Gb$ by conjugation), it is sufficient to 
prove that $(df)_{1*_\Pb x}$ is an isomorphism for every 
$x \in \OCB.\Vb$. 

For this, let $\Pb^-$ denote the \para of $\Gb$ opposed to 
$\Pb$ (with respect to $\Lb$), and let $\Vb^-$ denote its unipotent 
radical. Then $\Vb^- \times\OCB.\Vb$ is an open neighborhood of 
$1 *_\Pb x$ in $\Gb *_\Pb \OCB.\Vb$. Therefore, it is sufficient 
to prove that the differential of the map 
$$\fonction{f^-}{\Vb^- \times \OCB.\Vb}{\Gb}{(g,x)}{gxg^{-1}}$$
at $(1,x)$ is an isomorphism for every $x \in \OCB.\Vb$.

Let $\gGB$, $\vGB^-$, $\lGB$ and $\pGB$ denote the Lie 
algebras of $\Gb$, $\Vb^-$, $\Lb$ and $\Pb$ respectively. 
Since $\OCB$ is open in $\Lb$, we may identify the tangent space 
to $\OCB.\Vb$ at $x$ with $\pGB$ (using left translation by $x$). 
Similarly, we identify the tangent space to $\Gb$ 
at $x$ with $\gGB$ using left translation. Using these identifications, 
the differential of $f^-$ at $(1,x)$ may be identified with 
the map
$$\fonction{\d}{\vGB^- \oplus \pGB}{\gGB}{A \oplus B }{
(\ad x)^{-1}(A) - A + B .}$$
For dimension reasons, we only need to prove that $\d$ is injective. 

For this, let $\l \in Y(\Tb)$ be such that $\Lb=C_\Gb(\Im \l)$ and 
$$\Pb=\{g \in \Gb~|~\lim_{t \to 0} \l(t)g\l(t)^{-1}~{\mathrm{exists}}\}.$$
For the definition of $\lim_{t \to 0} \l(t)$, see \cite[Page 184]{DLM1}. 
We then define, for each $i \in \ZM$, 
$$\gGB(i)=\{X \in \gGB~|~(\ad \l(t))(X)=t^i X\}.$$
Then
$$\gGB=\mathop{\oplus}_{i \in \ZM}~ \gGB(i),$$
$$\pGB=\mathop{\oplus}_{i \ge 0}~ \gGB(i),$$
$$\vGB^-=\mathop{\oplus}_{i < 0}~ \gGB(i)\leqno{\mathrm{and}}$$
(see \cite[5.14]{DLM1}). For each $X \in \gGB$, we denote by $X_i$ its projection 
on $\gGB(i)$.

Now, let $l=\pi_\Lb(x)$. Then it is 
clear that we have, for any $i_0 \in \ZM$ and any $X \in \gGB(i_0)$, 
$$(\ad x)^{-1}(X) \in (\ad l)^{-1}(X) + 
(\mathop{\oplus}_{i > i_0} \gGB(i)).\leqno{(1)}$$

Now, let $A \oplus B \in \Ker \d$, and assume that $A \not= 0$. 
Then there exists $i_0 < 0$ minimal among all $i < 0$ such that 
$A_i \not= 0$. Then, by (1), the projection of $\d(A \oplus B)$ on 
$\gGB(i_0)$ is equal to $(\ad l)^{-1}(A_{i_0})-A_{i_0}$. But, 
$\d(A \oplus B)=0$, so $(\ad l)^{-1}(A_{i_0})-A_{i_0}=0$. 
Therefore, $C_\gGB(l) \notincl \lGB$. So 
$C_\gGB(s) \notincl \lGB$, where $s$ denotes the semisimple part 
of $l$. 
However, $l$ lies in $\OCB$, so its semisimple part $s$ also lies 
in $\OCB$ by Lemma \ref{equiv class}. But, by \cite[Proposition 9.1 (1)]{borel}, 
$C_\gGB(s)$ is the Lie algebra of $C_\Gb^\ci(s)$ which 
is contained in $\Lb$ by Lemma \ref{equiv class}. We get a contradiction.

So, this discussion shows that $A=0$. But $0=\d(A,B)=(\ad x)^{-1}(A) - A + B$, 
so $B=0$. This completes the proof of Lemma \ref{lissable}.\fin

\bi

\sub{Isolated class} An element $g \in \Gb$ is said to be {\it ($\Gb$-)isolated} 
if the centralizer of its semisimple part is not contained in 
a \levi of a proper \para of $\Gb$.

Let $\Lb_\isol$ denote the subset of $\Lb$ consisting of 
$\Lb$-isolated elements, and let $\Tb_\isol=\Tb \cap \Lb_\isol$. 
Then $\Tb_\isol$ is a closed subset of $\Tb$. Therefore 
$\na_\Lb(\Tb_\isol)$ is a closed subset of $\Tb/W_\Lb$. 
Moreover, we have $\Lb_\isol=\na_\Lb^{-1}(\na(\Tb_\isol))$, 
so $\Lb_\isol$ is a closed subset of $\Lb$. 

As a consequence of Lemma \ref{proj mor}, the image of the morphism 
$\Gb \times_\Pb \Lb_\isol.\Vb \to \Gb$, $g *_\Pb x \mapsto gxg^{-1}$ is 
a closed subvariety of $\Gb$~: we denote it by $\Xb_{\Gb,\Lb}$. 
On the other hand, by Lemma \ref{lissable}, the image of the morphism 
$\Gb \times_\Lb \OCB \to \Gb$, $g *_\Lb x \mapsto gxg^{-1}$ 
is an open subset of $\Gb$, which will be denoted by $\OCB_{\Gb,\Lb}$. 
Finally, we denote by $\Yb_{\Gb,\Lb}$ the intersection of $\Xb_{\Gb,\Lb}$ 
and $\OCB_{\Gb,\Lb}$.

Note that $W_\Gb(\Lb)=N_\Gb(\Lb)/\Lb$ acts (on the right) on 
the variety $\Gb \times_\Lb \Lb_{\isol,\reg}$ in the following way. 
If $w \in W_\Gb(\Lb)$ and if $g *_\Lb l \in \Gb \times_\Lb \Lb_{\isol,\reg}$, 
then 
$$(g *_\Lb l).w=g\wdo *_\Lb \wdo^{-1} l \wdo,$$
where $\wdo \in N_\Gb(\Lb)$ is any representative of $w$.

\bi

\begin{prop}\label{etalisable}
The map $\Gb \times_\Lb \Lb_{\isol,\reg} \to \Yb_{\Gb,\Lb}$, 
$g *_\Lb l \mapsto glg^{-1}$ 
is a Galois \'etale covering with group $W_\Gb(\Lb)$.
\end{prop}

\bi

\proof Set $\pi : \Gb \times_\Lb \Lb_{\isol,\reg} \to \Yb_{\Gb,\Lb}$, 
$g *_\Lb l \mapsto glg^{-1}$, and 
$\g : \Gb \times_\Lb \OCB \to \Gb$, $g *_\Lb l \mapsto glg^{-1}$. 
By Lemma \ref{lissable}, $\g$ is an \'etale morphism. If we prove 
that the square
$$\diagram 
\Gb \times_\Lb \Lb_{\isol,\reg} \rrto \ddto_{\DS{\pi}} && 
\Gb \times_\Lb \OCB \ddto^{\DS{\g}} \\
&&\\
\Yb_{\Gb,\Lb} \rrto && \Gb 
\enddiagram\leqno{(\#)}$$
is cartesian, then, by base change, we get that $\pi$ is an 
\'etale morphism. 

Since $\g$ is smooth, the fibred product (scheme) 
of $\Gb \times_\Lb \OCB$ and $\Yb_{\Gb,\Lb}$ over $\Gb$ 
is reduced (because $\Yb_{\Gb,\Lb}$ is), so it is enough 
to prove that $\Gb \times_\Lb \Lb_{\isol,\reg} =\g^{-1}(\Yb_{\Gb,\Lb})$. 

Let $g *_\Lb l \in \Gb \times_\Lb \OCB$ be such that 
$glg^{-1} \in \Yb_{\Gb,\Lb}$. Then there exist $h \in \Gb$, 
$m \in \Lb_\isol$ and $v \in \Vb$ such that $hmvh^{-1}=glg^{-1}$. 
Let $s$, $t$ and $t'$) denote the semisimple parts of $l$, $m$ 
and $mv$ respectively. By Lemma \ref{V conjugue}, there exists $x \in \Vb$ 
such that $t'=\lexp{x}{t}$. Thus $\lexp{hx}{t}=\lexp{g}{s}$. 

Since $t$ is $\Lb$-isolated, we have $\Zb(C_\Lb^\ci(t))^\ci=\Zb(\Lb)^\ci$. 
Therefore $\Zb(C_\Gb^\ci(t))^\ci \incl \Zb(\Lb)^\ci$. 
On the other hand, since $s \in \OCB$, we have $C_\Gb^\ci(s) \incl \Lb$, 
so $\Zb(\Lb)^\ci \incl \Zb(C_\Gb^\ci(s))$. This proves that 
$\lexp{g}{\Zb(\Lb)^\ci} \incl \lexp{hx}{\Zb(\Lb)^\ci}$. For dimension 
reasons, we have $\lexp{g}{\Zb(\Lb)^\ci}=\lexp{hx}{\Zb(\Lb)^\ci}$,  
so $\Zb(\Lb)^\ci=\Zb(C_\Gb^\ci(s))$. Hence, $l$ is isolated. 
So, we have proved that $\pi$ is \'etale. 

Now $W_\Gb(\Lb)$ acts freely on $\Gb \times_\Lb \Lb_{\isol,\reg}$. 
So the quotient morphism 
$\Gb \times_\Lb \Lb_{\isol,\reg} \to \Gb \times_{N_\Gb(\Lb)} \Lb_{\isol,\reg}$ 
is a Galois \'etale covering with group $W_\Gb(\Lb)$. 
Moreover, $\pi$ clearly factorizes through this quotient 
morphism. We get an \'etale morphism 
$\pi_0 : \Gb \times_{N_\Gb(\Lb)} \Lb_{\isol,\reg} \to \Yb_{\Gb,\Lb}$. 
To finish the proof of Proposition \ref{etalisable}, it is now enough to prove 
that $\pi_0$ is an isomorphism of varieties. Since $\pi_0$ is \'etale, 
we only need to prove that it is bijective. 

Note that $\pi_0$ is clearly surjective. We just need to prove that 
it is injective. Let $(g,l)$ and $(g',l')$ in $\Gb \times \Lb_{\isol,\reg}$ 
be such that $glg^{-1}=g'l'g^{\prime -1}$. Let $s$ and $s'$ be the semisimple 
parts of $l$ and $l'$ respectively. 
Then $\Zb(C_\Lb^\ci(s))^\ci=\Zb(\Lb)^\ci$ since $l$ is $\Lb$-isolated. 
Since $l \in \OCB$, we also have $C_\Gb^\ci(s)=C_\Lb^\ci(s)$, so 
$\Zb(C_\Gb^\ci(s))^\ci=\Zb(\Lb)^\ci$. Similarly, 
$\Zb(C_\Gb^\ci(s'))^\ci=\Zb(\Lb)^\ci$. So 
$\lexp{g^{-1}g'}{\Zb(\Lb)^\ci}=\Zb(\Lb)^\ci$. Therefore, $g^{-1}g' \in N_\Gb(\Lb)$. 
This completes the proof.\fin

\bi

\sub{Self-opposed Levi subgroups\label{distingue sub}} 
In this paper, we will be interested in the action of the relative Weyl 
group $W_\Gb(\Lb)$ on certain varieties. In this subsection (which is 
independent from the four previous ones), 
we deal with a particular case (where $\Lb$ is called self-opposed) 
for which this group is a reflection group for its action on $X(\Zb(\Lb)^\ci)$.  
The reader should note that most of the facts stated here will be used 
only in the next part. 

The \levi $\Lb$ of the \para $\Pb$ of $\Gb$ is said to be {\it $\Gb$-self-opposed} 
if, for every minimal \para $\Qb$ of $\Gb$ containing $\Pb$, we have $|W_\Mb(\Lb)|=2$, 
where $\Mb$ is the unique \levi of $\Qb$ containing $\Lb$. 
We recall in the next proposition some well-known basic properties of a 
$\Gb$-self-oposed \levi of a \para of $\Gb$. Most of them are due 
to Howlett \cite{howlett} and Lusztig \cite{lucox}.

\bi

\begin{prop}\label{propriete distingues}
Assume that $\Lb$ is $\Gb$-self-opposed. 
For every $\a \in \D-\D_\Lb$, let $\Qb_\a$ denote the \para 
of $\Gb$ generated by $\Pb$ and $\Ub_{-\a}$, $\Mb_\a$ the unique \levi of 
$\Qb_\a$ containing $\Lb$, and $s_{\Lb,\a}$ the unique non-trivial 
element of $W_{\Mb_\a}(\Lb)$. Then~:

\tete{a} For each $\a \in \D-\D_\Lb$, $s_{\Lb,\a}$ is a reflection on $X(\Zb(\Lb)^\ci)$.

\tete{b} $(W_\Gb(\Lb),(s_{\Lb,\a})_{\a \in \D-\D_\Lb})$ is a Coxeter system.

\tete{c} If $\Pb'$ is a \para of $\Gb$ having $\Lb$ as a Levi subgroup, then 
$\Pb'$ and $\Pb$ are conjugate in $\Gb$ (and in fact, they are conjugate under 
$N_\Gb(\Lb)$).

\tete{d} Let $\Mb$ be a \levi of a \para of $\Gb$ which contains $\Lb$ 
and let $g \in \Gb$ be such that $\lexp{g}{\Lb} \incl \Mb$. Then there 
exists $m \in \Mb$ such that $\lexp{g}{\Lb}=\lexp{m}{\Lb}$.
\end{prop}

\bi

\proof \cf \cite{howlett} and \cite{lucox} 
for (a), (b) and (c). We shall provide a general proof for (d).
For this, we may assume that $\Mb$ is standard 
with respect to $(\Tb,\Bb)$. Let $\Qb$ denote the \para of $\Gb$ 
which contains $\Pb$ and which has $\Mb$ as a Levi complement. 
By replacing $g$ by $mg$, where $m$ is a suitable element 
of $\Mb$, we may assume that $\lexp{g}{\Lb}$ is a standard 
Levi subgroup with respect to $(\Tb,\Bb \cap \Mb)$. 

Let $\Pb'$ be the standard \para of $\Mb$ having $\lexp{g}{\Lb}$ as a 
Levi complement. Then $\Pb$ and $\Pb' \Vb_\Qb$ are standard 
\paras having $\Lb$ and $\lexp{g}{\Lb}$ as respective Levi complements 
(here, $\Vb_\Qb$ denotes the unipotent radical of $\Qb$). 
By (c), $\Pb$ and $\Pb' \Vb_\Qb$ are conjugate. Since they both 
contain $\Bb$, they are equal. Hence $\lexp{g}{\Lb}=\Lb$.\fin

\bi

\remark{exemple} (1) It is well-known that $|W_\Mb(\Lb)| \le 2$ for every 
minimal \para $\Qb$ of $\Gb$ containing $\Pb$ and 
where $\Mb$ is the unique \levi of $\Qb$ containing $\Lb$.

\tete{2} It happens frequently that $W_\Gb(\Lb)$ is a reflection group 
for its action on $X(\Zb(\Lb)^\ci)$ and that $\Lb$ is not $\Gb$-self-opposed 
(see \cite{howlett} for the complete analysis of this question).

\tete{3} Compare Proposition \ref{propriete distingues} with 
\cite[Fact 1.1 (ii)]{DLM2}.\finl

\bi

A morphism $\s : \Gbh \to \Gb$ between two connected reductive groups 
is said to be {\it isotypic} if $\Ker \s$ is central in $\Gbh$ and 
$\Im \s$ contains the derived subgroup of $\Gb$. Note that in this case, 
the morphism $\s$ induces an isomorphism between the Dynkin diagram 
of $\Gbh$ and $\Gb$. 

The connected reductive group $\Gb$ is said to be 
{\it universally self-opposed} if, 
for every isotypic morphism $\Gbh \to \Gb$ such that $\Gbh$ 
is a \levi of a \para of a connected reductive group $\Gambh$, 
then $\Gbh$ is $\Gambh$-self-opposed. 

\bi

\example{T distingue} (1) A torus is universally self-opposed. 

\tete{2} If a unipotent class 
of $\Gb$ supports a cuspidal local system, then $\Gb$ is universally 
self-opposed \cite[Theorem 9.2]{luicc}.

\tete{3} The group $\Gb\Lb_2(\FM)$ is $\Sb\pb_4(\FM)$-self-opposed. 
However, it is not universally self-opposed. 
Indeed, $\Gb\Lb_2(\FM) \times \FM^\times$ is a \levi of a \para of $\Gb\Lb_3(\FM)$. 
But it is not $\Gb\Lb_3(\FM)$-self-opposed~: one can immediately check that 
$N_{\Gb\Lb_3(\FM)}(\Gb\Lb_2(\FM)\times\FM^\times)=\Gb\Lb_2(\FM)\times\FM^\times$.\finl

\bi

\rem In the statement (d) of Proposition \ref{propriete distingues}, if we assume 
only that $\Mb$ is a connected reductive subgroup of $\Gb$ (and not 
necessarily a \levi of a parabolic subgroup), then the conclusion does not 
hold. Indeed, 
if $\Gb=\Sb\pb_4(\FM)$,  
if $\Mb = \Sb\Lb_2(\FM) \times \Sb\Lb_2(\FM) \incl \Gb$, 
if $\Lb = \Sb\Lb_2(\FM) \times \FM^\times \incl \Mb$, 
then there exists $g \in \Gb$ such that 
$\lexp{g}{\Lb}= \FM^\times \times \Sb\Lb_2(\FM) \incl \Mb$, but $\Lb$ and 
$\lexp{g}{\Lb}$ are obviously not conjugate in $\Mb$. However, 
if $p\not= 2$, then the regular unipotent class of $\Lb$ supports 
a cuspidal local system, so $\Lb$ is universally self-opposed by 
Example \ref{T distingue} (2).\finl

\bi

\sec{Action of the relative Weyl group\label{sub relatif}}~ 

\med

\sub{The set-up} {\it From now on, and until the end of this paper, we 
denote by $\Sigb$ the inverse image of an $\Lb/\Zb(\Lb)^\ci$-isolated 
class of $\Lb/\Zb(\Lb)^\ci$ under the canonical projection 
$\Lb \to \Lb/\Zb(\Lb)^\ci$.} We also fix an element $v \in \Sigb$. 

Following \cite[\SEC\SEC 3 and 4]{luicc}, we consider the varieties
$$\Ybh=\Gb \times \Sigb_\reg,$$
$$\Ybt=\Gb \times_\Lb \Sigb_\reg,$$
$$\Xbh=\Gb \times \overline{\Sigb}\Vb$$
$$\Xbt=\Gb \times_\Pb \overline{\Sigb}\Vb.\leqno{\mathrm{and}}$$
In these definitions, the group $\Lb$ (\resp $\Pb$) acts on $\Gb$ by right 
translations, and acts on $\Sigb_\reg$ (\resp $\overline{\Sigb} \Vb$) 
by conjugation. We also set 
$$\Yb=\bigcup_{g \in \Gb} g \Sigb_\reg g^{-1}$$
$$\Xb=\bigcup_{g \in \Gb} g \overline{\Sigb} \Vb g^{-1}.\leqno{\mathrm{and}}$$
By Lemma \ref{V conjugue}, we have 
\equat\label{x}
\Xb \incl \na^{-1}\bigl(\na(\Zb(\Lb)^\ci)\bigr)
\endequat
and
\equat\label{y}
\Yb \incl \na^{-1}\bigl(\na(\Zb(\Lb)^\ci_\reg)\bigr).
\endequat
Moreover, $\Xb$ is the image of $\Gb \times_\Pb \overline{\Sigb}.\Vb$ under 
the morphism $\Gb \times_\Pb \Pb \to \Gb$, $g *_\Pb x \mapsto gxg^{-1}$ 
studied in Subsection \ref{sub fam}. So, by Lemma \ref{proj mor}, 
$\Xb$ is a closed irreducible subvariety of $\Gb$. We set 
$$\fonction{\piba}{\Gb \times_\Pb \overline{\Sigb}.\Vb}{\Xb}{g*_\Pb x}{gxg^{-1}.}$$
It is a projective morphism of varieties.

On the other hand, $\Yb^+=\bigcup_{g \in \Gb} g (\overline{\Sigb})_\reg g^{-1}$ 
is in fact the intersection of $\Xb$ with $\na^{-1}\bigl(\na(\Zb(\Lb)^\ci_\reg)\bigr)$, 
so it is an open subset of $\Xb$ (by Lemma \ref{nabla -2}). 
Moreover, the inverse image of $\Yb^+$ in $\Gb \times_\Lb \OCB$ is 
equal to $\bigcup_{w \in W_\Gb(\Lb)} \Gb \times_\Lb \lexp{w}{(\overline{\Sigb})_\reg}$. 
So, by 
Proposition \ref{etalisable}, 
the map $\Gb \times_\Lb (\overline{\Sigb})_\reg \to \Yb^+$, 
$g *_\Lb l \mapsto glg^{-1}$ is a Galois 
\'etale covering with group $W_\Gb(\Lb,\Sigb)=N_\Gb(\Lb,\Sigb)/\Lb$. 
Since $\Gb \times_\Lb \Sigb_\reg$ is open in $\Gb \times_\Lb (\overline{\Sigb})_\reg$, 
its image $\Yb$ under this \'etale morphism is open in $\Yb^+$. 
This proves that $\Yb$ is open in $\Xb$. Moreover, 
since the map 
$$\fonction{\pi}{\Gb \times_\Lb \Sigb_\reg}{\Yb}{g *_\Lb l}{glg^{-1}}$$
is a Galois \'etale covering with group $W_\Gb(\Lb,\Sigb)$, 
we get that $\Yb$ is smooth (indeed, $\Gb \times_\Lb \Sigb_\reg$ is smooth). 

Recall that $\Gb \times_\Lb \Sigb_\reg \to \Gb \times_\Pb \Sigb_\reg.\Vb$, 
$g *_\Lb l \mapsto g*_\Pb l$ is an isomorphism (see Corollary \ref{U conjugue bis}). 
Moreover, it is clear that $\piba^{-1}(\Yb)=\Gb \times_\Pb \Sigb_\reg.\Vb$. 
We summarize all these facts in the following proposition.

\bi

\begin{prop}[{\bf Lusztig \cite[\P 3.1, 3.2 and Lemma 4.3]{luicc}}]\label{resume truc}
With the above notation, we have~:

\tete{1} $\Xb$ is a closed irreducible subvariety of $\Gb$ and 
$\Yb$ is open in $\Xb$.

\tete{2} The natural map $\Ybt \to \Xbt$, $g *_\Lb x \mapsto g *_\Pb x$ is an 
open immersion and the square
$$\diagram 
\Ybt \rrto^{\DS{\pi}} \ddto && \Yb \ddto\\
&&\\
\Xbt 
\rrto^{\DS{\piba}}  && \Xb
\enddiagram$$
is cartesian.

\tete{3} $\piba$ (hence $\pi$) is a projective morphism.

\tete{4} $\pi$ is an \'etale Galois covering with group $W_\Gb(\Lb,\Sigb)$. 
In particular, $\Yb$ is smooth.
\end{prop}

\bi

Note that $\Gb$ acts on $\Ybh$, $\Ybt$, $\Xbh$ and $\Xbt$ 
by left translation on the first factor, and that it acts on $\Yb$ and $\Xb$ 
by conjugation. Also, the group $\Zb(\Gb) \cap \Zb(\Lb)^\ci$ acts 
on the varieties $\Ybh$, $\Ybt$, $\Xbh$ and $\Xbt$ 
by left translation on the second factor, and it acts on $\Yb$ and $\Xb$ 
by left translation. These actions of $\Gb$ and $\Zb(\Gb) \cap \Zb(\Lb)^\ci$ 
commute.

We have the following commutative diagram
\equat\label{diagramme central}
\diagram
\Sigb_\reg \ddto && \Ybh \llto_{\DS{\a}} 
\rrto^{\DS{\b}} \ddto && \Ybt \rrto^{\DS{\pi}} \ddto && \Yb \ddto\\
&&&&&&\\
\overline{\Sigb} && \Xbh \llto_{\DS{\alpba}} \rrto^{\DS{\betba}} && \Xbt 
\rrto^{\DS{\piba}}  && \Xb,
\enddiagram
\endequat
where $\a$ and $\alpba$ are the canonical projections, and $\b$ and $\betba$ 
are the canonical quotient morphisms. 
Moreover, the vertical maps are the natural ones. 
Note that the morphisms $\b$, $\betba$, $\pi$ and $\piba$ are $\Gb \times 
(\Zb(\Gb) \cap \Zb(\Lb)^\ci)$-equivariant, 
and that the vertical maps $\Ybh \to \Xbh$, $\Ybt \to \Xbt$ and $\Yb \to \Xb$ 
are also $\Gb \times (\Zb(\Gb) \cap \Zb(\Lb)^\ci)$-equivariant.

\bi

\sub{Extension of the action of $W_\Gb(\Lb,\Sigb)$} 
We observe that the action of the group $W_\Gb(\Lb,\Sigb)$ is defined only 
on $\Ybt$. 
However, we will see in this subsection that it is possible to 
extend it to an open subset $\Xbt_\mini$ of $\Xbt$ 
which, in general, strictly contains $\Ybt$.

We first need some preliminaries to construct this extension. 
If $l \in \overline{\Sigb}$, then 
$$\dim C_\Lb(l) \ge \dim C_\Lb(v)$$
and equality holds if 
and only if $l \in \Sigb$. Consequently, if $g \in \Xb$, then 
$$\dim C_\Gb(g) \ge \dim C_\Lb(v)$$ 
(\cf Lemma \ref{dim p} (2)). Moreover, if 
the equality holds, then, by Lemma \ref{dim p} and the previous 
remark, $\piba^{-1}(g)$ is contained in 
$$\Xbt_0=\Gb \times_\Pb \Sigb.\Vb,$$
which is a smooth open subset of $\Xbt$.

Let 
$$\Xb_\mini=\{g \in \Xb~|~ \dim C_\Gb(g) = \dim C_\Lb(v)\}$$
$$\Xbt_\mini=\piba^{-1}(\Xb_\mini).\leqno{\mathrm{and}}$$    
By \cite[Proposition 1.4]{humphreys}, $\Xb_\mini$ is an open subset of $\Xb$, and, 
by the previous discussion, 
$\Xbt_\mini \incl \Xbt_0$. Also, $\Yb \incl \Xb_\mini$, so $\Ybt \incl \Xbt_\mini$.
Now, let $\pi_\mini : \Xbt_\mini \to \Xb_\mini$ 
denote the restriction of $\piba$~: it is a projective morphism. 

Moreover, if $g \in \Xb_\mini$, then, by Lemma \ref{dim p} (1), 
$\pi_\mini^{-1}(g)$ is a finite set. 
Since the \mor $\pi_{\min}$ is projective and quasi-finite, it is finite 
\cite[Exercise III.11.2]{hartshorne}. We gather these facts in the 
next proposition.

\bi

\begin{prop}\label{tas}
With the above notation, we have~:

\tete{1} $\Xb_\mini$ is a $\Gb \times (\Zb(\Gb) \cap \Zb(\Lb)^\ci)$-stable 
open subset of $\Xb$ containing $\Yb$.

\tete{2} $\Xbt_\mini$ is a $\Gb \times (\Zb(\Gb) \cap \Zb(\Lb)^\ci)$-stable 
smooth open subset of $\Xbt$ containing $\Ybt$.

\tete{3} The \mor $\pi_\mini : \Xbt_\mini \to \Xb_\mini$ is finite.
\end{prop}

\bi

The Proposition \ref{tas} has the following immediate consequence~:

\bi

\begin{theo}\label{extension action}
$(a)$ The variety $\Xbt_\mini$ is the normalization of $\Xb_\mini$ in the variety 
$\Ybt$. Therefore, there exists a unique action of the finite group $W_\Gb(\Lb,\Sigb)$ 
on the variety $\Xbt_\mini$ extending its action on $\Ybt$. 

\tete{b} This action is $\Gb \times (\Zb(\Gb) \cap \Zb(\Lb)^\ci)$-equivariant, and 
the \mor $\pi_\mini$ factorizes through the quotient $\Xbt_\mini/W_\Gb(\Lb,\Sigb)$.

\tete{c} If the variety $\Xb_\mini$ is normal, then $\pi_\mini$ induces an \iso 
of varieties $\Xbt_\mini/W_\Gb(\Lb,\Sigb) \simeq \Xb_\mini$.
\end{theo}

\bi

\noindent{\sc Notation -} (1) Let $(\Sigb.\Vb)_\mini$ denote the open subset 
$\Sigb.\Vb \cap \Xb_\mini$ of $\Sigb.\Vb$. Then $\Xbt_\mini=\Gb \times_\Pb 
(\Sigb.\Vb)_\mini$. 

\med

\tete{2} If there is some ambiguity, we will denote by $?_\Lb^\Gb$ the 
object $?$ defined above (for instance, $\Ybh_\Lb^\Gb$, $\Xbt_{\min,\Lb}^\Gb$, 
$\Xb_\Lb^\Gb$, $\pi_{\min,\Lb}^\Gb$...).

\bi

\sub{Unipotent classes} {\it From now on, and until the end of this paper, 
$\Sigb$ is the inverse image of a unipotent class of $\Lb/\Zb(\Lb)^\ci$.} 
Note that a unipotent class is isolated. Let $\Cb$ denote 
the unique unipotent class contained in $\Sigb$. From now on, the 
element $v$ introduced in the previous section is chosen to be 
in $\Cb$. Note that $\Sigb=\Zb(\Lb)^\ci.\Cb \simeq \Zb(\Lb)^\ci \times \Cb$ 
and that $\Sigb_\reg=\Zb(\Lb)^\ci_\reg.\Cb$.

\bi

\noindent{\sc Notation - } 
(1) We denote by $\Cb^\Gb$ the induced unipotent class of $\Cb$ 
from $\Lb$ to $\Gb$ that is, the unique unipotent class $\Cb_0$ of $\Gb$ 
such that $\Cb_0 \cap \Cb.\Vb$ is dense in $\Cb.\Vb$.

\medskip

\tete{2} If $z \in \Zb(\Lb)^\ci$, the group $C_\Pb(z)=\Lb.C_\Vb(z)$ is connected 
(\cf Lemma \ref{spa connexe}). So it is a \para of $C_\Gb^\ci(z)$ 
by \cite[Proposition 1.11 (ii)]{dmgnc}, with unipotent radical $C_\Vb(z)$ and 
Levi factor 
$\Lb$. We denote by $u_z$ an element of 
$\Cb^{C_\Gb^\ci(z)} \cap vC_\Vb(z)$. We set $\uti_z=1*_\Pb zu_z 
\in \Xbt_\mini$. 

\med

\tete{3} For simplicity, the unipotent element $u_1$ will be denoted by $u$, 
and $\uti$ stands for $\uti_1$. 

\bi

\remark{description xmin} Let us investigate here what are the elements 
of $\Xb_{\min}$. Since $\Xbt_\mini \incl \Xbt_0$, we only need to determine 
which elements of $\Sigb.\Vb$ belong to $\Xb_{\min}$. 
Let $g \in \Sigb.\Vb$. Let $z$ (\resp $u'$) be the semisimple 
(\resp unipotent) part of $g$. By Lemma \ref{V conjugue}, we may assume 
that $z$ belongs to $\Zb(\Lb)^\ci$. Now, let $\Gb'=C_\Gb^\ci(z)$, 
$\Pb'=C_{\Pb}^\ci(z)$, and $\Vb'=C_\Vb(z)$. Then $\Gb'$ is a reductive subgroup 
of $\Gb$ containing $\Lb$, $\Pb'$ is a \para of $\Gb'$, and $\Vb'$ is its 
unipotent radical. Then, by 
\cite[Corollary 11.12]{borel}, we have $u' \in \Gb'$. On the other hand, 
$C_\Gb^\ci(g) = C_{\Gb'}^\ci(u')$. Now, by Lemma \ref{dim p} (2) and by 
\cite[Proposition II.3.2 (b) and (e)]{spaltenstein}, $g \in \Xb_{\min}$ if and only 
if $u' \in \Cb^{\Gb'}$.

Hence, we have proved that
\equat\label{une seule}
\na^{-1}(\na(z)) \cap \Xb_\mini =(zu_z)_\Gb
\endequat
for every $z \in \Zb(\Lb)^\ci$, and that
$$\Xb_\mini=\bigcup_{z \in \Zb(\Lb)^\ci} (zu_z)_\Gb.~\SS{\square}$$

\bi

If $z \in \Zb(\Lb)^\ci$, 
we denote by $H_\Gb(\Lb,\Sigb,z)$ the stabilizer 
of $\uti_z$ in $W_\Gb(\Lb,\Sigb)$. 
We first investigate the group $H_\Gb(\Lb,\Sigb,1)$. 
Recall that the group $C_\Gb^\ci(u)$ is 
contained in $\Pb$ \cite[Proposition II.3.2 (e)]{spaltenstein}, so that 
$C_\Gb(u)/C_\Pb(u)$ is a finite set.

\bi

\begin{lem}\label{agu general}
We have $\pi_\mini^{-1}(u)=\{g *_\Pb u~|~g \in C_\Gb(u)\}$. In particular, 
$$|\pi_\mini^{-1}(u)|=|C_\Gb(u)/C_\Pb(u)|=|A_\Gb(u)/A_\Pb(u)|.$$
\end{lem}

\bi

\proof It is clear that $\{g *_\Pb u~|~g \in C_\Gb(u)\}$ is contained in 
$\pi_\mini^{-1}(u)$. Conversely, let $g *_\Pb x \in \pi_\mini^{-1}(u)$. 
By replacing $(g,x)$ by $(gl^{-1},lx)$ for a suitable choice of 
$l \in \Lb$, we may assume that $\pi_\Lb(x)=v$. Since $gxg^{-1} =u$, 
this implies that $x \in v\Vb \cap \Cb^\Gb$. By 
\cite[Proposition II.3.2 (d)]{spaltenstein}, there exists $y \in \Pb$ such 
that $yxy^{-1}=u$. Therefore $gy^{-1} \in C_\Gb(u)$ and 
$g*_\Pb x= gy^{-1} *_\Pb u$.\fin

\bi

\begin{coro}\label{agu 1}
If $C_\Gb(u) \incl \Pb$, then $\pi_\mini^{-1}(u)=\{\uti\}$. In particular, 
$W_\Gb(\Lb,\Sigb)$ stabilizes $\uti$, that is $H_\Gb(\Lb,\Sigb,1)=W_\Gb(\Lb,\Sigb)$. 
\end{coro}

\bi

Next, let us consider the general case. The second projection $\Xbh
\simeq \Gb \times \Zb(\Lb)^\ci \times \overline{\Cb} \times \Vb \to \Zb(\Lb)^\ci$ 
factors through the quotient morphism $\Xbh \to \Xbt$. We 
denote by $\nablat : \Xbt \to \Zb(\Lb)^\ci$ the morphism obtained after 
factorization. The group $W_\Gb(\Lb)$ acts on $\Zb(\Lb)^\ci$ by conjugation, 
and it is easy to check that the restriction $\nablat_\reg : \Ybt \to \Zb(\Lb)^\ci_\reg$ 
of $\nablat$ to $\Ybt$ is $W_\Gb(\Lb,\Sigb)$-equivariant. Hence, 
the morphism $\nablat_\mini : \Xbt_\mini \to \Zb(\Lb)^\ci$ obtained by 
the restriction of $\nablat$ is $W_\Gb(\Lb,\Sigb)$-equivariant.

As a consequence, we get 
\equat\label{ZZZ}
\Stab_{W_\Gb(\Lb,\Sigb)}(\gti) \incl \Stab_{W_\Gb(\Lb,\Sigb)}(\nablat(\gti))
\endequat
for every $\gti \in \Xbt_\mini$. Also, note that $\nablat(\gti)$ is conjugate 
in $\Gb$ to the semisimple part of $g=\pi_\mini(\gti)$ (\cf Proposition 
\ref{V conjugue}).
In fact, one can easily get a better result~:

\bi

\begin{prop}\label{stabilisateur z}
Let $z \in \Zb(\Lb)^\ci$. Then~:

\tete{1} $H_\Gb(\Lb,\Sigb,z)=H_{C_\Gb^\ci(z)}(\Lb,\Sigb,1)$.

\tete{2} If $C_{C_\Gb^\ci(z)}(u_z) \incl \Pb$, then 
$H_\Gb(\Lb,\Sigb,z) = W_{C_\Gb^\ci(z)}(\Lb,\Sigb)$.
\end{prop}

\bi

\noindent{\sc Sketch of the proof -} 
For the proof of (1), we refer the reader to the proof of 
Proposition \ref{z} below. The situation here is quite analogous to that 
of Proposition \ref{z} and the arguments involved are very similar. 
We have decided to present them in detail only once, for Proposition 
\ref{z}, since the situation there is a little bit more complicated.
(2) follows from (1) and from Corollary \ref{agu 1}.\fin

\bi

\sub{An example} 
In this subsection, and only in this subsection, we assume that $\Lb=\Tb$. 
Then $\Cb=1$, $\Sigb=\Tb$, $\Xbt = \Gb \times_\Bb \Bb$, 
$\Xb=\Gb$ and $\piba : \Xbt \to \Gb$ is the well-known Grothendieck map. Also, 
$W_\Gb(\Lb,\Sigb)=W$ in this case. 
Moreover, $\Xb_\mini$ is the open subset of $\Gb$ consisting of regular elements. 
As an open subset of $\Gb$, it is smooth. So the action 
of $W$ on $\Ybt$ extends to $\Xbt_\mini$ and $\Xbt_\mini/W = \Xb_\mini$.

Now, let $\gti \in \Xbt_\mini$, $g=\pi_\mini(\gti)$, and $t=\nablat(\gti) \in \Tb$. 
We denote by $W^\ci(t)$ the Weyl group of $C_\Gb^\ci(t)$ relative to $\Tb$. The fibre 
$\pi_\mini^{-1}(g)$ may be identified with the set of Borel subgroups 
of $\Gb$ containing $g$. Since $\Xbt_\mini/W = \Xb_\mini$, $W$ acts transitively 
on $\pi_\mini^{-1}(g)$. But, $u_t$ is a regular unipotent element of $C_\Gb^\ci(t)$. 
Therefore, $C_{C_\Gb^\ci(t)}(u_t) \incl \Bb$. So, by Proposition \ref{stabilisateur z}, 
we have
$$\Stab_W(\gti)=W^\ci(t).$$
As a consequence, we get the well-known result 
\equat\label{nombre borel}
|\{x\Bb \in \Gb/\Bb~|~g \in \lexp{x}{\Bb}\}|=|W|/|W^\ci(t)|.
\endequat

\bi

\rem It is not the easiest way to prove \ref{nombre borel}~!

\bi

{\small
\example{sl2 2} Assume in this example that $\Gb=\Gb\Lb_2(\FM)$, that 
$$\Lb=\Tb=\{\diag(a,b)~|~a,b \in \FM^\times\},$$ 
and that 
$\Sigb=\Tb$. Let $\Pb^1$ denote the projective line. Then
$$\Xbt\simeq\{(\matrice{a & b \\ c & d},[x,y]) \in \Gb \times \Pb^1~|~ 
[ax+by,cx+dy]=[x,y]\},$$
$$\Xb=\Gb,$$ 
and $\pi : \Xbt \to \Xb$ is identified with the first projection. 
Moreover, $\Xb_\mini$ is the open subset of $\Gb$ consisting of 
non-central elements. We shall give a precise formula for describing 
the action of $W$ on $\Xbt_\mini$ in this little example.

Let $w$ denote the unique non-trivial element of $W$. 
It has order $2$. We define the right action of $w$ on 
$(g,[x,y]) \in \Xbt_\mini$ by
$$(g,[x,y]).w=\left\{\begin{array}{ll}
(g,[bx,(d-a)x-by]) & \quad{\mathrm{if}}~(bx,(d-a)x-by) 
\not=(0,0) \\
&\\
(g,[(a-d)y-cx,cy]) & \quad{\mathrm{if}}~((a-d)y-cx,cy)\not=(0,0),
\end{array}\right.$$
where $g=\matrice{a & b \\ c & d}$. 
One can check that, if $\Xbt_1$ (\resp $\Xbt_2$) is the open subset 
of $\Xbt$ defined by the first condition (\resp the second condition), 
then $\Xbt_1 \cup \Xbt_2 = \Xbt_\mini$, and that the formulas given above 
coincide on $\Xbt_1 \cap \Xbt_2$. So we have defined a \mor of varieties. 
The fact that it is an \auto of order $2$ is obvious, and the reader can 
check that it extends the action of $W$ on $\Ybt$. 

One can also check, as it is expected from \ref{nombre borel}, 
that $W$ acts trivially on the elements 
$(g,[x,y]) \in \Xbt_\mini$ such that $g$ is not semisimple.\finl}

\bi

\sec{A morphism $W_\Gb(\Lb,\Sigb) \to A_\Lb(v)$\label{mor phi}}~ 

\med

The restriction of an $\Lb$-equivariant local system on $\Cb$ 
through the morphism $\Lb/C_\Lb^\ci(v) \to \Cb$, 
$lC_\Lb^\ci(v) \mapsto lvl^{-1}$ is constant. This fact makes this 
morphism interesting when one is working with character 
sheaves (which are equivariant intersection cohomology complexes). 
This morphism 
can be followed through Diagram \ref{diagramme central}, 
and it gives rise to new varieties on which the group 
$W_\Gb(\Lb,v)=N_\Gb(\Lb,v)/C_\Lb^\ci(v)$ acts 
(note that $W_\Gb(\Lb,v)/A_\Lb(v) \simeq W_\Gb(\Lb,\Sigb)$). 
Following the method of the previous section, these actions 
can be extended to some variety $\Xbt_\mini^\prime$ lying over $\Xbt_\mini$. 
We show in this section that some stabilizers under this action
can be described in terms of a morphism of groups 
$W_\Gb(\Lb,\Sigb) \to A_\Lb(v)$ 
(under certain conditions which are fulfilled if $\Cb$ supports 
a cuspidal local system). In \SEC\ref{ele sec}, elementary 
properties of this morphism will be investigated.

\bi

\sub{Notation} Let $\Sigb'=\Lb/C_\Lb^\ci(v) \times \Zb(\Lb)^\ci$, and 
let $\Sigb_\reg^\pr = \Lb/C_\Lb^\ci(v) \times \Zb(\Lb)_\reg^\ci $. We denote by $f : 
\Sigb' \to \Sigb$, $(lC_\Lb^\ci(v),z) \mapsto lzvl^{-1}=zlvl^{-1}$. Then $f$ is a 
finite surjective $\Lb$-equivariant morphism (here, $\Lb$ acts 
on $\Sigb'$ by left translation on the first factor). We denote by $f_\reg : 
\Sigb_\reg^\pr \to \Sigb_\reg$ the restriction of $f$. 

Now, let 
$$\Ybh'= \Gb \times \Sigb_\reg^\pr,$$
$$\Ybt'=\Gb \times_\Lb \Sigb_\reg^\pr = \Gb/C_\Lb^\ci(v) \times \Zb(\Lb)^\ci_\reg.
\leqno{\mathrm{and}}$$
We then get a commutative diagram
\equat\label{diagramme}
\diagram
\Sigb_\reg^\pr \ddto_{\DS{f_\reg}} && \Ybh' \llto_{\DS{\a'}} \rrto^{\DS{\b'}} 
\ddto_{\DS{\fha}} &&
 \Ybt' \ddto_{\DS{\fti}} \ddrrto^{\DS{\pi'}} & \\
&&&&&&\\
\Sigb_\reg && \Ybh \llto^{\DS{\a}} \rrto_{\DS{\b}} && \Ybt \rrto_{\DS{\pi}} && \Yb,
\enddiagram
\endequat
where the vertical maps are induced by $f_\reg$, $\a'$ is the projection 
on the first factor, $\b'$ is the quotient morphism, and $\pi'=\pi \ci \fti$. 
The group $\Gb$ acts on $\Ybh'$ and $\Ybt'$ by left translation on the 
first factor, and acts on $\Yb$ by conjugation. The group $\Zb(\Gb) \cap \Zb(\Lb)^\ci$ 
acts on $\Sigb_\reg^\pr$ by translation on the second factor~: it induces 
an action on $\Ybh^\pr$ and $\Ybt'$. The morphisms $\fha$, $\fti$, 
$\b'$ and $\pi'$ are $\Gb \times (\Zb(\Gb) \cap \Zb(\Lb)^\ci)$-equivariant. 
Moreover, all the squares of diagram \ref{diagramme} are cartesian.

Now, we define  
$$W_\Gb(\Lb,v) = N_\Gb(\Lb,v)/C_\Lb^\ci(v)$$ 
(note that $N_\Gb(\Lb,v)^\ci=C_\Lb^\ci(v)$).
The group $C_\Lb(v)$ is a normal subgroup of $N_\Gb(\Lb,v)$ so $A_\Lb(v)$ is a normal 
subgroup of $W_\Gb(\Lb,v)$. Note that
$$W_\Gb(\Lb,v)/A_\Lb(v) \simeq W_\Gb(\Lb,\Sigb).$$
The group $N_\Gb(\Lb,v)$ acts freely on the right 
on the variety $\Ybh'$ in the following way~: if $w \in N_\Gb(\Lb,v)$ 
and if $(g,lC_\Lb^\ci(v),z) \in \Ybh'$, then
$$(g,lC_\Lb^\ci(v),z).w=(gw,w^{-1}lwC_\Lb^\ci(v),w^{-1}zw).$$
This induces a free right $\Gb \times (\Zb(\Gb) \cap \Zb(\Lb)^\ci)$-equivariant 
action of $W_\Gb(\Lb,v)$ on $\Ybt'$. 
Moreover, the fibres of the \mor $\pi'$ are $W_\Gb(\Lb,v)$-orbits.

\bi

\remark{action A} If $a \in A_\Lb(v)$ and $g *_\Lb (lC_\Lb^\ci(v),z) \in \Ybt'$, then 
$$(g *_\Lb (lC_\Lb^\ci(v),z)).a=g*_\Lb (laC_\Lb^\ci(c),z).~\SS{\square}$$

\bi

\sub{Normalization\label{sub nor}} 
Let $\Xbt'$ be the normalization of the variety $\Xbt$ in $\Ybt'$. We denote by 
$\fba : \Xbt' \to \Xbt$ the corresponding \mor of varieties. Let 
$\Xbt_0^\pr$ (\resp $\Xbt_\mini^\pr$) denote the inverse image, in $\Xbt'$, 
of the variety $\Xbt_0$ (\resp $\Xbt_\mini$). We denote by $\fti_0 : 
\Xbt_0^\pr \to \Xbt_0$ (\resp $\fti_\mini : \Xbt_\mini^\pr \to \Xbt_\mini$) 
the restriction of $\fba$ to $\Xbt_0^\pr$ (\resp $\Xbt_\mini^\pr$). 
Then $\Xbt_0^\pr$ (\resp $\Xbt_\mini^\pr$) is the normalization of $\Xbt_0$ 
(\resp $\Xbt_\mini$) in $\Ybt^\pr$. 
We summarize the notation in the following commutative diagram
$$\diagram
\Ybt' \rrto \ddto_{\DS{\fti}}\xto'[2,1]^{\DS{\pi'}}[4,2]&& 
\Xbt_\mini^\pr \ddto_{\DS{\fti_\mini}}\xto'[2,1]^{\DS{\pi_\mini^\pr}}[4,2]
\rrto && \Xbt_0^\pr \ddto^{\DS{\fti_0}} \rrto && \Xbt' \ddto^{\DS{\fba}} 
\xto[4,2]^{\DS{\piba}'}\\
&&\\
\Ybt \xto[0,2] \ddrrto_{\DS{\pi}} && \Xbt_\mini \ddrrto_{\DS{\pi_\mini}} \rrto 
&& \Xbt_0  \rrto && \Xbt \ddrrto^{\DS{\piba}}\\
&&&&\\
&& \Yb \rrto && \Xb_\mini \xto[0,4] &&&&\Xb. \\
\enddiagram$$
In this diagram, all the horizontal maps are open immersions, 
and all the squares are cartesian. 
Since $\Xbt_\mini$ is the normalization of $\Xb_\mini$ in $\Ybt$, we get~:

\begin{theo}\label{tango}
We have~:

\tete{1} The variety $\Xbt_\mini^\pr$ is the normalization of $\Xb_\mini$ in $\Ybt^\pr$. 
Therefore, the action of $W_\Gb(\Lb,v)$ on $\Ybt'$ extends uniquely to 
an action of $W_\Gb(\Lb,v)$ on $\Xbt_\mini^\pr$.

\tete{2} $\Xbt'$ inherits from $\Ybt'$ an action of $\Gb \times (\Zb(\Gb) 
\cap \Zb(\Lb)^\ci)$, and this action commutes with the action of $W_\Gb(\Lb,v)$ 
on $\Xbt_\mini$.
\end{theo}

\bi

\remark{xtilde 0} We do not know how to determine the variety 
$\Xbt'$ in general. However, it is possible to give an explicit description of 
$\Xbt_0^\pr$. 
This can be done as follows. The \para $\Pb$ acts on $\Sigb' \times \Vb$ 
by the following action~: 
if $l,l_0 \in \Lb$, $x, x_0 \in \Vb$, and $z_0 \in \Zb(\Lb)^\ci$, then 
$$\lexp{lx}{(l_0 C_\Lb^\ci(v),z_0,x_0)}=(ll_0C_\Lb^\ci(v),z_0,
l(\lexp{l_0z_0^{-1}v^{-1}l_0^{-1}}{x})x_0 x^{-1}l^{-1}).$$
It is easy to check that this defines an action of $\Pb$. Moreover, the 
\mor $f \times \Id_\Vb : 
\Sigb' \times \Vb \to \Sigb\Vb$, $(l_0C_\Lb^\ci(v),z_0,x_0) \mapsto 
l_0z_0vl_0^{-1}x_0$ induced by $f$ is $\Pb$-equivariant. 

By Corollary \ref{U conjugue bis}, 
$\Gb \times_\Lb \Sigb_\reg^\pr \simeq \Gb \times_\Pb (\Sigb_\reg^\pr \times \Vb)$ 
is an open subset of $\Gb \times_\Pb (\Sigb' \times \Vb)$ isomorphic 
to $\Ybt'$. Moreover, the morphism 
$$\Gb \times_\Pb (\Sigb' \times \Vb) \longto \Xbt=\Gb \times_\Pb \Sigb\Vb$$
induced by $f$ is finite, as it can be checked by restriction to an open subset 
of the form $g\Vb^- \Pb \times_\Pb (\Sigb' \times \Vb)$, where $\Vb^-$ 
is the unipotent radical of the opposite \para $\Pb^-$ of $\Pb$ 
with respect to $\Lb$. Finally, by the same argument, 
$\Gb \times_\Pb (\Sigb' \times \Vb)$ is smooth. Hence
\equat\label{normalization x0}
\Xbt_0^\pr=\Gb \times_\Pb (\Sigb' \times \Vb).
\endequat
Since $\Xbt'$ is the normalization of $\Xbt$ in $\Ybt'$, it inherits an action 
of the group $A_\Lb(v)$. It is very easy to describe this action on $\Xbt_0^\pr$ 
using \ref{normalization x0}.
It acts on $\Xbt_0^\pr$ by right translation on the 
factor $\Lb/C_\Lb^\ci(v)$ of $\Sigb'$. This is a free action and 
the fibres of $\fti_0$ are $A_\Lb(v)$-orbits.\finl

\bi

We will denote by $(\Sigb' \times \Vb)_\mini$ the inverse image, under 
$f \times \Id_\Vb$, of the open subset $(\Sigb\Vb)_\mini$ of $\Sigb\Vb$. 
Then $\Xbt_\mini^\pr=\Gb \times_\Pb (\Sigb' \times \Vb)_\mini$. The action 
of $W_\Gb(\Lb,v)$ is quite mysterious, but the action of its subgroup $A_\Lb(v)$ 
is understandable. It is obtained by restriction from its action on $\Xbt_0^\pr$ 
which is described at the end of Remark \ref{xtilde 0}.

\bi

\sub{Stabilizers\label{02}} For $z \in \Zb(\Lb)^\ci$, let $\uti_z^\pr=
1 *_\Pb (C_\Lb^\ci(v),z,v^{-1}u_z) \in \Xbt_\mini^\pr$. Recall 
that $u_z$ is an element of $v C_\Vb(z) \cap \Cb^{C_\Gb^\ci(z)}$. Note that 
$\fti_\mini^\pr(\uti^\pr_z)=\uti_z$, so that 
$\pi_\mini^\pr(\uti^\pr_z)=u_z$. For simplification, we denote the element 
$\uti_1^\pr$ by $\uti'$. The stabilizer of the element $\uti_z^\pr$ in 
$W_\Gb(\Lb,v)$ is denoted 
by $H_\Gb(\Lb,v,z)$. The goal of this subsection is to get some information 
about these stabilizers.

The first result comes from the fact that $A_\Lb(v)$ acts freely on $\Xbt_\mini^\pr$~:
\equat\label{H cap A}
H_\Gb(\Lb,v,z) \cap A_\Lb(v) =\{1\}.
\endequat
The second one is analogous to Proposition \ref{stabilisateur z}~: it may be viewed 
as a kind of Jordan decomposition.

\bi

\begin{prop}\label{z}
If $z \in \Zb(\Lb)^\ci$, then $H_\Gb(\Lb,v,z)=H_{C_\Gb^\ci(z)}(\Lb,v,1)$.
\end{prop}

\bi

\proof Let $\nablat'$ denote the composite morphism of varieties 
$\Xbt' \longmapright{\fba} \Xbt \longmapright{\nablat} \Zb(\Lb)^\ci$, and let 
$\nablat_\mini^\pr : \Xbt_\mini^\pr \to \Zb(\Lb)^\ci$ denote the restriction of 
$\nablat'$. Then $\nablat_\mini^\pr$ is a $W_\Gb(\Lb,v)$-equivariant morphism 
(as can be 
verfied by restriction to $\Ybt'$). So, the group $H_\Gb(\Lb,v,z)$ is contained 
in $\WC_z=W_{C_\Gb(z)}(\Lb,v)$. 

Let $\Ab_z=\{t \in \Zb(\Lb)^\ci~|~C_\Gb^\ci(t) \incl C_\Gb^\ci(z)\}$. $\Ab_z$ 
is an open subset of $\Zb(\Lb)^\ci$ containing $z$ and $\Zb(\Lb)^\ci_\reg$. 
Now let $\Sigb_z=\Ab_z .\Cb$ and let $\Sigb_z^\pr=\Lb/C_\Lb^\ci(v) \times \Ab_z$. 
Then 
$$\Xbt_z^\pr=\Gb\times_\Pb (\Sigb_z^\pr \times \Vb)_{\min,\Lb}^\Gb$$ 
is an open subset of $\Xbt_\mini$ containing $\uti_z$, and it is stable under 
the action of $\WC_z$, since $\Ab_z$ is and since 
$\Xbt_z^\pr=\nablat_\mini^{\pr -1}(\Ab_z)$. Now, let 
$$\Xbt'(z)=
C_\Gb(z) \times_{C_\Pb(z)} (\Sigb_z^\pr \times C_\Vb(z))_{\mini,\Lb}^{C_\Gb^\ci(z)}.$$
The natural morphism $\Xbt'(z) \to \Xbt_z^\pr$ is injective and $\WC_z$-equivariant. 
This proves that the stabilizer $H_\Gb(\Lb,v,z)$ is equal to the 
stabilizer of $1 *_{C_\Pb(z)} (C_\Lb^\ci(v),z,v^{-1}u_z) \in \Xbt'(z)$ 
in $\WC_z$. But this stabilizer must stabilize the connected component 
of $1 *_{C_\Pb(z)} (C_\Lb^\ci(v),z,v^{-1}u_z)$, which is equal to 
$C_\Gb^\ci(z) \times_{C_\Pb(z)} 
(\Sigb_z^\pr \times \Vb)_{\min,\Lb}^{C_\Gb^\ci(z)}$ (because $C_\Pb(z)$ 
is connected). Hence, it is contained 
in $\WC_z^\ci=W_{C_\Gb^\ci(z)}(\Lb,v)$, thus it is equal 
to $H_{C_\Gb^\ci(z)}(\Lb,v,z)$ because this last variety is an open 
subset of $(\Xbt_{\min}^\pr)_\Lb^{C_\Gb^\ci(z)}$. 

Now, the action of $\WC_z^\ci$ on $(\Xbt_{\min}^\pr)_\Lb^{C_\Gb^\ci(z)}$ 
commutes with the translation by $z$. Hence $H_{C_\Gb^\ci(z)}(\Lb,v,z)=
H_{C_\Gb^\ci(z)}(\Lb,v,1)$.\fin

\bi

\rem The reader may be surprised by the fact that $C_\Gb(z)$ is not necessarily 
connected. However, 
they can check directly that the previous constructions ($\Ybh$, $\Xbt_\mini$, 
$W_\Gb(\Lb,v)$...) remain valid whenever $\Gb$ is not connected, provided 
that the \para $\Pb$ of $\Gb$ is connected.\finl

\bi

In order to determine the stabilizers $H_\Gb(\Lb,v,z)$, 
Proposition \ref{z} shows that it is necessary and sufficient 
to compute the stabilizer $H_\Gb(\Lb,v,1)$. 
However, we are only able to get a satisfying result 
when the centralizer of $u$ in $\Gb$ is contained in the \para $\Pb$.

\bi

\begin{prop}\label{agu}
If $C_\Gb(u) \incl \Pb$, then~:

\tete{1} $\pi_\mini^{\pr -1}(u)$ is the $A_\Lb(v)$-orbit of $\uti'$. In particular, 
$|\pi_\mini^{\pr -1}(u)|=|A_\Lb(v)|$.

\tete{2} $W_\Gb(\Lb,v)=A_\Lb(v) \rtimes H_\Gb(\Lb,v,1)$.
\end{prop}

\bi

\proof (1) follows immediately from Corollary \ref{agu 1}~: indeed, 
$\pi_\mini^{\pr -1}(u)=\fti^{-1}(\uti)$. 
By (1), $A_\Lb(v)$ acts freely and transitively on 
$\pi_\mini^{\pr -1}(u)$, so $W_\Gb(\Lb,v)$ acts transitively on 
$\pi_\mini^{\pr -1}(u)$. (2) follows from this remark and from 
\ref{H cap A}.\fin

\bi

\sub{Further investigations\label{01}} 
The group $C_\Gb^\ci(v) \cap \Lb=C_{C_\Gb^\ci(v)}(\Zb(\Lb)^\ci)$ 
is connected, because it is the centralizer of a torus in a connected 
group \cite[Corollary 11.12]{borel}. Therefore, we have the well-known equality 
\equat\label{cgo}
C_\Gb^\ci(v) \cap \Lb =C_\Lb^\ci(v).
\endequat
Thus the natural morphism $C_\Lb(v) \injto C_\Gb(v)$ induces an injective 
\mor 
\equat\label{310}
A_\Lb(v) \injto A_\Gb(v).
\endequat

\bi

\example{contre alv} Let $\Gb \simeq \Sb\pb_4(\FM)$, $\Lb \simeq \Gb\Lb_2(\FM)$ 
and $v$ be a regular unipotent element of $\Lb$. Then $A_\Lb(v)=\{1\}$ 
and $|A_\Gb(v)|=2$. This shows that the morphism \ref{310} is in general 
not surjective.\finl

\bi

Let $W_\Gb^\ci(\Lb,v)=N_\Gb(\Lb,v) \cap C_\Gb^\ci(v)/C_\Lb^\ci(v)$. 
Since $C_\Gb^\ci(v) \cap \Lb=C_\Lb^\ci(v)$, we have $W_\Gb^\ci(\Lb,v) \cap 
A_\Lb(v)=1$. Moreover, since $W_\Gb^\ci(\Lb,v)$ and $A_\Lb(v)$ are normal subgroups 
of $W_\Gb(\Lb,v)$, $W_\Gb^\ci(\Lb,v) \times A_\Lb(v)$ is naturally 
a subgroup of $W_\Gb(\Lb,v)$. This discussion has the following immediate 
consequence~:

\bi

\begin{lem}\label{alv agv alv agv}
If $A_\Lb(v)=A_\Gb(v)$, then $W_\Gb(\Lb,v)=W_\Gb^\ci(\Lb,v) \times A_\Lb(v)$, 
so $W_\Gb^\circ(\Lb,v) \simeq W_\Gb(\Lb,\Sigb)$.
\end{lem}

\bi

\begin{coro}\label{mor fti}
Assume that $C_\Gb(u) \incl \Pb$, and that $A_\Lb(v)=A_\Gb(v)$. 
Then there exists a unique \mor of groups 
$\ph_{\Lb,v}^\Gb : W_\Gb^\ci(\Lb,v) \to A_\Lb(v)$ such that 
$$H_\Gb(\Lb,v,1)=\{(w,a) \in W_\Gb^\ci(\Lb,v) \times 
A_\Lb(v)~|~a=\ph_{\Lb,v}^{\Gb}(w)\}.$$
\end{coro}

\bi

\proof This follows from Proposition \ref{agu} (2) and from Lemma 
\ref{alv agv alv agv}.\fin

\bi

Recall that the unipotent element $v$ of $\Lb$ is said to be {\it distinguished 
in $\Lb$} if $v$ is not contained in a \levi of a proper \para of $\Lb$.

\bi

\begin{coro}\label{mor odd}
Assume that $v$ is distinguished in $\Lb$, that 
$C_\Gb(u) \incl \Pb$, that $A_\Lb(v)=A_\Gb(v)$, and that $|A_\Lb(v)|$ is odd. 
Then the morphism $\ph_{\Lb,v}^\Gb$ 
is trivial and $H_\Gb(\Lb,v,1)=W_\Gb^\ci(\Lb,v)$.
\end{coro}

\bi

\proof If $v$ is distinguished in $\Lb$, then $\Zb(\Lb)^\ci$ is a maximal 
torus of $C_\Gb^\ci(v)$. So $W_\Gb^\ci(\Lb,v)$ is the Weyl group 
of $C_\Gb^\ci(v)$ relative to $\Zb(\Lb)^\ci$. This shows that 
$W_\Gb^\ci(\Lb,v)$ is generated by elements of order $2$. So, since 
$|A_\Lb(v)|$ is supposed to be odd, the morphism 
$\ph_{\Lb,v}^\Gb : W_\Gb^\circ(\Lb,v) \to A_\Lb(v)$ is trivial. Therefore, 
by Corollary \ref{mor fti}, $H_\Gb(\Lb,v,1)=W_\Gb^\circ(\Lb,v)$.\fin 

\bi

\example{cuspidalite cuspidalite} We will see in Theorem \ref{alv} that, 
if the class $\Cb$ supports a cuspidal 
local system, then $A_\Lb(v)=A_\Gb(v)$ and $C_\Gb(u) \incl \Pb$. 
Consequently, the morphism $\ph_{\Lb,v}^\Gb$ is then well-defined.\finl

\bi

The morphism $\ph_{\Lb,v}^\Gb$ is the central object of this paper. 
In Part II, we will compute it explicitly whenever $v$ is a 
regular unipotent element under some restriction on $\Lb$. 

\bi

\sub{Separability} Let $\Cb^\et$ denote the separable closure of $\Cb$ in 
$\Lb/C_\Lb^\ci(v)$ (under the morphism $\Lb/C_\Lb^\ci(v) \to \Cb$, $l \mapsto 
lvl^{-1}$). Note that $\Cb^\et$ is smooth. 
The variety $\Cb^\et$ inherits from $\Lb/C_\Lb^\ci(v)$ 
the action of $\Lb$ by left translation, and the action of $A_\Lb(v)$ 
by right translation. Thus we have a sequence of $\Lb \times A_\Lb(v)$-equivariant 
morphisms
$$\diagram
\Lb/C_\Lb^\ci(v) \rrto && \Cb^\et \rrto && \Cb.
\enddiagram$$
The first morphism is bijective and purely inseparable, the second one 
is a Galois \'etale covering with group $A_\Lb(v)$. We define
$\Sigb^\et=\Zb(\Lb)^\ci \times \Cb^\et$, $\Sigb^\et_\reg=\Zb(\Lb)_\reg^\ci 
\times \Cb^\et$, $\Ybh^\et=\Gb \times \Sigb^\et_\reg$, and 
$\Ybt^\et=\Gb \times_\Lb \Sigb_\reg^\et$. We have a commutative diagram with 
cartesian squares
\equat\label{diagramme etale}
\diagram
\Sigb_\reg^\pr \ddto_{\DS{f_\reg^\ins}} && \Ybh' \llto_{\DS{\a'}} \rrto^{\DS{\b'}} 
\ddto_{\DS{\fha^\ins}} &&
 \Ybt' \ddto_{\DS{\fti^\ins}} \xto[4,2]^{\DS{\pi'}} & \\
 &&&&&&\\
\Sigb_\reg^\et \ddto_{\DS{f_\reg^\et}} && \Ybh^\et \llto_{\DS{\a^\et}} 
\rrto^{\DS{\b^\et}} 
\ddto_{\DS{\fha^\et}} &&
 \Ybt^\et \ddto_{\DS{\fti^\et}} \ddrrto_{\DS{\pi^\et}} & \\
&&&&&&\\
\Sigb_\reg && \Ybh \llto^{\DS{\a}} \rrto_{\DS{\b}} && \Ybt \rrto_{\DS{\pi}} && \Yb.
\enddiagram
\endequat
Here the maps $?^\et$ and $?^\ins$ are induced by the maps $?$ or $?'$. Moreover, 
all the morphisms $?^\et$ are Galois \'etale coverings, and all 
the morphisms $?^\ins$ are bijective purely inseparable morphisms. 
Note that $f^\et_\reg$, $\fha^\et$ and $\fti^\et$ are Galois coverings 
with group $A_\Lb(v)$ and $\pi^\et$ is a Galois covering 
with group $W_\Gb(\Lb,v)$. 

By the same argument as in Remark \ref{xtilde 0}, the group $\Pb$ acts 
on the variety $\Sigb^\et \times \Vb$ and the quotient $\Xbt_0^\et=\Gb \times_\Pb 
(\Sigb^\et \times \Vb)$ exists~: it is the separable closure of $\Xbt_0$ 
in $\Xbt_0^\pr$. If $(\Sigb^\et \times \Vb)_\mini$ denotes the inverse 
of $(\Sigb.\Vb)_\mini$ under the morphism $f^\et \times \Id_\Vb$, then 
$\Xbt_\mini^\et=\Gb \times_\Pb (\Sigb^\et \times \Vb)_\mini$ is the normalization 
of $\Xb_\mini$ in $\Ybt^\et$. So it inherits an action of $W_\Gb(\Lb,v)$ 
and the bijective purely inseparable morphism $\fti_\mini^\ins : \Xbt_\mini^\pr 
\to \Xbt_\mini^\et$ induced by $\fti_\mini$ 
is $W_\Gb(\Lb,v)$-equivariant. Moreover, the \mor $\fti_\mini^\et 
: \Xbt_\mini^\et \to \Xbt_\mini$ induced by $\fti_\mini$ 
is a Galois \'etale covering with group $A_\Lb(v)$. We summarize the notation 
in the following diagram.

$$
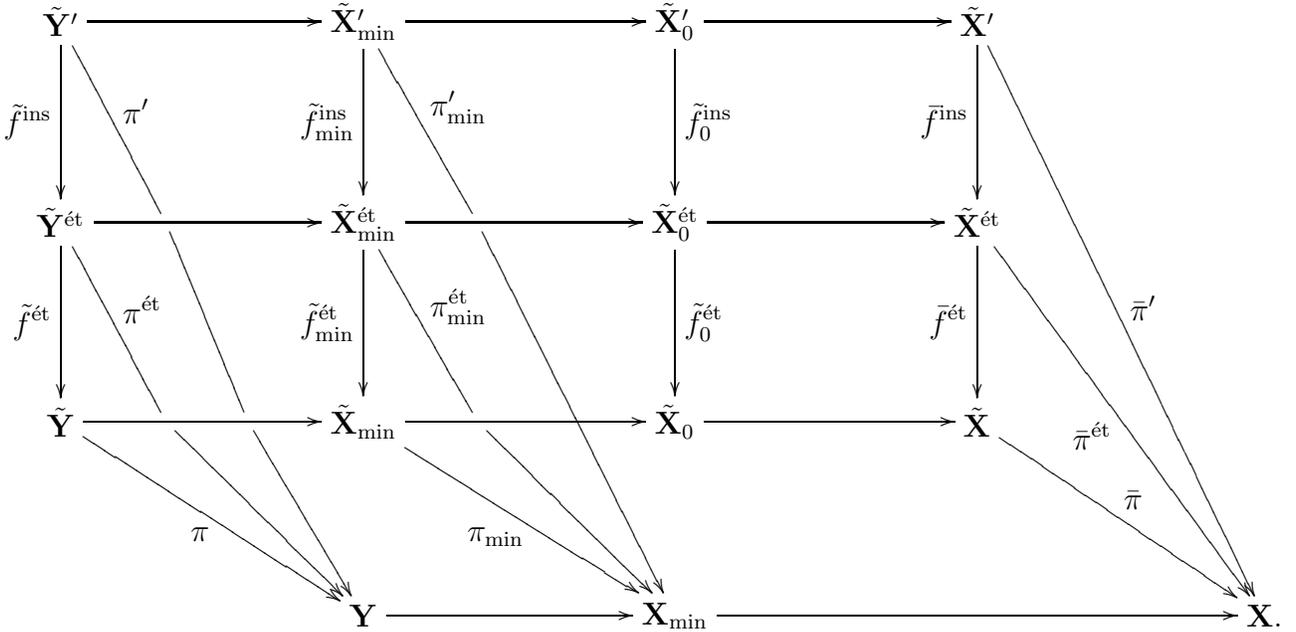
\diagram
\Ybt' \xto[0,3] \ddto_{\DS{\fti^\ins}}\xto'[2,1]^{\DS{\pi'}}'[4,2][6,3]&&& 
\Xbt_\mini^\pr \ddto_{\DS{\fti_\mini^\ins}}\xto'[2,1]^{\DS{\pi_\mini^\pr}}[6,3]
\xto[0,3] &&& \Xbt_0^\pr \ddto^{\DS{\fti_0^\ins}} \xto[0,3] &&& 
\Xbt' \ddto_{\DS{\fba^\ins}} \xto[6,3]^{\DS{\piba'}}\\
&&\\
\Ybt^\et \xto[0,3] \ddto_{\DS{\fti^\et}}\xto'[2,1]^{\DS{\pi^\et}}[4,3]&&& 
\Xbt_\mini^\et \ddto_{\DS{\fti_\mini^\et}}\xto'[2,1]^{\DS{\pi_\mini^\et}}[4,3]
\xto[0,3] &&& \Xbt_0^\et \ddto^{\DS{\fti_0^\et}} \xto[0,3] &&& 
\Xbt^\et \ddto_{\DS{\fba^\et}} 
\xto[4,3]_{\DS{\piba^\et}}&&& \\
&&\\
\Ybt \xto[0,3] \xto[2,3]_{\DS{\pi}} &&& \Xbt_\mini \xto[2,3]_{\DS{\pi_\mini}} \xto[0,3] 
&&& \Xbt_0  \xto[0,3] &&& \Xbt \xto[2,3]^{\DS{\piba}}\\
&&&\\
&&& \Yb \xto[0,3] &&& \Xb_\mini \xto[0,6] &&&&&&\Xb. \\
\enddiagram$$

\bi

\remark{stab etale} 
For $z \in \Zb(\Lb)^\ci$, we denote by $\uti_z^\et$ the image of $\uti_z^\pr \in 
\Xbt_\mini^\pr$ in $\Xbt_\mini^\et$ under the morphism $\fti_\mini^\ins$. 
Since $\fti_\mini^\ins$ is bijective and $W_\Gb(\Lb,v)$-equivariant, 
the stabilizer of $\uti_z^\et$ in $W_\Gb(\Lb,v)$ is equal to $H_\Gb(\Lb,v,z)$.\finl

\bi

\example{exemple etale} It may happen that the variety $\Cb^\et$ is different 
from $\Lb/C_\Lb^\ci(v)$, so that the construction above is not trivial. 
Of course, it only occurs in positive characteristic. The smallest example 
is given by the group $\Lb=\Gb=\Sb\Lb_2(\FM)$, whenever $p=2$ and
$$v=\matrice{1 & 1 \\ 0 & 1}.$$
Nevertheless, this is quite an unusual phenomenon. Indeed, 
if $\Gb=\Sb\Lb_n(\FM)$ and if $p$ does not divide $n$, then 
the morphism $\Lb/C_\Lb^\ci(v) \to \Cb$ is always \'etale. 
Also, if $\Gb$ is a quasisimple group of type different from $A$ and 
if $p$ is good for $\Gb$, then again the morphism $\Lb/C_\Lb^\ci(v) \to \Cb$ 
is always \'etale.\finl

\bi

\sec{Elementary properties of $\ph_{\Lb,v}^\Gb$\label{ele sec}}~

\med

As will be shown in \SEC\ref{sec endo}, knowing the morphism 
$\ph_{\Lb,v}^\Gb$ will be of fundamental 
importance in the description of the endomorphism algebra of an induced 
cuspidal character sheaf (\cf Corollary \ref{coro zeta phi}). For this reason 
we devote a section to gathering the properties of this morphism. These 
properties help to reduce the computations to small groups. 

\bi

\sub{Product of groups\label{produit}} We assume in this subsection, and 
only in this subsection, 
that $\Gb=\Gb_1 \times \Gb_2$. Let $\Lb=\Lb_1 \times \Lb_2$, $v=(v_1,v_2)$, 
$\Pb=\Pb_1 \times \Pb_2$. Then, for every $z =(z_1,z_2) \in \Zb(\Lb)^\ci$, we have
\equat
H_\Gb(\Lb,v,z)=H_{\Gb_1}(\Lb_1,v_1,z_1) \times H_{\Gb_2}(\Lb_2,v_2,z_2).
\endequat
Moreover, if $A_\Lb(v) =A_\Gb(v)$ and if $C_\Gb(u) \incl \Pb$, then 
\equat\label{produit H}
\ph_{\Lb,v}^\Gb = \ph_{\Lb_1,v_1}^{\Gb_1} \times \ph_{\Lb_2,v_2}^{\Gb_2}.
\endequat

\bi

\sub{Changing the group\label{changement}} Let $\s : \Gb_1 \to \Gb$ 
be an isotypic morphism between connected reductive groups. 
We put $\Lb_1=\s^{-1}(\Lb)$ and $\Pb_1=\s^{-1}(\Pb)$, so that $\Lb_1$ 
is a \levi of the \para $\Pb_1$ of $\Gb_1$. Let $v_1$ denote the unique 
unipotent element of $\Gb_1$ such that $\s(v_1)=v$ and let $\Cb_1$ 
be the class of $v_1$ in $\Gb_1$. Note that $\s(\Cb_1)=\Cb$. 

\bi

\begin{lem}\label{G1 G}
$({\mathrm{a}})$ $\s(C_{\Lb_1}(v_1)) \Zb(\Gb)^\ci=C_\Lb(v)$ and 
$\s^{-1}(C_\Lb(v)) = C_{\Lb_1}(v_1)$.

\tete{b} The morphism $A_{\Lb_1}(v_1) \to A_\Lb(v)$ induced by $\s$ is surjective.

\tete{c} The morphism $W_{\Gb_1}(\Lb_1,\Sigb_1) \to W_\Gb(\Lb,\Sigb)$ induced by $\s$ 
is an isomorphism.

\tete{d} The morphism $W_{\Gb_1}(\Lb_1,v_1) \to W_\Gb(\Lb,v)$ induced by $\s$ 
is surjective.

\tete{e} The morphism $W_{\Gb_1}^\ci(\Lb_1,v_1) \to W_\Gb^\ci(\Lb,v)$ induced by $\s$ 
is an isomorphism.

\tete{f} $C_\Gb(u) \incl \Pb$ if and only if $C_{\Gb_1}(u_1) \incl \Pb_1$.

\tete{g} $A_\Lb(v)=A_\Gb(v)$ if and only if $A_{\Gb_1}(v_1)=A_{\Lb_1}(v_1)$. 
\end{lem}

\bi

\proof First, note that (b) is an immediate consequence of (a). 
To prove (a), 
let $l \in C_\Lb(v)$. Then there exists $l_1$ in $\Lb_1$ and $z \in \Zb(\Gb)^\ci$ such 
that $\s(l_1)=lz$. So there exists $x \in \Ker \s$ such that $l_1v_1l_1^{-1}=xv_1$. 
But $l_1v_1l_1^{-1}$ is unipotent and $\Ker \s$ consists of central 
semisimple elements of $\Lb_1$, so $x=1$, that is 
$l_1$ centralizes $v_1$. Hence the first assertion of (a) follows. 
The second follows by the same argument. 

(c), (d), (e), (f) and (g) follow by similar arguments (note that 
$\s(C_{\Gb_1}^\ci(v_1)).\Zb(\Gb)^\ci=C_\Gb^\ci(v)$).\fin 

\bi

We will denote by a subscript $?_1$ the object associated to the datum $(\Lb_1,v_1)$ 
and defined in the same way as the object $?$ in $\Gb$ (for instance, $\Sigb_1=
\Zb(\Lb_1)^\ci .\Cb_1$, $\Xbt_1$, $\Ybt_1^\pr$...).

The reader may check that $\s((\Sigb_1.\Vb_1)_\mini) \incl (\Sigb.\Vb)_\mini$, so 
that $\s$ induces a morphism $\Xbt_{1,\min} \to \Xbt_\mini$. This morphism 
is $W_{\Gb_1}(\Lb_1,\Sigb_1)$-equivariant as it can be verified by restriction to 
$\Ybt_1$ (here, $W_{\Gb_1}(\Lb_1,v_1)$ and $W_\Gb(\Lb,\Sigb)$ are identified 
via the morphism $\s$ by Lemma \ref{G1 G} (c)). 
Similarly, $\s$ induces a $W_{\Gb_1}(\Lb_1,v_1)$-equivariant 
morphism $\Xbt_{1,\min}^\pr \to \Xbt_\mini^\pr$ (here, $W_{\Gb_1}(\Lb_1,v_1)$ 
acts on $\Xbt_\mini^\pr$ via the surjective morphism $W_{\Gb_1}(\Lb_1,v_1)
\to W_\Gb(\Lb,v)$). 

\bi

\begin{prop}\label{G1 G H}
If $z_1 \in \Zb(\Lb_1)^\ci$, then~:

\tete{a} $H_{\Gb_1}(\Lb_1,\Sigb_1,z_1)=
H_\Gb(\Lb,\Sigb,\s(z_1))$. 

\tete{b} $\s$ induces an isomorphism 
$H_{\Gb_1}(\Lb_1,v_1,z_1)\simeq H_\Gb(\Lb,v,\s(z_1))$.

Moreover,

\tete{c} If $C_\Gb(u) \incl \Pb$ and $A_\Lb(v)=A_\Gb(v)$, then 
the diagram
$$\diagram
W_{\Gb_1}^\ci(\Lb_1,v_1) \rrto_{\sim}^{\DS{\s}} \ddto_{\DS{\ph_{\Lb_1,v_1}^{\Gb_1}}}
&& W_\Gb^\ci(\Lb,v) \ddto^{\DS{\ph_{\Lb,v}^\Gb}} \\
&&\\
A_{\Lb_1}(v_1) \rrto^{\DS{\s}} && A_\Lb(v) 
\enddiagram$$
is commutative.
\end{prop}

\bi

\proof It is clear that $H_{\Gb_1}(\Lb_1,\Sigb_1,z_1) \incl 
H_\Gb(\Lb,\Sigb,\s(z_1))$. To prove the reverse inclusion, we may 
and do assume that $z_1=1$ (by Proposition \ref{stabilisateur z}). Then if an element 
$w \in W_\Gb(\Lb,\Sigb)$ stabilizes $1 *_\Pb u$, this proves that there exists 
$a \in \Ker \s$ such that $(1 *_{\Pb_1} u_1).w=1*_{\Pb_1} au_1$. But, by 
\ref{ZZZ}, $a=1$. This proves (a).

To prove (b), we notice that 
$\s(H_{\Gb_1}(\Lb_1,v_1,z_1) )\incl H_\Gb(\Lb,v,\s(z_1))$, so $\s$ induces 
a \mor $H_{\Gb_1}(\Lb_1,v_1,z_1) \longto  H_\Gb(\Lb,v,\s(z_1))$. This morphism 
is injective by \ref{H cap A}. It is surjective by the same argument as the 
one used in (a). The proof of (b) is complete.

By Lemma \ref{G1 G} (a), the diagram described in (c) is well-defined. 
Its commutativity follows from (b).\fin

\bi

\sub{Parabolic restriction\label{para sub red}} In this subsection, we show that 
the above constructions are compatible with ``restrictions'' to parabolic subgroups 
of $\Gb$. We need some notation. Let $\Pb'$ denote a \para of $\Gb$ containing 
$\Pb$, and let $\Lb'$ denote the unique \levi of $\Pb'$ containing $\Lb$. 
Let $\Vb'$ denote the unipotent radical of $\Pb'$, and let $v' =
\pi_{\Lb'}(u)$. Then $v' \in v(\Vb \cap \Lb')$. We start with some 
elementary properties.

\bi

\begin{prop}\label{fourre-tout}
$(1)$ $v' \in \Cb^{\Lb'}$.

\tete{2} If $A_\Lb(v) =A_\Gb(v)$, then $A_\Lb(v)=A_{\Lb'}(v)$. 

\tete{3} If $C_\Gb(u) \incl \Pb$, then $C_{\Lb'}(v') \incl \Pb \cap \Lb'$.
\end{prop}

\proof By Lemma \ref{dim p}, we have $\dim C_\Pb(u) \ge \dim C_{\Pb \cap \Lb'}(v') \ge 
\dim C_\Lb(v)$. But 
$$\dim C_\Pb(u)=\dim C_\Lb(v)$$ 
because $u \in \Cb^\Gb$. 
So, $\dim C_{\Pb \cap \Lb'}(v') = \dim C_\Lb(v)$. This proves 
that $\dim~(v')_{\Pb \cap \Lb'} = \dim \Cb + \dim \Vb \cap \Lb'$. 
So $(v')_{\Pb \cap \Lb'}$ is dense in $\Cb.(\Vb \cap \Lb')$, that is $v' \in \Cb^{\Lb'}$. 
Hence (1) is proved. 

(2) follows from \ref{cgo} applied to the \levi $\Lb'$. Let us now prove 
(3). Let $m \in C_\Lb'(v')$. We only need to prove that $m \in \Pb$. 
But $\lexp{m}{u} \in v'.\Vb' \cap \Cb^\Gb \incl v\Vb \cap \Cb^\Gb$. So, 
by \cite[Proposition II.3.2 (d)]{spaltenstein}, there exists $x \in \Pb$ 
such that $\lexp{m}{u}=\lexp{x}{u}$. So $x^{-1}m \in C_\Gb(u)$. 
But $C_\Gb(u) \incl \Pb$ by hypothesis, so $m \in \Pb$.\fin

\bi

\begin{prop}\label{prop compa para}
If $C_\Gb(u) \incl \Pb$, then $H_\Gb(\Lb,v,1) \cap W_{\Lb'}(\Lb,v)=H_{\Lb'}(\Lb,v,1)$.
\end{prop}

\bi

\proof By Proposition \ref{agu} (2), the subgroups 
$H_\Gb(\Lb,v,1) \cap W_{\Lb'}(\Lb,v)$ and $H_{\Lb'}(\Lb,v,1)$ 
have the same index in $W_{\Lb'}(\Lb,v)$ (this index is equal to $|A_\Lb(v)|$). 
Consequently, it is sufficient to prove that 
$$H_\Gb(\Lb,v,1) \cap W_{\Lb'}(\Lb,v) \incl H_{\Lb'}(\Lb,v,1).\leqno{(\#)}$$
Let 
$$\Fbt'=\Pb' \times_\Pb (\Sigb^\pr \times \Vb)_{\min,\Lb}^\Gb.$$
$\Fbt'$ is an irreducible closed subvariety of $\Xbt'$, and it is stable 
under the action of $W_{\Lb'}(\Lb,v)$ (indeed, the open subset $\Obt'=\Pb' \times_\Lb 
\Sigb_\reg^\pr$ is obviously $W_{\Lb'}(\Lb,v)$-stable). 

By Lemma \ref{dim p} (2), the projection $\pi_{\Lb'} : \Pb' \to \Lb'$ sends an element of 
$(\Sigb.\Vb)_{\mini,\Lb}^\Gb$ to an element of 
$(\Sigb.(\Vb \cap \Lb'))_{\min,\Lb}^{\Lb'}$. 
So, it induces a map $\g : \Fbt' \to (\Xbt_\mini^\pr)_\Lb^{\Lb'}$. 
Moreover, the diagram
$$\diagram
\Obt' \rrto \ddto &&\Fbt' \ddto^{\DS{\g}} \\
&&\\
(\Ybt')_\Lb^{\Lb'} \rrto && (\Xbt_\mini^\pr)_\Lb^{\Lb'}
\enddiagram$$
is commutative. The first vertical map is $W_{\Lb'}(\Lb,v)$-equivariant, so, by density, 
the second vertical map is also $W_{\Lb'}(\Lb,v)$ equivariant. This 
proves $(\#)$.\fin

\bi

\begin{coro}\label{coro compa para}
If $C_\Gb(u) \incl \Pb$ and if $A_\Lb(v)=A_\Gb(v)$, then
$$\ph_{\Lb,v}^{\Lb'}=\Res_{W_{\Lb'}^\ci(\Lb,v)}^{W_\Gb^\ci(\Lb,v)} \ph_{\Lb,v}^\Gb.$$
\end{coro}

\bi

\proof By Proposition \ref{fourre-tout} (2) and (3), 
$\ph_{\Lb,v}^{\Lb'}$ is well-defined. 
So Corollary \ref{coro compa para} follows from 
Proposition \ref{prop compa para}.\fin

\bi

\remark{qwerty} If $\Gb'$ is a connected reductive subgroup of $\Gb$ containing 
$\Lb$, then it may happen that 
$$\ph_{\Lb,v}^{\Gb'}\not=
\Res_{W_{\Gb'}^\ci(\Lb,v)}^{W_\Gb^\ci(\Lb,v)} \ph_{\Lb,v}^\Gb.$$
An example is provided in Part II of this paper.\finl

\bi

\sec{Endomorphism algebra of induced cuspidal character sheaves\label{sec endo}}~

\med

\noindent{\bf Hypothesis and notation :} {\it From now on, and until the end 
of this first part, we assume that $\Cb$ supports an irreducible cuspidal 
\cite[Definition 2.4]{luicc} local system $\EC$. 
To this local system is associated an irreducible character $\z$ of 
$A_\Lb(v)$, via the Galois \'etale covering $\Cb^\et \to \Cb$. 
Let $\FC = \qlb \boxtimes \EC$ ($\FC$ is a local system on $\Sigb$) and let 
$\FC_\reg$ denote the restriction of $\FC$ to $\Sigb_\reg$.
Let $K$ be the perverse sheaf on $\Gb$ obtained from the triple $(\Lb,v,\z)$ 
by parabolic induction \cite[4.1.1]{lucs} and let $\AC$ denote its endomorphism 
algebra.}

\bi

In \cite[Theorem 9.2]{luicc}, Lusztig constructed an isomorphism 
$\Th : \qlb W_\Gb(\Lb) \to \AC$. This \iso is very 
convenient for computing the generalized Springer correspondence. 
On the other hand, Lusztig's construction is canonical but not explicit. 
The principal aim of this paper is to construct 
an explicit \iso $\Th' : \qlb W_\Gb(\Lb) \to \AC$ by another 
method. It turns out that this \iso differs from Lusztig's one 
by a linear character $\g$ of $W_\Gb(\Lb)$. Knowledge of $\g$ 
would allow us to combine the advantages of both 
\isos $\Th$ and $\Th'$. However, we are not able to determine 
it explicitly in general, although we can relate it 
to the morphism $\ph_{\Lb,v}^\Gb$ defined in the previous section. 
Note that $\g=1$ whenever $\Lb=\Tb$. 

In this section we recall some well-known facts about cuspidal 
local systems, parabolic induction and endomorphism algebra. 
Most of these results may be found in \cite{luicc} or \cite{lucs}. 
However, Theorem \ref{alv}, which is a fundamental step 
for constructing the \iso $\Th'$ defined above, 
was proved in full generality in \cite{Bonnafe}. 
The \iso $\Th'$ will be constructed in \SEC\ref{section endo}. 
We will also prove in \SEC\ref{section endo} the existence of $\g$ 
and its relationship with $\ph_{\Lb,v}^\Gb$. 
We will provide in \SEC\ref{section 1 gamma} some properties of $\g$ which 
allow us to reduce its computation to the case where $\Gb$ is semisimple, 
simply-connected, quasi-simple, and $\Pb$ is a maximal \para of $\Gb$. 

\bi

\sub{Parabolic induction\label{subsec para}} For the convenience of the reader, 
we reproduce here the diagram \ref{diagramme central}.
$$\diagram
\Sigb_\reg \ddto && \Ybh \llto_{\DS{\a}} 
\rrto^{\DS{\b}} \ddto && \Ybt \rrto^{\DS{\pi}} \ddto && \Yb \ddto\\
&&&&&&\\
\overline{\Sigb} && \Xbh \llto_{\DS{\alpba}} \rrto^{\DS{\betba}} && \Xbt 
\rrto^{\DS{\piba}}  && \Xb.
\enddiagram$$
We define $\FCh_\reg = \a^* \FC_\reg$~: it is a local system on $\Ybh$. 
Moreover, since $\EC$ is $\Lb$-equivariant, there exists a local system $\FCt_\reg$ 
on $\Ybt$ such that $\b^* \FCt_\reg \simeq \FCh_\reg$.
By \cite[3.2]{luicc}, 
the \mor $\pi : \Ybt \to \Yb$ is a Galois covering with Galois group $W_\Gb(\Lb)
=N_\Gb(\Lb)/\Lb$, 
so $\pi_*\FCt_\reg=\pi_!\FCt_\reg$ (because $\pi$ is finite hence proper) is a 
local system on $\Yb$. We denote by $K$ the following perverse sheaf on $\Gb$~:
\equat
K=IC(\overline{\Yb},\KC)[\dim \Yb].
\endequat
where $\pi_* \FCt_\reg=\KC$. Recall that $\overline{\Yb}=\Xb$, so $\dim \Yb
 = \dim \Xb$.

We shall give, following \cite[\SEC 4]{luicc}, 
an alternative description of the perverse sheaf $K$. 
Let $A$ be the following perverse sheaf on $\Lb$~:
$$A=IC(\overline{\Sigb},\FC)[\dim \Sigb].$$
Note that $\overline{\Sigb}=\Zb(\Lb)^\ci \overline{\Cb}$ so $A=\qlb[\dim \Zb(\Lb)^\ci] 
\boxtimes IC(\overline{\Cb},\EC)[\dim \Cb]$. 
Since $A$ is $\Lb$-equivariant, 
there exists a perverse sheaf $\Kti$ on $\Xbt$ such that 
$$\alpba^*A[\dim \Gb + 
\dim \Vb] = \betba^*\Kti[\dim \Pb].$$
The perverse sheaf $\Kti$ is in fact 
equal to $IC(\Xbt,\FCt_\reg)[\dim \Xbt]$. By \cite[Proposition 4.5]{luicc}, we have
\equat\label{K}
K={\mathrm{R}}\piba_*\Kti.
\endequat

\bi

The fact that $\Cb$ admits a cuspidal 
local system has a lot of consequences. We gather some of them in the next 
two theorems.

\bi

\begin{theo}[{\bf Lusztig}]\label{rappellusztig}
\tete{a} $v$ is a distinguished unipotent element of $\Lb$.

\tete{b} $N_\Gb(\Lb)$ stabilizes $\Cb$ and $\EC$, so $W_\Gb(\Lb,\Sigb)=W_\Gb(\Lb)$.

\tete{c} $\Lb$ is universally self-opposed. 
\end{theo}

\bi

\proof \cf \cite[Proposition 2.8]{luicc} for (a) and \cite[Theorem 9.2]{luicc}  
for (b) and (c).\fin

\bi

The first assertion 
of the next theorem has been proved by Lusztig for $p$ 
large enough by using the classification of cuspidal pairs \cite{luhecke}. 
In \cite[Corollary to Proposition 1.1]{Bonnafe}, the author gave a proof 
of the general case without using the classification. The second assertion 
is proved in \cite[Theorem 3.1 (4)]{bonnafe cgu}.

\bi

\begin{theo}\label{alv}
We have $A_\Lb(v)=A_\Gb(v)$ and $C_\Gb(u) \incl \Pb$. 
\end{theo}

\bi

By Theorem \ref{alv} and Lemma \ref{alv agv alv agv}, we have 
\equat\label{W produit}
W_\Gb(\Lb,v)=A_\Lb(v) \times W_\Gb^\ci(\Lb,v),
\endequat
and, by Theorem \ref{rappellusztig} (b), we have 
\equat\label{W}
W_\Gb^\ci(\Lb,v) \simeq W_\Gb(\Lb)=N_\Gb(\Lb)/\Lb.
\endequat

\bi

\sub{Lusztig's description of $\AC$\label{sub endo lusztig}} For each $w$ in
$W_\Gb^\ci(\Lb,v)$, we choose a representative $\wdo$ of $w$ in $N_\Gb(\Lb) 
\cap C_\Gb^\ci(v)$. By Theorem \ref{rappellusztig} (b), the local 
systems $\FC$ and $(\INT \dot{w})^* \FC$ are isomorphic. Let $\th_w$ denote an 
\iso of $\Lb$-equivariant local systems $\FC \to (\INT \dot{w})^* \FC$.
Then $\th_w$ induces an \iso $\thet_w : \FCt_\reg \to \g_w^* \FCt_\reg$, where 
$\g_w : \Ybt \to \Ybt$, $(g,x\Lb) \mapsto (g,x\wdo^{-1}\Lb)$ 
(\cf \cite[proof of Proposition 3.5]{luicc}). 
Since $\pi_* \g_w^* = \pi_*$, $\pi_*\thet_w$ is an \auto of $\KC$. 
By applying the functor $IC(\Xb,.)[\dim \Yb]$, $\pi_* \thet_w$ induces 
an \auto $\Th_w$ of $K$. The \auto $\Th_w$, as well as $\th_w$, is 
defined up to multiplication by an element of $\qlb^\times$. \

By \cite[Proposition 3.5]{luicc}, $(\Th_w)_{w \in W_\Gb^\ci(\Lb,v)}$ is a basis 
of the endomorphism algebra $\AC$ of $K$~; 
moreover, by \cite[Remark 3.6]{luicc}, 
there exists a family $(a_{w,w'})_{w,w' \in W_\Gb^\ci(\Lb,v)}$ 
of elements of $\qlb^\times$ 
such that $\Th_w \Th_{w'}= a_{w,w'} \Th_{ww'}$ for all $w$ and $w'$ in 
$W_\Gb^\ci(\Lb,v)$. 
Lusztig proved that it is possible to choose in a 
canonical way the family $(\th_w)_{w \in W_\Gb^\ci(\Lb,v)}$ 
such that $\Th_w \Th_{w'}= \Th_{ww'}$ for all $w$ and $w'$ in 
$W_\Gb^\ci(\Lb,v)$. The next theorem \cite[Theorem 9.2]{luicc} 
explains his construction. 

\bi

\begin{theo}[{\bf Lusztig}]\label{lusztig theta}
There exists a unique family of \isos of local systems 
$(\th_w : \FC \to (\INT \dot{w})^* \FC)_{w \in W_\Gb^\ci(\Lb,v)}$ 
such that the following condition holds~:  
for each $w \in W_\Gb^\ci(\Lb,v)$, $\Th_w$ acts as the identity  
on $\HC^{- \dim \Yb}_u K$, where $u$ is any element of $\Cb^\Gb$. 
\end{theo}

\bi

In the previous theorem, the uniqueness of the family 
$(\th_w)_{w \in W_\Gb^\ci(\Lb,v)}$ follows from the fact that 
$\HC^{- \dim \Yb}_u K \not= 0$ for each $u \in \Cb^\Gb$. 
As a consequence, one gets that the linear mapping 
$$\fonction{\Th}{\qlb W_\Gb^\ci(\Lb,v)}{\AC=\End_{\MC\Gb}(K)}{w}{\Th_w}$$
is an \iso of algebras. 

\bigskip

\rem In \cite[\SEC 3]{luicc}, Lusztig starts with an element $w \in W_\Gb(\Lb)$ 
and a representative $\wdo$ of $w$ in $N_\Gb(\Lb)$. Here, we have slightly modified 
his argument using the fact that $W_\Gb^\circ(\Lb,v) \simeq W_\Gb(\Lb)$ 
which implies that every $w \in W_\Gb(\Lb)$ has a representative in 
$N_\Gb(\Lb) \cap C_\Gb^\circ(v)$. Such a choice of the representatives will allow 
us to provide another interesting family of isomorphisms 
$(\th_w' : \FC \to (\INT \dot{w})^* \FC)_{w \in W_\Gb^\circ(\Lb,v)}$ 
(see \SEC\ref{sous endo}).\finl 

\bigskip

If $\ch$ is an \irr \car of $W_\Gb^\ci(\Lb,v)$, we denote by $K_\ch$ 
an irreducible component of $K$ associated to $\ch$ via the 
\iso $\Th$. 

\bi

\begin{coro}[{\bf Lusztig}]\label{lusztig coro}
For each $u \in \Cb^\Gb$, we have~:

\tete{a} $\HC^{-\dim \Yb}_u K_1 = \HC^{-\dim \Yb}_u K$. 

\tete{b} $\HC^{-\dim \Yb}_u K_\ch = 0$ for every non-trivial \irr \car $\ch$ of
$W_\Gb^\ci(\Lb,v)$.
\end{coro}

\bi

\sub{Restriction to the open subset $\Xb_\mini$} 
The restriction $\KCt_0$ of $\Kti[-\dim \Xb]$ to $\Xbt_0$ 
is a local system \cite[4.4]{luicc}, 
that is, a complex concentrated in degree $0$, whose $0$th term 
is a local system. In fact, $\KCt_0$ is the local system on $\Xbt_0$ 
associated to the Galois \'etale covering $\Xbt_0^\et \to \Xbt_0$ and to 
the character $\z$ of $A_\Lb(v)$. Therefore, the restriction $\KCt_\mini$ 
of $\Kti[-\dim \Xb]$ to $\Xbt_\mini$ is a local system. More precisely, 
it is the local system associated to the Galois \'etale covering 
$\Xbt_\mini^\et \to \Xbt_\mini$ and to the character $\z$.
Let $K_{\min}$ denote the restriction of $K[- \dim \Xb]$ to $\Xb_{\min}$. 
We have the following result.

\bi

\begin{prop}\label{k min}
We have $K_{\min}=\pi_{\min,*} \KCt_{\min}$. So, $K_{\min}$ is a 
constructible sheaf, that is a complex concentrated in degree $0$.
\end{prop}

\bi

\proof Since $\pi_{\min}$ is finite, the functor $\pi_{\min,*}$ is exact. 
The proposition follows from this remark and the Proper Base Change Theorem.\fin

\bi

\sec{Another isomorphism between $\AC$ and 
$\qlb W_\Gb^\ci(\Lb,v)$\label{section endo}} 

\bi

The aim of this section is to construct an explicit \iso 
$\Th'$ between the endomorphism algebra $\AC$ and 
the group algebra $\qlb W_\Gb^\ci(\Lb,v)$. Our strategy is 
the following. First, note that the endomorphism algebra $\AC$ of $K$ is 
canonically isomorphic to the endomorphism algebra of the local 
system $\KC$ on $\Yb$. To this local system is associated 
a representation of the fundamental group $\pi_1(\Yb,y)$ of $\Yb$ 
(where $y$ is any point of $\Yb$). This representation and its endomorphism 
algebra are easy to describe (\cf \ref{iso psi}). 

\bi

\sub{Representations of the fundamental group\label{sous endo}} 
Let $V_\z$ denote an \irr left $\qlb A_\Lb(v)$-module affording the character $\z$. 
We may, and we will, assume that
$$\FC=(f^\et)_* \qlb \otimes_{\qlb A_\Lb(v)} V_\z,$$
$$\FC_\reg=(f_\reg^\et)_* \qlb \otimes_{\qlb A_\Lb(v)} V_\z,$$
$$\FCh_\reg=(\fha_\reg^\et)_* \qlb \otimes_{\qlb A_\Lb(v)} V_\z,$$
$$\FCt_\reg=(\fti_\reg^\et)_* \qlb \otimes_{\qlb A_\Lb(v)} V_\z.\leqno{\mathrm{and}}$$
Here, $V_\z$ is identified with the constant sheaf with values in $V_\z$. 
From the fourth equality, we deduce that
$$\KC=\pi_* \FCt_\reg=\pi_*^\et \qlb \otimes_{\qlb W_\Gb(\Lb,v)} 
\Ind_{A_\Lb(v)}^{W_\Gb(\Lb,v)} V_\z.$$
Therefore, the \endo algebra of $\KC$ is canonically isomorphic to the endomorphism 
algebra of the $\qlb W_\Gb(\Lb,v)$-module $\Ind_{A_\Lb(v)}^{W_\Gb(\Lb,v)} V_\z$. 
But, by \ref{W produit}, this endomorphism algebra is canonically isomorphic 
to $\qlb W_\Gb^\ci(\Lb,v)$. 
Since the functor $IC(\overline{\Yb}, .)[\dim \Yb]$ 
is fully faithful, it induces an \iso 

\equat\label{iso psi}
\Th' : \qlb W_\Gb^\ci(\Lb,v) \longto \AC.
\endequat

\sma

This \iso may be constructed in another way. 
The action of an element $\wdo \in N_\Gb(\Lb) \cap C_\Gb^\ci(v)$ on 
$\Sigb_\reg^\et$, $\Ybh^\et$, and $\Ybt^\et$ commutes with the action of 
$A_\Lb(v)$. Therefore, there exists an \iso $\th_w^\pr : 
\FC \to (\INT \dot{w})^* \FC$ (\resp $\t_w^\pr : 
\FCt_\reg \to (\INT \dot{w})^* \FCt_\reg$) which induces the identity on the stalks 
at $zv$ (\resp $1 *_\Lb zv$) for every $z \in \Zb(\Lb)^\ci$ (\resp for every 
$z \in \Zb(\Lb)^\ci_\reg$). Then 
\equat\label{tau theta}
\b^*\t_w^\pr = \a^* \th_w^\pr|_{\Sigb_\reg}.
\endequat
Now, let $\Th_w^\pr =IC(\overline{\Yb},\pi_* \t_w^\pr)[\dim \Yb] : K 
\longmapright{\sim} K$. By \ref{tau theta}, 
there exists an element $\g_w \in \qlb^\times$ such that
\equat\label{gamma g}
\Th(w)=\g_w \Th_w^\pr.
\endequat
By looking at the action on the stalk at $zv \in \Yb$, one can immediately get 
the following result.

\bi 

\begin{prop}\label{wouaw}
With the above notation, we have $\Th^\pr_w=\Th'(w)$ for every 
$w \in W_\Gb^\ci(\Lb,v)$.
\end{prop}

\bi

\begin{coro}\label{epsilon}
There exists 
a linear character $\g_{\Lb,v,\z}^\Gb$ of the Weyl group $W_\Gb^\ci(\Lb,v)$ such that
$$\Th^\pr(w)=\g_{\Lb,v,\z}^\Gb(w) \Th(w)$$
for every $w \in W_\Gb^\ci(\Lb,v)$.
\end{coro}

\bi

\proof This follows from Theorem \ref{lusztig theta}, from \ref{gamma g}, 
and from Proposition \ref{wouaw}.\fin

\bi

\example{wa} Whenever $\Gb$ is symplectic or orthogonal, $\Lb$ is quasisimple 
and $p \not= 2$, Waldspurger 
\cite[\SEC VIII.8]{waldspurger} has considered an explicit subgroup 
of $W_\Gb(\Lb,v)$ (isomorphic to $W_\Gb^\ci(\Lb,v)$ by the canonical projection) 
constructed by ad hoc methods and has computed explicitly its action 
on the perverse sheaf $K$. It turns out 
that this subgroup coincides with $W_\Gb^\ci(\Lb,v)$ whenever $\Gb$ 
is symplectic, but does not whenever $\Gb$ is orthogonal. 
In other words, he computed explicitly the linear character 
$\g_{\Lb,v,\z}^\Gb$, up to the difference between our conventions. 
We give here his result. 

We assume that $\Lb \not= \Tb$ (if $\Lb=\Tb$, then $\g_{\Lb,v,\z}^\Gb=1$ 
by Corollary \ref{tore} below).
First, note that $W_\Gb^\ci(\Lb,v)$ is always a Weyl group of type $C$. We denote by 
$\g$ the linear character of $W_\Gb^\ci(\Lb,v)$ defined as follows~: 
it is the non-trivial linear character different from the sign character 
with value $-1$ on the reflection corresponding to the minimal Levi 
subgroup containing $\Lb$ having an irreducible root system 
(see Proposition \ref{propriete distingues} (a)).  
Then, by \cite[Lemma VIII.9]{waldspurger}, we have, for $p \not= 2$~:

\sma

\tete{a} {\it If $\Gb \simeq \Sb\Ob_n(\FM)$ and if $\Lb \simeq \Sb\Ob_{k^2}(\FM) \times 
(\FM^\times)^{n-k^2 \over 2}$, then $\g_{\Lb,v,\z}^\Gb=\g$.} 

\tete{b} {\it If $\Gb\simeq \Sb\pb_{2n}(\FM)$ and if $\Lb \simeq 
\Sb\pb_{k(k+1)} \times (\FM^\times)^{n-{k(k+1) \over 2}}$, 
then $\g_{\Lb,v,\z}^\Gb=\g^k$.}

\sma

In the second part of this paper, the case where $v$ is a regular unipotent 
element of $\Lb$ and $p$ is good for $\Gb$ will be treated. 
If $\Gb$ is symplectic and $\Lb$ is 
of type $C_1=A_1$, or if $\Gb$ is special orthogonal in even dimension 
and $\Lb$ is of type $D_2=A_1 \times A_1$, then $v$ is a regular 
unipotent element of $\Lb$~: the reader can check that 
the result given by Waldspurger 
\cite[Lemma VIII.9]{waldspurger} (and stated above with our convention) 
coincides with ours.\finl

\bi

\remark{class 2} 
By the classification of cuspidal local systems in good characteristic
\cite[\SEC 10-15]{luicc}, using Waldspurger's result \cite[Lemma VIII.9]{waldspurger} 
and our result whenever $v$ is regular (see Part II), the only 
case which remains is when $\Gb$ is a spin group and $v$ is not regular.\finl

\bi
 
If $\ch$ is an \irr \car of $W_\Gb^\ci(\Lb,v)$, we 
denote by $K_\ch^\pr$ an irreducible component of $K$ associated to $\ch$ 
via the \iso $\Th'$. By Corollary \ref{epsilon}, we have
\equat\label{gamma idiot}
K_\ch^\pr=K_{\g_{\Lb,v,\z}^\Gb \ch}.
\endequat

\bi

\sub{Links between $\g_{\Lb,v,\z}^\Gb$ and $\ph_{\Lb,v}^\Gb$} 
The following proposition is an immediate consequence of 
\ref{gamma idiot} and Corollary \ref{lusztig coro}~:

\bi

\begin{prop}\label{carac gamma}
The linear character $\g_{\Lb,v,\z}^\Gb$ of $W^\ci_\Gb(\Lb,v)$ is the 
unique \irr \car $\g$ of $W_\Gb^\ci(\Lb,v)$ 
satisfying $\HC_u^{-\dim \Yb} K_\g^\pr \not= 0$ for some (or any) $u \in \Cb^\Gb$.
\end{prop}

\bi

If $\ch$ is an irreducible \car of $W_\Gb^\ci(\Lb,v)$, we denote by 
$K_{\mini,\ch}$ (\resp $K_{\mini,\ch}^\pr$) the irreducible component 
of $K_{\mini}$ associated to $\ch$ via the \iso $\Th$ (\resp $\Th'$). 
Now, let $V_\ch$ be an irreducible $\qlb W_\Gb^\ci(\Lb,v)$-module 
affording $\ch$ 
as character. Then, since $W_\Gb(\Lb,v)$ acts on $\Xbt_\mini^\et$, 
it also acts on the constructible sheaf $(\pi_\mini^\et)_* \qlb$ and 
by construction we have 
\equat\label{equation K}
K_{\min,\ch}^\pr=(\pi^\et_\mini)_* \qlb \otimes_{\qlb W_\Gb(\Lb,v)} 
(V_\ch \otimes V_\z).
\endequat
If $x \in \Xb_\mini$, then $((\pi^\et_\mini)_* \qlb)_x$ 
is isomorphic, as a right $\qlb W_\Gb(\Lb,v)$-module, to 
the permutation module associated to the set $H_x \backslash W_\Gb(\Lb,v)$ 
(here, $H_x$ denotes the stabilizer, in $W_\Gb(\Lb,v)$, of 
a preimage of $x$ in $\Xbt_\mini'$). 
Since $u \in \Xb_\mini$, the stalk of $\HC_u^{-\dim \Yb} K_\ch^\pr$ at $u$ may easily 
be computed from this remark and from \ref{equation K}~: we have
$$\dim_{\qlb} \HC_u^{-\dim \Yb} K_\ch^\pr= \dim (K_{\min,\ch}^\pr)_u = 
\langle 
\Res_{H_\Gb(\Lb,v,1)}^{W_\Gb(\Lb,v)} (\ch \otimes \z), \boldsymbol{1}_{H_\Gb(\Lb,v,1)} 
\rangle$$
From this and from Proposition \ref{carac gamma}, we deduce immediately the 
following proposition.

\bi

\begin{prop}\label{zeta phi}
The linear character $\g_{\Lb,v,\z}^\Gb$ is the unique irreducible 
character $\g$ of $W_\Gb^\ci(\Lb,v)$ such that $\langle 
\Res_{H_\Gb(\Lb,v,1)}^{W_\Gb(\Lb,v)} (\g \otimes \z), \boldsymbol{1}_{H_\Gb(\Lb,v,1)} 
\rangle \not= 0$. 
\end{prop}

\bi

\begin{coro}\label{coro zeta phi}
We have  
$$\g_{\Lb,v,\z}^\Gb={1 \over \z(1)} \z \ci \ph_{\Lb,v}^\Gb.$$
\end{coro}

\bi

\proof Since $A_\Lb(v)=A_\Gb(v)$ and $C_\Gb(u) \incl \Pb$ 
(see Theorem \ref{alv}), the morphism $\ph_{\Lb,v}^\Gb : W_\Gb^\ci(\Lb,v) 
\to A_\Lb(v)$ is well-defined. Hence, the corollary follows from 
Proposition \ref{zeta phi} and Corollary \ref{mor fti}.\fin

\bi

\begin{coro}\label{odd}
If $|A_\Lb(v)|$ is odd, then $\g_{\Lb,v,\z}^\Gb=1$.
\end{coro}

\bi

\begin{coro}\label{tore}
$\g_{\Tb,1,1}^\Gb=1$.
\end{coro}

\bi

\example{regular example} Assume in this example that $v$ is a {\it regular} 
unipotent element of $\Lb$. In this case, $A_\Lb(v)$ is abelian \cite{springer}
so $\z$ is a linear character. 
By Corollary \ref{coro zeta phi}, we get
$$\g_{\Lb,v,\z}^\Gb=\z \ci \ph_{\Lb,v}^\Gb.$$
This case will be studied in full details in Part II.\finl

\bi

\noindent{\bf Questions :} Corollary \ref{coro zeta phi} shows that the image 
of the morphism $\ph_{\Lb,v}^\Gb$ is contained in the center of the character 
$\z$. On the other hand, the image of the canonical morphism 
$\jmath_v : \Zb(\Lb)/\Zb(\Lb)^\circ \to A_\Lb(v)$ is obviously contained in the center 
of $A_\Lb(v)$. Does there always exist a morphism 
$\psi_{\Lb,v}^\Gb : W_\Gb^\circ(\Lb,v) \simeq W_\Gb(\Lb) \to \Zb(\Lb)/\Zb(\Lb)^\circ$ 
such that $\ph_{\Lb,v}^\Gb = \jmath_v \circ \psi_{\Lb,v}^\Gb$~? 
Is it exists, does it really depend on $v$~? 

\bigskip

\sec{Elementary properties of the character $\g_{\Lb,v,\z}^\Gb$\label{section 1 gamma}}~

\med

\sub{Product of groups} We assume in this subsection that $\Gb=\Gb_1 \times \Gb_2$ 
where $\Gb_1$ and $\Gb_2$ are reductive groups. Let $\Lb=\Lb_1 \times 
\Lb_2$, $v=(v_1,v_2)$, $\Cb=\Cb_1 \times \Cb_2$ and $\z=\z_1 \otimes \z_2$. 
Then it is clear that 
\equat
W_\Gb^\ci(\Lb,v) = W_{\Gb_1}^\ci(\Lb_1,v_1) \times W_{\Gb_2}^\ci(\Lb_2,v_2)
\endequat
and that 
\equat\label{gamma produit}
\g_{\Lb,v,\z}^\Gb = \g_{\Lb_1,v_1,\z_1}^{\Gb_1} \otimes 
\g_{\Lb_2,v_2,\z_2}^{\Gb_2}.
\endequat

\bi

\sub{Changing the group\label{subsection change}} We use here the notation introduced 
in \SEC\ref{changement}. In particular, we still denote by a subscript $?_1$ 
the object in 
$\Gb_1$ corresponding to $?$ in 
$\Gb$ (e.g. $\Lb_1$, $\FC_1$, $\FC_{\reg,1}$, $\Ybh_1$, $K_1$, $\AC_1$, $\Th_1$\dots). 
We have $\s^{-1}(\Yb)=\Yb_1$. 
Note that the groups $W_{\Gb_1}^\ci(\Lb_1,v_1)$ and $W_\Gb^\ci(\Lb,v)$ are 
isomorphic via $\s$ by Lemma \ref{G1 G} (e). 

\bi

\begin{lem}
The local system $\s^*\EC_1$ on $\Cb_1$ is cuspidal. It is associated 
to the irreducible character $\z_1$ of $A_{\Lb_1}(v_1)$ obtained from $\z$ by composing 
with the surjective \mor $A_{\Lb_1}(v_1) \to A_\Lb(v)$ $($\cf Lemma 
\ref{G1 G} ${\mathrm{(b))}}$.
\end{lem}

\bi

\proof 
This is immediate from the alternative definition of a cuspidal local system given 
in terms of permutation representations \cite[Introduction]{luicc}.\fin

\bi

\begin{prop}\label{K1}
\tete{a} The restriction of $\s^*(K)$ to $\Yb_1$ is isomorphic to $K_1$.

\tete{b} $\s$ induces an isomorphism $\dot{\s} : \AC_1 \simeq \AC$.

\tete{c} The diagrams
$$\diagram
\qlb W_{\Gb_1}^\ci(\Lb_1,v_1) \rrto^{\qquad\DS{\Th_1}} \ddto_{\DS{\s}}&& \AC_1
\ddto^{\DS{\dot{\s}}}\\
&&\\
\qlb W_\Gb^\ci(\Lb,v) \rrto^{\qquad\DS{\Th}} && \AC\\
\enddiagram$$
$$\diagram
\qlb W_{\Gb_1}^\ci(\Lb_1,v_1) \rrto^{\qquad\DS{\Th^\pr_1}} \ddto_{\DS{\s}}&& \AC_1
\ddto^{\DS{\dot{\s}}}\\
&&\\
\qlb W_\Gb^\ci(\Lb,v) \rrto^{\qquad\DS{\Th'}} && \AC\\
\enddiagram\leqno{\mathrm{\it and}}$$
are commutative.
\end{prop}

\proof (a) follows from the Proper Base Change Theorem 
\cite[Chapter VI, Corollary 2.3]{Milne} applied 
to the Cartesian square 
$$\diagram
\Ybt_1 \rrto^{\DS{\pi_1}} \ddto_{\DS{\sigt}} && \Yb_1 \ddto^{\DS{\s}} \\
&&\\
\Ybt \rrto_{\DS{\pi}} &&\Yb.
\enddiagram$$
(b) follows from Lusztig's description of $\AC$. The commutativity of the first 
diagram in (c) follows from the fact that
$\s(\Cb_1^{\Gb_1})=\Cb^\Gb$ and from Theorem \ref{lusztig theta},  
while the commutativity of the second one follows from Proposition \ref{wouaw}.\fin

\bi

\begin{coro}\label{coro gamma}
We have $\g_{\Lb,v,\z}^\Gb\ci \s=\g_{\Lb_1,v_1,\z_1}^{\Gb_1}$.
\end{coro}

\bi

\sub{Parabolic restriction} Let 
$\Qb$ be a \para of $\Gb$ containing $\Pb$ and let 
$\Mb$ be the \levi of $\Qb$ containing $\Lb$. It follows from 
\cite[Theorem 8.3 (b)]{luicc} that~:

\bigskip

\begin{prop}\label{restriction parabolique}
$\g_{\Lb,v,\z}^\Mb=\Res_{W_\Mb^\ci(\Lb,v)}^{W_\Gb^\ci(\Lb,v)} \g_{\Lb,v,\z}^\Gb$.
\end{prop}

\bi

\remark{non restriction} 
If $\Gb'$ is a connected reductive subgroup of $\Gb$ which contains 
$\Lb$, then it may happen that 
$\g_{\Lb,v,\z}^{\Gb'}\not=
\Res_{W_{\Gb'}^\ci(\Lb,v)}^{W_\Gb^\ci(\Lb,v)} \g_{\Lb,v,\z}^\Gb$. 
An example is provided by the group $\Gb=\Sb\pb_4(\FM)$, as it will be shown in 
Part II of this paper.\finl

\bi

\remark{reduction maxi} Using \ref{gamma produit}, 
Proposition \ref{restriction parabolique}, Corollary \ref{coro gamma} 
and Theorem \ref{rappellusztig} (b), the computation of $\g_{\Lb,v,\z}^\Gb$ can be 
essentially reduced to the following case~: $\Gb$ is semisimple, simply 
connected, quasi-simple and $\Lb$ is a \levi of a maximal \para of $\Gb$.\finl

\bi

\rem Corollary \ref{coro zeta phi} can be used to give alternative 
proofs of \ref{gamma produit}, Proposition \ref{K1}, and Proposition 
\ref{restriction parabolique} (as consequences of \ref{produit H}, 
Proposition \ref{G1 G H}, and Proposition \ref{prop compa para} respectively).\finl

\bi

\sec{Introducing Frobenius\label{part finite}}

\bi

\sub{Hypothesis and notation} 
In this section, and only in this section, we assume that $\FM$ is an 
algebraic closure of a finite field. In particular, $p > 0$. We fix a power 
$q$ of $p$ and we denote by $\fq$ the subfield of $\FM$ with $q$ elements. 
We assume also that $\Gb$ is defined over $\fq$ and we denote by 
$F : \Gb \to \Gb$ the corresponding Frobenius endomorphism. 
If $g \in \Gb^F$, we denote by $[g]$ (or $[g]_{\Gb^F}$ if 
necessary) the $\Gb^F$-conjugacy class of $g$.

We keep the notation introduced in \SEC\ref{sec endo} ($\Lb$, 
$\Cb$, $v$, $\EC$, $K$, $\Th$, $\g_{\Lb,v,\z}^\Gb$...). We assume that 
$\Lb$ is $F$-stable. Then, by Theorem \ref{rappellusztig} (e), there exists 
$n \in N_\Gb(\Lb)$ such that $F(\Pb)=\lexp{n}{\Pb}$. Now, by Lang's  
theorem, we can pick an element $g \in \Gb$ such that $g^{-1}F(g)=n^{-1}$. 
Then $\lexp{g}{\Lb}$ and $\lexp{g}{\Pb}$ are $F$-stable. Since we are 
interested in the family of all $F$-stable subgroups of $\Gb$ which are 
conjugate to $\Lb$ under $\Gb$, we may, and we will assume that $\Lb$ and 
$\Pb$ are both $F$-stable. Without loss of generality, me may also 
assume that $\Bb$ and $\Tb$ are $F$-stable.

We also assume that $v$ and $\EC$ are $F$-stable.
Let $w \in W_\Gb^\ci(\Lb,v)$. We choose an element $g_w \in \Gb$ such 
that $g_w^{-1}F(g_w)=\dot{w}^{-1}$ (recall that $\dot{w}$ is a 
representative of $w$ in $N_\Gb(\Lb) \cap C_\Gb^\ci(v)$). We then put~:
$$\Lb_w=\lexp{g_w}{\Lb},\quad\quad v_w=\lexp{g_w}{v},
\quad\quad \Cb_w=\lexp{g_w}{\Cb},$$
$$\EC_w=(\ad g_w^{-1})^*\EC,\quad\quad{\mathrm{and}}
\quad\quad \FC_w=(\ad g_w^{-1})^*\FC.$$
Then $\Lb_w$ is an $F$-stable \levi of a \para of $\Gb$, $v_w \in \Lb_w^F$ 
is conjugate to $v$ in $\Gb^F$ (because $g_w^{-1}F(g_w) \in C_\Gb^\ci(v)$), 
$\Cb_w$ is the conjugacy class of $v_w$ in 
$\Lb_w$, $\EC_w$ is an $F$-stable cuspidal local system on $\Cb_w$ and 
$\FC_w=\qlb \boxtimes \EC_w$ (as a local system on 
$\Sigb_w=\Zb(\Lb_w)^\ci \times \Cb_w$).

\bi

\sub{Two conjugacy results} In \cite[Proposition 2.1]{Bonnafe}, we proved the 
following result~:

\bi

\begin{prop}\label{conjugue M}
Let $\Mb$ and $\Mb'$ be two $F$-stable \levis of $($non necessarily $F$-stable$)$ 
\paras of $\Gb$ which are geometrically conjugate and let $u'$ be a unipotent element of 
$\Mb^F$. Assume that the following conditions hold~:

\tete{a} $u'$ is a distinguished element of $\Mb$,

\tete{b} $N_\Gb(\Mb)$ stabilizes the class $(u')_\Mb$, and

\tete{c} $A_\Mb(u') = A_\Gb(u')$.

\noindent Then $[u']_{\Gb^F} \cap \Mb'$ is a single $\Mb^{\pr F}$-conjugacy class.
\end{prop}

\bi

\begin{coro}\label{conjugue w}
If $w \in W_\Gb^\ci(\Lb,v)$, then $[v]_{\Gb^F} \cap \Lb_w^F=[v_w]_{\Lb_w^F}$.
\end{coro}

\bi

\proof This follows from Theorem \ref{rappellusztig} (a) and (c), from 
Theorem \ref{alv} and from the previous Proposition \ref{conjugue M}.\fin

\bi

The next result implies \cite[II, 9.10.2]{lucs}.

\bigskip

\begin{coro}\label{ngm}
Let $w \in W_\Gb^\ci(\Lb,v)$. Then $N_{\Gb^F}(\Lb_w)$ stabilizes $[v_w]_{\Lb_w^F}$.
\end{coro}

\bigskip

\sub{Characteristic functions} We choose once and for all an isomorphism of local 
systems $\ph : F^*\EC \to \EC$ and we 
denote by $\XC_{\EC,\ph}$ the class function on $\Lb^F$ defined by
$$\XC_{\EC,\ph}(l)=\left\{\begin{array}{ll}
\Tr(\ph_l,\EC_l) & {\mathrm{if~}}l \in \Cb^F, \\
0 & {\mathrm{otherwise}},
\end{array}\right.$$
for any $l \in \Lb^F$. Using the isomorphism $\Th : \qlb W_\Gb^\ci(\Lb,v) \to \AC$, 
Lusztig defined  an isomorphism of local systems $\ph_w : F^*\EC_w \to \EC_w$. 
We recall his construction \cite[9.3]{lugf}. Let $\th_w : 
\EC \to (\ad \dot{w})^* \EC$ 
be the isomorphism of local systems defined in Theorem \ref{lusztig theta}. 
Then $\th_w$ induces an isomorphism of local systems 
$$F^* \ci (\ad g_w^{-1})^* \th_w : F^* \EC_w \longto (\ad g_w^{-1})^* \ci F^* \EC.$$
Moreover, $\ph$ induces an isomorphism 
$$(\ad g_w^{-1})^* \ph : (\ad g_w^{-1})^* \ci F^* \EC \longto \EC_w.$$
By the composition of the two previous isomorphisms, we get an isomorphism
$$\ph_w : F^* \EC_w \longto \EC_w.$$
Once the isomorphism $\ph$ is chosen, the isomorphism $\ph_w$ depends only 
on the construction of the isomorphism $\th_w$. By Corollary 
\ref{epsilon}, the knowledge 
of the isomorphism $\th_w$ is equivalent to the knowledge of the linear character 
$\g_{\Lb,v,\z}^\Gb$.

\bi

If $\ch$ is an $F$-stable 
irreducible character of $W_\Gb^\ci(\Lb,v)$, then we denote by $\chit$ 
the {\it preferred extension} of $\ch$ to $W_\Gb^\ci(\Lb,v) \rtimes <F>$ 
(the preferred extension has been defined by Lusztig \cite{lucs}). 
The choice of $\ph$ and $\chit$ induces a well-defined isomorphism $\ph_\ch : 
F^*K_\ch \longmapright{\sim} K_\ch$. We set, for every $g \in \Gb^F$, 
$$\XC_{K_\ch,\ph_\ch}(g)=\left\{\begin{array}{ll}
\DS{\sum_{i \in \ZM} (-1)^i \Tr(\HC_g^i(\ph_\ch),\HC_g^i K_\ch)} & {\mathrm{if}}~g
~{\mathrm{is}~} unipotent,\\
0 & {\mathrm{otherwise}}
\end{array}\right.$$
The importance of the knowledge of the characteristic functions $\XC_{\EC_w,\ph_w}$ 
(equivalently, of the linear character $\g_{\Lb,v,\z}^\Gb$) is given by 
the following theorem of Lusztig.

\bi

\begin{theo}[{\bf Lusztig}]\label{formule lusztig}
If $p$ is almost good for $\Gb$ and if $q$ is large enough, we have
$$\XC_{K_\ch,\ph_\ch} = {1 \over |W_\Gb^\ci(\Lb,v)|} 
\sum_{w \in W_\Gb^\ci(\Lb,v)} \chit(wF) R_{\Lb_w}^\Gb (\XC_{\EC_w,\ph_w}).$$
\end{theo}

\bi

We conclude this section by explaining how the characteristic 
functions $\XC_{\EC_w,\ph_w}$ can be computed explicitly 
once the linear character $\g_{\Lb,v,\z}^\Gb$ is known 
(see \cite[\SEC VIII]{waldspurger} for a particular case).
It follows from Lang's theorem that the set of rational conjugacy classes 
contained in $\Cb_w^F$ is in one-to-one correspondence with $H^1(F,A_\Lb(v_w))
\simeq H^1(\wdo F,A_\Lb(v)) =H^1(F,A_\Lb(v))$ (the last equality 
follows from the fact that $W_\Gb^\ci(\Lb,v)$ acts trivially 
on $A_\Lb(v)$). Let $a \in H^1(F,A_\Lb(v))$. We denote by $\aha$ 
a representative of $a$ in $A_\Lb(v)$ and by $v_{w,a}$ a representative 
of the rational conjugacy class contained in $\Cb_w^F$ parameterized by $a$. 
If $w=1$, we denote by $v_a$ the element $v_{w,a}$. We have $v_a \in \Lb^F$.  
It must be noticed that $[v_{w,a}]_{\Lb_w^F}=[v_a]_{\Gb^F} \cap \Lb^F$ 
(\cf Corollary \ref{conjugue w}).

By following step by step the construction of the \isos 
$\ph_w$, we obtain that the link between the class functions 
$\XC_{\EC,\ph}$ and $\XC_{\EC_w,\ph_w}$ 
is given in terms of the linear character $\g_{\Lb,v,\z}^\Gb$. 
More precisely, we get~:

\bi

\begin{prop}\label{fonction}
Let $w \in W_\Gb^\ci(\Lb,v)$ and let $a \in H^1(F,A_\Lb(v))$. Then 
$$\XC_{\EC_w,\ph_w}(v_{w,a})=\XC_{\EC,\ph}(v_a)\g_{\Lb,v,\z}^\Gb(w).$$
\end{prop}

\bi

Assume now until the end of this section that $\z$ is a linear character. 
In this case, we have
\equat\label{je sais pas}
\XC_{\EC_w,\ph_w}(v_{w,a})=\XC_{\EC_w,\ph_w}(v_w)\z(\aha).
\endequat
Note that $\z(\aha)$ does not depend on the choice of $\aha$ 
because $\z$ is $F$-stable. Hence, we deduce from Proposition 
\ref{fonction} the following result~:

\bi

\begin{coro}\label{description chi}
Assume that $\z$ is a linear character, and let $w \in W_\Gb^\ci(\Lb,v)$ 
and $a \in H^1(F,A_\Lb(v))$. Then 
$$\XC_{\EC_w,\ph_w}(v_{w,a})=
\XC_{\EC,\ph}(v)\g_{\Lb,v,\z}^\Gb(w)\z(\aha).$$
On the other hand, if $l \not\in \Cb_w^F$, then
$$\XC_{\EC_w,\ph_w}(l)=0.$$
\end{coro}

\bi

\remark{importance} In the theory of character sheaves applied to finite 
reductive groups, the characteristic functions $\XC_{\EC_w,\ph_w}$ play a 
crucial role, as it is shown in Theorem \ref{formule lusztig}. 
The Corollary \ref{description chi} shows the importance 
of the determining of the linear character $\g_{\Lb,v,\z}^\Gb$. 
We will show in Part II 
how the knowledge of the linear character $\g_{\Lb,v,\z}^\Gb$ 
whenever $v$ is regular and $p$ is good for $\Lb$ leads to an 
improvement of Digne, Lehrer and Michel's theorem on Lusztig 
restriction of Gel'fand-Graev characters \cite[Theorem 3.7]{DLM2}.\finl

\bi

\remark{scalaire} The characteristic function $\XC_{\EC_w,\ph_w}$ depends 
on the choice of the isomorphism $\ph$. Since $\EC$ is an irreducible 
local system, two \isos between $F^*\EC$ and $\EC$ differ only by 
a scalar. Hence the two characteristic functions they define 
differ also by the same constant. This shows that the formula 
in Corollary \ref{description chi} cannot be improved, because 
the factor $\XC_{\EC,\ph}(v)$ depends on the choice of the 
\iso $\ph$. We can give it, by multiplying $\ph$ by a scalar, 
any non-zero value we want.\finl

\bi

\bi

\newpage

%%%%%%%%%%%%%%%%%%
%%%%%%%%%%%%%%%%%%
%%%%%%%%%%%%%%%%%%
%%%%%%%%%%%%%%%%%%
%%%%%%%%%%%%%%%%%%
%%%%%%%%%%%%%%%%%%
%%%%%%%%%%%%%%%%%%
%%%%%%%%%%%%%%%%%%
%%%%%%%%%%%%%%%%%%
%%%%%%%%%%%%%%%%%%
%%%%%%%%%%%%%%%%%%
%%%%%%%%%%%%%%%%%%
%%%%%%%%%%%%%%%%%%
%%%%%%%%%%%%%%%%%%
%%%%%%%%%%%%%%%%%%
%%%%%%%%%%%%%%%%%%
%%%%%%%%%%%%%%%%%%

%\pagenumbering{arabic}

\begin{centerline}{\Large \bf Actions of relative Weyl groups II}\end{centerline}

\bi

\begin{centerline}{\sc C\'edric Bonnaf\'e}\end{centerline}

\bi

\begin{centerline}{\today}\end{centerline}

\bi

\bi

\begin{quotation}{\small \noindent {\sc Abstract - } 
The aim of this second part is to compute explicitly the Lusztig restriction 
of the characteristic function of a regular unipotent class of a 
finite reductive group, improving slightly a theorem of Digne, Lehrer and Michel. 
We follow their proof but we introduced a new information, namely the 
computation of some morphisms $\ph_{\Lb,v}^\Gb$ defined in the first 
part when $v$ is a regular unipotent element. This new information 
is obtained by studying generalizations of the variety of companion 
matrices which were introduced by Steinberg.}\end{quotation}

\bi

\begin{quotation}
{\small \noindent {\bf MSC Classification.} 20G40, 14L99}
\end{quotation}

\bi

\section*{Introduction to this part}~

\med

\noindent{\bf Hypothesis and notation :} {\it We work under the hypothesis 
of the beginning of the first part. More precisely, we keep all the notation 
and hypothesis introduced until the end of 
Section \ref{ele sec} only. For instance, unless otherwise specified, we 
assume neither that $\Cb$ supports a cuspidal local system, nor that $\FM$ 
is an algebraic closure of a finite field.
However, the following assumption will be added and will remain in force all 
along this second part~: we suppose that $v$ is a regular unipotent element of $\Lb$.}

\bi 

The aim of this second part is to compute explicitly the Lusztig restriction 
of the characteristic function of a regular unipotent class of a 
finite reductive group. Digne, Lehrer and Michel proved that 
this is again the characteristic function of a regular unipotent class
\cite[Theorem 3.7]{DLM2}, but they were not able to determine which one. 
In this paper, we follow their strategy (using extensively the theory 
of character sheaves) and introduce the new ingredients developed 
in the first part. The key point here is that we are able to compute 
explicitly the morphism $\ph_{\Lb,v}^\Gb$  because $v$ is regular 
(see Table \ref{tabletable}) under some restriction on $\Lb$. 
This new information allows us 
to obtain the desired precision on Digne-Lehrer-Michel's theorem
(see Theorem \ref{but dans la vie}). 

The computation of the morphism $\ph_{\Lb,v}^\Gb$ whenever $v$ is regular 
in $\Lb$, $p$ is good for $\Gb$ and $\Lb$ is cuspidal 
(see \SEC\ref{cuspidal sous} or \cite[\SEC 1]{Bonnafe} 
for the definition of cuspidality) constitutes the hardest part 
of this paper. 
For this purpose, we introduce the variety 
$\Xb_\cdo=\Xb \cap \Ub\cdo \cap \cdo \Ub^-$ (the variety $\Ub\cdo \cap \cdo\Ub^-$ 
was studied 
by Steinberg \cite{steinreg} in connection with regular elements) 
and its inverse image $\Xbt_\cdo$ and $\Xbt_\cdo^{\pr}$ under 
the morphisms $\piba$ and $\piba'$. Here, $\cdo$ is a representative 
of a standard Coxeter element of $W$, and $\Ub^-$ is the unipotent 
radical of the Borel subgroup $\Bb^-$ opposed to $\Bb$ relatively to 
$\Tb$. It turns out that $\Xb_\cdo$ is contained 
in the variety $\Xb_\mini$, that $C^\Gb \cap \Xb_\cdo$ consists in a single element, 
and that $\Yb \cap \Xb_\cdo$ 
is a dense open subset of $\Xb_\cdo$. 
This shows that the computation of the stabilizer $H_\Gb(\Lb,v,1)$ 
may be performed in the variety $\Xbt_\cdo^\pr$. 

On the other hand, using and improving slightly 
Steinberg's work, we can give a complete 
description of the varieties $\Xbt_\cdo$ and $\Xbt_\cdo^\pr$ 
in terms of the root datum of $\Gb$. However, it is not clear how 
$W_\Gb(\Lb,v)$ acts on $\Xbt_\cdo^\pr$ 
through this description (the action of its subgroup 
$A_\Lb(v)$ is well-understood). But, using some reduction arguments 
coming from Part I, the computation may be reduced to five cases~: 
one lives inside the special linear groups and the four others 
live in classical groups of small ranks (less than or equal to $5$). 
In the special linear group, the computation is done by brute force 
(see \SEC\ref{cas A}). 
In the remaining cases, another reduction argument 
shows that the result can be deduced from a very small information, 
namely the number of irreducible components of the variety $\Xbt_\cdo^\prime$ 
(see \SEC\ref{cas B}.A). 
This variety is described as a subvariety of $\Tb$ given by very 
explicit equations (see Theorem \ref{iso eta}). 
So the number of irreducible components can be easily 
obtained by a case-by-case analysis (see \SEC\ref{cas B}.B, C and D). 

\medskip

This part is organized as follows. In Section \ref{preliminaire}, we recall some 
well-known preliminary results. 
In Sections \ref{companion section} and \ref{sec jor com}, we introduce the companion 
variety and give the description of the varieties $\Xbt_\cdo$ and $\Xbt_\cdo^\prime$ 
defined above. In Section \ref{phi sec} we give a table providing 
the values of $\ph_{\Lb,v}^\Gb$ whenever $\Gb$ is semisimple, simply 
connected and quasi-simple and $\Lb$ is cuspidal (the general case may be deduced 
easily from this table using results of Part I). We also 
give in the same section a useful compatibility property of the 
family of morphisms $\ph_{\Lb,v}^\Gb$ (see Corollary 
\ref{compatible prop}), and we explain how 
the results of this table may be reduced to five cases. 
Sections \ref{cas A} and \ref{cas B} 
are devoted to the proof of the result for these five cases. Section 
\ref{cas A} deals with the case of the special linear group while 
Section \ref{cas B} deals with the four remaining ``small'' cases.
The last section \SEC\ref{section fini} concerns the application 
to finite reductive groups (Lusztig restriction of characteristic functions 
of regular unipotent classes) of all the material developed 
in these two papers. This was the original motivation of this work. 
It must be noticed that the final result (see Theorem \ref{but dans la vie}) is 
in contradiction with \cite[Conjecture 2]{Bonnafe}~!

\bi

\bi

\section*{Notation}~

\med

\bi

For the purpose of this part, we need some further notation regarding 
the root system and the Lie algebra of $\Gb$. First, let $r=|\D|$ 
and let $r'=|\D_\Lb|$. Note that $r$ and $r'$ are the semisimple 
ranks of $\Gb$ and $\Lb$ respectively. We set $\D=\{\a_1,\dots,\a_r\}$. 
Let $w_0$ denote the longest element of $W$ (with respect to the 
length function on $W$ defined by $\D$) and let $c$ denote the Coxeter element 
$s_{\a_1} \dots s_{\a_r}$ of $W$. 
Let $\Bb^-$ be the \borel of $\Gb$ opposite to $\Bb$ with respect to $\Tb$ and 
let $\Ub^-$ be its unipotent radical. 
We denote by $\gGB$, $\tGB$, $\bGB$, $\bGB^-$, $\uGB$, and $\uGB^-$ the Lie algebras 
of $\Gb$, $\Tb$, $\Bb$, $\Bb^-$, $\Ub$, and $\Ub^-$ respectively. 
For each $\a \in \Phi$, 
we set $e_\a =(dx_\a)_1(1) \in \gGB$ so that $e_\a$ is a non-zero element of the 
$1$-dimensional Lie algebra $\uGB_\a$ of $\Ub_\a$. 
Let $\Phi^\vee$ denote the coroot system of $\Gb$ relative to $\Tb$, 
let $\Phi \to \Phi^\vee$, $\a \mapsto \a^\vee$ denote 
the canonical bijection, and let $\D^\vee=\{\a^\vee~|~\a \in \D\}$.

Since $-w_0$ permutes the elements of $\D$, there exists a unique sequence 
$1 \le i_1 < \dots < i_{r'} \le r$ such that $\D_\Lb=\{-w_0(\a_{i_1}),\dots,
-w_0(\a_{i_{r'}})\}$. With this notation, we can define
$$\Ub_*=(\Ub_{-w_0(\a_1)}-\{1\}) \dots (\Ub_{-w_0(\a_r)}-\{1\}) \incl \Ub$$
$$\Ub_{\Lb,*}=(\Ub_{-w_0(\a_{i_1})}-\{1\}) \dots (\Ub_{-w_0(\a_{i_{r'}})}-\{1\})
=\pi_\Lb(\Ub_*).
\leqno{\mathrm{and}}$$
Then $\Ub_*$ and $\Ub_{\Lb,*}$ are locally closed subvarieties of $\Gb$.

Finally, we denote by $\ZC(\Gb)$ the group of components of the center 
of $\Gb$, that is $\ZC(\Gb)=\Zb(\Gb)/\Zb(\Gb)^\ci$. 

\bi 

\sec{Preliminaries on regular elements\label{preliminaire}}~

\med

This section is a recollection of well-known facts about regular 
elements in reductive groups. In \SEC\ref{first sous}, we gather 
some consequences of the regularity of $v$. In \SEC\ref{sub reg B}, we 
provide some conditions for regular elements lying in $\Bb$ subgroup 
to be conjugated in $\Gb$. 

\bi

\sub{Regularity of $v$\label{first sous}} 
First, since there is only one regular unipotent 
class in $\Lb$, we have
$$W_\Gb(\Lb,\Sigb)=W_\Gb(\Lb).$$
Moreover, $\overline{\Cb}=\Lb_\uni$. Therefore, 
$$\overline{\Sigb}.\Vb=\Zb(\Lb)^\ci.\Pb_\uni.$$
This shows that $u$ is a regular unipotent element of $\Gb$. 
Without loss of generality, we may, and we will, assume that
$$u=x_{-w_0(\a_1)}(1) \dots x_{-w_0(\a_s)}(1) \in \Ub_*$$
$$v=x_{-w_0(\a_{i_1})}(1) \dots x_{-w_0(\a_{i_s})}(1) \in \Ub_{\Lb,*}.
\leqno{\mathrm{and}}$$
Note that $v=\pi_\Lb(u)$. 
We recall the following well-known results about the structure of the centralizer 
of $u$ in $\Gb$~:

\bi

\begin{prop}\label{abelianite}
${\mathrm{(a)}}$ $C_\Gb(u)=\Zb(\Gb) \times C_\Ub(u)$.

\tete{b} If $p$ is good for $\Gb$, then $A_\Ub(u)=\{1\}$.

\tete{c} If $p$ is bad for $\Gb$ and if $\Gb$ is quasi-simple, 
then $u \not\in C_\Ub^\ci(u)$ and 
$A_\Ub(u)$ is cyclic and generated by the image of $u$.

\tete{d} $C_\Gb(u)$ is abelian, therefore 
$A_\Gb(u)=\ZC(\Gb) \times A_\Ub(u)$ is also abelian.
\end{prop}

\bi

The following corollary shows that many results of Part I are 
valid in this situation.

\bi

\begin{coro}\label{CGU}
$C_\Gb(u) \incl \Pb$.
\end{coro}

\bi

%\remark{cgu p} In \cite{Bonnafe I..}, many results were subject to the condition that 
%$C_\Gb(u) \incl \Pb$. In this paper, this assumption is satisfied by the 
%previous corollary.\finl
%
%\bi
%
%\remark{clv l} 

We now provide a description of the variety $\Xb$ using the 
Steinberg map $\na$.

\bi

\begin{prop}\label{X et Y}
$\Xb=\na^{-1}\bigl(\na(\Zb(\Lb)^\ci)\bigr)$.
\end{prop}

\bi

\proof By definition, 
$$\Xb=\bigcup_{g \in \Gb} \lexp{g}{(\Zb(\Lb)^\ci\Pb_\uni)}
=\bigcup_{g \in \Gb} \lexp{g}{(\Zb(\Lb)^\ci\Lb_\uni\Vb)}.$$
It follows from Lemma \ref{V conjugue} that 
$\Xb$ is contained in $\na^{-1}\bigl(\na(\Zb(\Lb)^\ci)\bigr)$. 

Conversely, let $y \in 
\na^{-1}\bigl(\na(\Zb(\Lb)^\ci)\bigr)$. Let $z$ denote its semisimple part. 
We may assume 
that $z \in \Zb(\Lb)^\ci$. 
Now let $\Gb'=C_\Gb^\ci(z)$ and $\Pb'=C_\Pb^\ci(z)$. Then, by 
\cite[Corollary 11.12]{borel}, $y \in \Gb'$. Moreover, by 
\cite[Proposition 1.11]{dmgnc}, $\Pb'$ is a \para of $\Gb'$. 
Hence there exists $g \in \Gb'$ such 
that $\lexp{g}{y} \in \Pb'$. But, since $g$ centralizes $z$, the semisimple part 
of $\lexp{g}{y}$ is equal to $z$. Consequently, $\lexp{g}{y} \in z \Pb^\pr_\uni 
\incl \Zb(\Lb)^\ci \Pb_\uni$. Hence $y \in \Xb$.\fin

\bi

\begin{prop}\label{pfff}
$\Xb_\mini$ is the set of elements of $\Xb$ which are regular 
in $\Gb$.
\end{prop}

\proof Indeed, $\dim C_\Lb(v)=\dim \Tb$ because $v$ is regular in $\Lb$.\fin

\bi

\sub{Regular elements in $\Bb$\label{sub reg B}} If two elements of the Borel subgroup 
$\Bb$ are regular in $\Gb$ and conjugate under some element $g \in \Gb$, 
then it is possible, in many cases, to conclude that $g \in \Bb$. The next 
lemma provides three different conditions under which this result holds~:

\bi

\begin{lem}\label{conjugue regulier}
Let $b$ and $b'$ be two $\Gb$-regular elements lying in $\Bb$ and let $g \in \Gb$. 
Assume that $\pi_\Tb(b)=\pi_\Tb(b')$ and that 
$b'=\lexp{g}{b}$. Assume that at least one of the following conditions holds~:

\tete{1} $b$ is unipotent.

\tete{2} $g$ is unipotent.

\tete{3} The derived group of $\Gb$ is simply connected.

\noindent Then $g \in \Bb$.
\end{lem}

\bi

\proof Let $s$ (\resp $s'$) be the semisimple part of $b$ (\resp $b'$). 
By Lemma \ref{V conjugue} applied to the case where $\Lb=\Tb$ and $\Pb=\Bb$, 
$s$ and $s'$ are conjugate under $\Ub$. Hence we may assume that $s=s'$. 
In particular, $g \in C_\Gb(s)$. 

Let $x$ (\resp $x'$) denote the unipotent part of $b$ (\resp $b'$). Then $x$ and 
$x'$ are regular unipotent elements of $C_\Gb^\ci(s)$ and they both belong to 
$C_\Bb(s)$~: this group is connected \cite[Theorem 10.6 (5)]{borel}  
and is a Borel subgroup of $C_\Gb^\ci(s)$ 
\cite[Theorem 1.8 (ii) and (iii)]{dmgnc}. Consequently, 
$C_\Bb(s)$ is the unique \borel of $C_\Gb^\ci(s)$ containing 
$x$ (\resp $x'$). This implies 
that $g \in N_{C_\Gb(s)}(C_\Bb(s))$ (indeed, $x'=\lexp{g}{x}$).

Since $N_{C_\Gb^\ci(s)}(C_\Bb(s))=C_\Bb(s)$, 
the conclusion follows from the fact that, if either (1), (2) or (3) holds, 
then $g \in C_\Gb^\ci(s)$. Indeed, in case (1), we have $s=1$ so $C_\Gb(s)=\Gb$ 
is connected.
In case (2), $g \in C_\Gb^\ci(s)$ by \cite[Corollary 11.12]{borel}. In case (3), 
since the derived group of $\Gb$ is simply connected, 
the group $C_\Gb(s)$ is connected \cite[Corollary 8.5]{steinendo}.\fin

\bi

\rem In general, the conclusion of Lemma \ref{conjugue regulier} does not hold. 
For instance, if $\Gb=\Pb\Gb\Lb_2(\FM)$, if $\Bb$ consists of upper triangular 
matrices and $\Tb$ of diagonal matrices, and if $p \not=2$, 
then the conclusion fails for the following 
elements~:
$$b=b'=\matricecrochet{1 & 0 \\ 0 & -1}\quad\quad\quad{\mathrm{and}}\quad\quad\quad
g=\matricecrochet{0 & 1 \\ 1 & 0}.$$
Here $\matricecrochet{a & b \\ c & d}$ denotes the image in $\Pb\Gb\Lb_2(\FM)$ 
of the invertible $2 \times 2$ matrix $\matrice{a & b \\ c & d}$ 
under the canonical 
\mor $\Gb\Lb_2(\FM) \to \Pb\Gb\Lb_2(\FM)$.\finl

\bi

\sec{Companion variety\label{companion section}}~

\bi

\sub{The Bruhat cell associated to a Coxeter 
element\label{bruhat}}~
The Bruhat cell $\Bb c \Bb$ is a central object in Steinberg's paper 
\cite{steinreg} (recall that $c=s_{\a_1}\dots s_{\a_r}$ 
is a standard Coxeter element in $W$). Here are two results taken 
from \cite{steinreg} which are of fundamental importance for 
our purpose.

\bi

\begin{theo}[{\bf Steinberg}]\label{rappel}
${\mathrm{(a)}}$ All elements of $\Bb c \Bb$ are regular \cite[Remark 8.8]{steinreg}.

\tete{b} If $\Gb$ is semisimple and simply connected, then 
the restriction $\na_\cdo : \Ub\cdo \cap \cdo\Ub^- \to \Tb/W$ of 
$\na$ is an \iso of varieties \cite[Theorem 1.4]{steinreg}.
\end{theo}

\bi

\rem Whenever $\Gb=\Sb\Lb_{r+1}(\FM)$, it is possible 
to choose $\Ub$ and $\cdo$ in such a way that $\Ub\cdo\cap\cdo\Ub^-$ 
is the set of {\it companion matrices}. In this case, 
Theorem \ref{rappel} (b) is well-known. So, in an arbitrary 
reductive group $\Gb$, the closed subvariety $\Ub\cdo\cap\cdo\Ub^-$ 
of $\Gb$ may be viewed as a generalization of the variety of companion 
matrices. We call it the {\it companion variety} of $\Gb$.\finl

\bi

In the next proposition, we extend slightly \cite[Proposition 8.9]{steinreg} 
and, in the 
rest of this section, we provide some immediate consequences of this result. 

\bi

\begin{prop}\label{isomorphisme}
The map
$$\fonction{\n_0}{(\Ub \cap \lexp{c}{\Ub}) \times (\Ub \cap \lexp{c}{\Ub^-})}{\Ub}{(u,
v)}{\dot{c}^{-1}u\dot{c}vu^{-1}}$$
is an isomorphism of varieties.
\end{prop}

\bi

\proof By \cite[Proposition 8.9]{steinreg}, the map $\n_0$ 
is bijective. Moreover, the varieties $\Ub$ and 
$(\Ub \cap \lexp{c}{\Ub}) \times (\Ub \cap \lexp{c}{\Ub^-})$ 
are smooth, hence normal. 
Therefore, by \cite[Theorems AG.17.3 and AG.18.2]{borel}, 
we only need to prove that the differential $(d\n_0)_{(1,1)}$ is surjective. 

The linear map $(d\n_0)_{(1,1)}$ can be described as follows~:
$$\fonction{(d\n_0)_{(1,1)}}{\bigl(\uGB \cap (\Ad \dot{c})(\uGB) \bigr) 
\oplus \bigl(\uGB \cap (\Ad \dot{c})(\uGB^-) \bigr)}{\uGB}{u \oplus 
v}{\Ad \dot{c}^{-1}(u) + v - u.}$$
Since the dimension of both $\FM$-vector spaces is equal to $\dim \Ub$, we only need 
to prove that $(d\n_0)_{(1,1)}$ is injective. So let $u \oplus v \in 
\Ker (d\n_0)_{(1,1)}$. Write 
$$u=\sum_{\a \in \Phi_c^+} a_\a e_\a\quad\quad{\mathrm{and}}\quad\quad 
v=\sum_{\a \in \Phi_c^-} b_\a e_\a$$
where $\Phi_c^\pm=\Phi^+ \cap c(\Phi^\pm)$. If $\a$ is any root, there exists 
$c_\a \in\FM^\times$ such that $(\Ad \cdo)^{-1}(e_\a)=c_\a e_{c^{-1}(\a)}$. We have 
$(\Ad \cdo)^{-1}(u)+v-u = 0$ therefore,
$$\sum_{\a \in \Phi_c^+} a_\a c_\a e_{c^{-1}(\a)} + \sum_{\a \in \Phi_c^-} b_\a e_\a -
\sum_{\a \in \Phi_c^+} a_\a e_\a = 0.\leqno{(\sharp)}$$

Let $\a \in \Phi_c^-$. We denote 
by $k_\a$ the smallest natural number such that $c^{k_\a}(\a) \in \Phi^-$ and we 
define
$$A_\a=\{\a,c(\a),\dots,c^{k_\a-1}(\a)\}$$
$$B_\a=\{c(\a),\dots,c^{k_\a-1}(\a)\}.\leqno{\mathrm{and}}$$
Note that $B_\a=\vide$ if $k_\a=1$. 
Then, by \cite[chapter VI, ${\mathrm{n}}^\ci$ 1.11, Proposition 33 (iv)]{bourbaki}, 
we have
$$\Phi^+=\coprod_{\a \in \Phi_c^-} A_\a.$$
It implies that
$$\Phi_c^+=\coprod_{\a \in \Phi_c^-} B_\a.$$
Let $\a \in \Phi_c^-$. If $k_\a=1$, then $b_\a=0$ by $(\sharp)$. If $k_\a \ge 2$, 
then it follows from $(\sharp)$ that 
\eqna
b_\a&=&-a_{c(\a)}c_{c(\a)}\\
a_{c(\a)}&=& a_{c^2(\a)}c_{c^2(\a)} \\
&\dots& \\
a_{c^{k_\a-1}(\a)}&=&0.
\endeqna
From this we deduce that $a_{\a'}=0$ for all $\a' \in B_\a$ and that $b_\a=0$. 
Finally, we get $u=v=0$.\fin

\bi

In view of our applications, the most important result of this subsection 
is given by the following corollary, which is a straightforward consequence 
of Proposition \ref{isomorphisme}.

\bi

\begin{coro}\label{gamma}
The map
$$\fonction{\n}{\Ub \times (\Ub \cdo \cap \cdo \Ub^-)}{\Ub \cdo \Ub}{(u,g)}{ugu^{-1}}$$
is an isomorphism of varieties.
\end{coro}

\bi

\proof By Proposition \ref{isomorphisme}, the map
$$\fonction{\n_1}{(\Ub \cap \lexp{c^{-1}}{\Ub}) \times (\Ub \cdo \cap \cdo
\Ub^-)}{\Ub\cdo}{(u,g)}{ugu^{-1}}$$
is an isomorphism of varieties. By \cite[Theorem 14.12 (a)]{borel}, the map
$$\fonction{\n_2}{\Ub\cdo \times (\Ub \cap \lexp{c^{-1}}{\Ub^-})}{\Ub\cdo 
\Ub}{(x,v)}{x v^{-1}}$$
is an isomorphism of varieties. Hence, the map 
$$\fonctio{(\Ub \cap \lexp{c^{-1}}{\Ub^-})\times 
(\Ub \cap \lexp{c^{-1}}{\Ub}) \times (\Ub \cdo \cap 
\cdo\Ub^-)}{\Ub\cdo\Ub}{(v,u,g)}{\n_2(v\n_1(u,g),v)=vugu^{-1}v^{-1}}$$
is an isomorphism of varieties. The result follows then from the fact that 
the map
$$\fonctio{(\Ub \cap \lexp{c^{-1}}{\Ub^-}) \times (\Ub \cap 
\lexp{c^{-1}}{\Ub})}{\Ub}{(v,u)}{vu}$$
is also an isomorphism of varieties \cite[Remark following Theorem 14.4]{borel}.\fin

\bi

\begin{coro}
The maps 
$$\fonctio{(\Ub \cap \lexp{c}{\Ub}) \times (\Bb \cap \lexp{c}{\Bb^-})}{\Bb}{(u,
b)}{\dot{c}^{-1}u\dot{c}bu^{-1}}$$
$$\fonctio{\Ub \times (\Bb c \cap c 
\Bb^-)}{\Bb c \Bb}{(u,g)}{ugu^{-1}}\leqno{\mathit{and}}$$
are isomorphisms of varieties.
\end{coro}

\bi

\sub{A technical lemma\label{techno}} By \cite[Lemma 4.5]{steinreg}, 
$\lexp{w_0}{(\Tb\Ub_*)}$ is contained in $\Bb c \Bb$. In particular, all its elements 
are regular. Let $\o_0 : 
\lexp{w_0}{(\Tb\Ub_*)} \to \Tb$ denote the morphism of varieties that sends 
$\lexp{w_0}{(tu)}$ to $t$ for all $t \in \Tb$ and $u \in \Ub_*$. 
Finally, let 
$$\o : \lexp{w_0}{(\Tb\Ub_*)} \cap \Ub\cdo\Ub \longto \Tb$$ 
denote the restriction of $\o_0$. 

\bi

\begin{lem}\label{omega iso}
If $\Gb$ is semisimple and simply connected, then $\o$ is an isomorphism.
\end{lem}

\bi

\proof First, let us prove the result 
whenever $\Gb=\Sb\Lb_2(\FM)$ and whenever we have made 
the following choices~:
$$\fonction{x_{\a}}{\FM}{\Gb}{a}{\matrice{1 & a \\ 0 & 1},}$$
$$\fonction{x_{-\a}}{\FM}{\Gb}{a}{\matrice{1 & 0 \\ a & 1},}$$
$$\fonction{\a^\ve}{\FM^\times}{\Gb}{a}{\matrice{a & 0 \\ 0 & a^{-1}}},$$
$$\sdo_{\a}=\matrice{0 & -1 \\ 1 & 0}.\leqno{\mathrm{and}}$$
Note that $w_0=c=s_\a$. Here, $\a=\a_1$ denotes the unique positive root. 
Then we have, for all $a \in \FM^\times$, 
$$x_{\a}(a)\sdo_{\a}x_{\a}(a^{-1})=\a^\ve(a)x_{-\a}(a).$$
Moreover, $x_{\a}(a)\sdo_\a x_{\a}(a') \in \lexp{w_0}{(\Tb\Ub_*)}$ if and only if 
$aa'=1$. This proves Lemma \ref{omega iso} in this particular 
case.

\med

Now, we prove the Lemma \ref{omega iso} in the general case. First, note 
that if the result holds for some representative $\cdo$ of $c$, then 
it holds for any. We will prove Lemma \ref{omega iso} by using 
the following choice. For each $1 \le i \le r$, let $\Gb_i$ 
denote the subgroup of $\Gb$ generated by $\Ub_{\a_i}$ and $\Ub_{-\a_i}$. Then 
$\Gb_i\simeq \Sb\Lb_2(\FM)$ since $\Gb$ is semisimple and simply connected. We 
choose a representative $\sba_{\a_i}$ of $s_{\a_i}$ in $\Gb_i$. 
To prove (a), we may (and we will) assume that $\cdo=\sba_{\a_1} 
\dots \sba_{\a_r}$. By the previous 
computation in $\Sb\Lb_2(\FM)$, there exists a morphism of varieties $f_i : 
\FM^\times \to \Ub_{-\a_i} - \{1\}$ such that
$$\forall a \in \FM^\times,~ \a_i^\ve(a)f_i(a) \in \Ub_{\a_i}\sba_{\a_i}\Ub_{\a_i}.$$
Now let
$$\fonction{f}{(\FM^\times)^r}{\lexp{w_0}{(\Tb\Ub_*)}}{(a_1,\dots,a_r)}{
\a_1^\ve(a_1)f_1(a_1)\dots\a_r^\ve(a_r)f_r(a_r).}$$
Since $\lexp{w_0}{(\Tb\Ub_*)}=\Tb(\Ub_{-\a_1}-\{1\})\dots 
(\Ub_{-\a_r}-\{1\})$, the image of $f$ is indeed contained in 
$\lexp{w_0}{(\Tb\Ub_*)}$. Moreover, 
$$f(a_1,\dots,a_r) \in 
(\Ub_{\a_1}\sba_{\a_1}\Ub_{\a_1})(\Ub_{\a_2}\sba_{\a_2}\Ub_{\a_2})
\dots (\Ub_{\a_r}\sba_{\a_r}\Ub_{\a_r}) \incl \Ub\cdo\Ub,$$
$$\o(f(a_1,\dots,a_r))= \lexp{w_0}{(\a_1^\ve(a_1)\dots\a_r^\ve(a_r))}
\leqno{\mathrm{and}}$$
for any $(a_1,\dots,a_r) \in (\FM^\times)^r$. Since $\Gb$ is semisimple and 
simply connected, the morphism of algebraic groups
$$\fonction{\alpb}{(\FM^\times)^r}{\Tb}{(a_1,\dots,
a_r)}{\lexp{w_0}{(\a_1^\ve(a_1)\dots\a_r^\ve(a_r))}}$$
is an isomorphism of algebraic groups. Hence, $\fb=f \ci \alpb^{-1} : \Tb \to 
\lexp{w_0}{(\Tb\Ub_*)} \cap \Ub \cdo \Ub$ is a morphism of varieties and 
$$\o \ci \fb=\Id_\Tb.$$
To complete the proof of (a), we now only need to show that 
$\o$ is injective. Let $g$ and $h$ be two elements of $\lexp{w_0}{(\Tb\Ub_*)} 
\cap \Ub \cdo \Ub$ such that $\o(g)=\o(h)$. Then, by Lemma \ref{V conjugue} 
applied to the case where $\Lb=\Tb$ and $\Pb=\Bb$, we have $\na(g)=\na(h)$. 
Since $\na_\cdo : \Ub\cdo \cap \cdo \Ub^- \to \Tb/W$ is an isomorphism (\cf 
Theorem \ref{rappel} (b)), $g$ and 
$h$ are conjugate under $\Ub$ (\cf Corollary \ref{gamma}). In other words, 
$\lexp{\wdo_0}{g}$ and $\lexp{\wdo_0}{h}$ are conjugate under an element 
$x \in \Ub^-$. But 
$\lexp{\wdo_0}{g}$ and $\lexp{\wdo_0}{h}$ are regular elements of $\Bb$ (indeed, 
$g$ and $h$ belong to $\Ub\cdo \Ub \incl \Bb c \Bb$ so they are regular by  
Theorem \ref{rappel} (a)). Hence,  
by Lemma \ref{conjugue regulier} (2) or (3), $x \in \Bb$. Therefore $x=1$ and 
the injectivity of $\o$ follows.\fin

\bi

For further investigation, we will need an explicit formula for the 
inverse of the morphism $\o$. This formula can be obtained from the proof 
of Lemma \ref{omega iso} as follows. 
By definition, the group $\Gb$ is semisimple and simply connected if and only if 
$(\a^\ve)_{\a \in \D}$ is a basis of $Y(\Tb)$. We denote 
by $(\varpi_\a)_{\a \in \D}$ the basis of $X(\Tb)$ dual to 
$(\a^\ve)_{\a \in \D}$. 
Then, for a suitable choice of the family $(\sdo_\a,x_\a,x_{-\a})_{\a \in \D}$ 
and of $\wdo_0$, and if $\cdo=\sdo_{\a_1} \dots \sdo_{\a_r}$, we have
$$\lexp{\wdo_0^{-1}}{(\o^{-1}(t))}=
\a_1^\ve\bigl(\varpi_{\a_1}(t)\bigr)x_{-w_0(\a_1)}\bigl(\varpi_{\a_1}(t)\bigr)
\dots  \a_r^\ve\bigl(\varpi_{\a_r}(t)\bigr)x_{-w_0(\a_r)}\bigl(\varpi_{\a_r}(t)\bigr)$$
for any $t \in \Tb$.
(This follows from the proof of Lemma \ref{omega iso}.)
We set, for each $1 \le i \le r$, 
\equat\label{defi chi}
\ch_i=\varpi_{\a_i}-\sum_{j=i+1}^r \langle -w_0(\a_i),\a_j^\ve \rangle 
\varpi_{\a_j}.
\endequat
Then $(\ch_i)_{1 \le i \le r}$ is a basis of $X(\Tb)$ and, from 
the previous formula, we get~:

\begin{coro}\label{inverse omega}
Assume that $\Gb$ is semisimple and simply connected. 
Then, for a suitable choice of the family $(\sdo_\a,x_\a,x_{-\a})_{\a \in \D}$ 
and of $\wdo_0$, and if $\cdo=\sdo_{\a_1} \dots \sdo_{\a_r}$, we have
$$\lexp{\wdo_0^{-1}}{(\o^{-1}(t))}=t 
x_{-w_0(\a_1)}\bigl(\ch_1(t)\bigr)\dots x_{-w_0(\a_r)}\bigl(\ch_r(t)\bigr)$$
for any $t \in \Tb$.
\end{coro}

\bi

\sec{Companion variety and Jordan decomposition\label{sec jor com}}~

\med

This section is concerned with the study of some varieties 
related to the companion variety. In order to use the 
strongest results of \cite{steinreg} and of the previous section, 
we need however to make the following hypothesis.

\bi

\noindent{\bf Hypothesis :} {\it Until the end of \SEC\ref{sec jor com}, 
we assume that $\Gb$ is semisimple and simply connected.}

\bi

For the convenience of the reader, we will recall, in each result stated 
until the end of \SEC\ref{sec jor com}, that $\Gb$ is assumed to be semisimple 
and simply connected. 

\bi

\sub{Companion variety and parabolic subgroups\label{compa para subsub}} 
Let $\l : \Tb \to \Ub$ and $\m : \Tb \to \Ub\cdo \cap \cdo \Ub^-$ 
denote the morphisms of varieties defined, for any $t \in \Tb$, by
$$\n^{-1}(\o^{-1}(t))=(\l(t),\m(t))$$
or, equivalently, by 
$$\l(t)\m(t)\l(t)^{-1}=\o^{-1}(t).$$

\bi

\remark{nabla nabla} By Lemma \ref{V conjugue}, the semisimple part of $\o^{-1}(t)$ 
is conjugate to $t$. Therefore, $\m(t)=\na_\cdo^{-1}(\na(t))$, because 
$\na_\cdo : \Ub\cdo \cap \cdo \Ub^- \to \Tb/W$ is an \iso of varieties 
(\cf Theorem \ref{rappel} (b)).\finl

\bi

We consider here the following varieties 
\eqna
\Xb_\cdo &=& \Xb \cap (\Ub\cdo \cap \cdo \Ub^-) \\
&=&\na^{-1}_\cdo\bigl(\na(\Zb(\Lb)^\ci)\bigr) 
\endeqna
$$\Xbt_\cdo = \piba^{-1}(\Xb_\cdo).\leqno{\mathrm{and}}$$ 
The equality $\Xb \cap (\Ub\cdo \cap \cdo \Ub^-) =
\na^{-1}_\cdo\bigl(\na(\Zb(\Lb)^\ci)\bigr)$ comes from Proposition \ref{X et Y} 
and Theorem \ref{rappel} (b).
Then, $\Xb_\cdo$ is a closed subvariety of $\Xb$ isomorphic to 
$\na(\Zb(\Lb)^\ci)$, and $\Xbt_\cdo$ is a closed 
subvariety of $\Xbt$. Since all the elements of $\Xb_\cdo$ are regular 
(\cf Theorem \ref{rappel} (a)), we have
$$\Xb_\cdo \incl \Xb_\mini\quad\quad\quad{\mathrm{and}} \quad\quad\quad 
\Xbt_\cdo \incl \Xbt_\mini.$$
We denote by $\pi_\cdo : \Xbt_\cdo \to \Xb_\cdo$ the restriction of $\piba$. 
The square
$$\diagram
\Xbt_\cdo \rrto \ddto_{\DS{\pi_\cdo}} && \Xbt_\mini \ddto^{\DS{\pi_\mini}} \\
&&\\
\Xb_\cdo \rrto && \Xb_\mini
\enddiagram$$
is cartesian, so $\pi_\cdo$ is a finite morphism.

We denote by $\nablat_\cdo : \Xbt_\cdo \to \Zb(\Lb)^\ci$ the restriction of $\nablat : 
\Xbt \to \Zb(\Lb)^\ci$. We also set
$$\fonction{\xbt_\cdo}{\Zb(\Lb)^\ci}{\Xbt_\cdo}{z}{\l(z)^{-1}\wdo_0
*_\Pb \wdo_0^{-1} \o^{-1}(z) \wdo_0.}$$
note that $\xbt_\cdo$ is well-defined and that it does not depend on a 
choice of a representative $\wdo_0$ of $w_0$. 

\bi

\begin{theo}\label{key}
Assume that $\Gb$ is semisimple and simply connected. 
Then, the maps $\nablat_\cdo$ and $\xbt_\cdo$ are $W_\Gb(\Lb)$-equivariant 
isomorphisms of varieties inverse of each other. Moreover, the 
following diagram 
$$\diagram
\Xbt_\cdo \rrto^{\DS{\nablat_\cdo}} 
\ddto_{\DS{\pi_\cdo}} && \Zb(\Lb)^\ci \ddto^{\DS{\na}} \\
&\\
\Xb_\cdo \rrto_{\DS{\na_\cdo}} && \na(\Zb(\Lb)^\ci)
\enddiagram$$
is commutative.
\end{theo}

\bi

\proof First, the commutativity of the diagram is clear. Also, 
the $W_\Gb(\Lb)$-equivariance of $\nablat_\cdo$ follows from the 
same property for $\nablat$. Therefore, we only need to show that
$\nablat_\cdo$ and $\xbt_\cdo$ are inverse of each other.
Moreover, 
one can easily check that $\nablat_\cdo \ci \xbt_\cdo=\Id_{\Zb(\Lb)^\ci}$. 
To prove that $\xbt_\cdo \ci \nablat_\cdo = \Id_{\Xbt_\cdo}$, 
we only need to prove that $\nablat_\cdo$ is injective. 

Let $g *_\Pb x$ and $g' *_\Pb x'$ in $\Xbt_\cdo$ 
be such that $\nablat_\cdo(g*_\Pb x)=\nablat_\cdo(g' *_\Pb x')=z$. We may, and we will, 
assume that $x$ and $x'$ belong to $\Bb$. 
By Corollary \ref{V conjugue}, 
the semisimple part of $x$ and $x'$ are conjugate. But, by Theorem 
\ref{rappel} (b), $\na_\cdo : \Ub\cdo \cap \cdo\Ub^- \to \Tb/W$ is an isomorphism. 
Hence $gxg^{-1}=g'x' g^{\pr -1}$. Now, let $h=g^{-1}g'$. Then $x=hx'h^{-1}$, 
and $\pi_\Tb(x)=\pi_\Tb(x')=z$. Moreover, $x$ and $x'$ are regular. 
It follows from Lemma \ref{conjugue regulier} 
(3) that $h \in \Bb$. In particular, $g*_\Pb x=gh *_\Pb h^{-1}xh=g' *_\Pb x'$. 
The injectivity of $\nablat_\cdo$ follows. The proof of Theorem \ref{key} 
is complete.\fin

\bi

\example{P=B} Let us assume here, and only here, that $\Pb=\Bb$, 
so that $\Lb=\Tb$. In this case, $\Xb_\cdo=\Ub\cdo \cap \cdo \Ub^-$ and 
$$\Xbt_\cdo=\{g *_\Bb x \in \Gb \times_\Bb \Bb~|~
gxg^{-1} \in \Ub\cdo \cap \cdo\Ub^- \}.$$
Theorem \ref{key} proves that, if $\Gb$ is semisimple and simply connected, 
then $\Xbt_\cdo$ is isomorphic to $\Tb$~; the \iso is given 
by $\Xbt_\cdo \to \Tb$, $g *_\Bb x \mapsto \pi_\Tb(x)$. Moreover, the 
diagram 
$$\diagram 
\Xbt_\cdo \ddto_{\DS{\pi_\cdo}} \rrto^{\DS{\nablat_\cdo}} && \Tb \ddto^{\DS{\na}} \\
&& \\
\Ub\cdo \cap \cdo\Ub^- \rrto_{\DS{\na_\cdo}} && \Tb/W 
\enddiagram$$
is commutative.\finl 

\bi

\remark{restriction perverse} Let us assume in this remark, and only in this remark, 
that $\Cb$ supports a cuspidal local system $\EC$, and let us use here 
the notation introduced in Sections \ref{sec endo}, \ref{section endo}, and 
\ref{section 1 gamma} ($\z$, $\AC$, $K$, $\KC$...). 

Let $K_\cdo$ denote the restriction of $K[-\dim \Xb]$ to $\Xb_\cdo$ 
and let $\KC_\cdo$ denote the restriction of $\KC$ to $\Yb_\cdo$. Since $\Xb_\cdo$ 
is contained in $\Xb_\mini$, $K_\cdo$ is the restriction of $K_\mini$ 
to $\Xb_\cdo$. In particular, it is a constructible sheaf. Moreover~:

\bi

{\it If $\Gb$ is semisimple and simply connected, then 
$K_\cdo = IC(\Xb_\cdo,\KC_\cdo)$.}

\bi

Indeed, by the proper base change theorem \cite[Chapter VI, Corollary 2.3]{Milne}, 
we have
$$K_\cdo=(\pi_\cdo)_* \Kti_\cdo,$$
where $\Kti_\cdo$ denotes the restriction of $\Kti[-\dim \Xb]$ to $\Xbt_\cdo$. 
By Theorem \ref{key}, the \mor $\pi_\cdo : \Xbt_\cdo \to \Xb_\cdo$ 
is a finite morphism and the variety $\Xbt_\cdo$ is smooth. Therefore, 
by \cite[Example following Definition 5.4.4]{luicc}, it is sufficient to 
prove that $\Kti_\cdo$ is a local system, that is a complex concentrated in 
degree $0$ whose $0$th term is a local system. But $\Xb_\cdo 
\incl \Xb_{\min}$, so $\Xbt_\cdo \incl \Xbt_{\min}$. So $\Kti_\cdo$ 
is the restriction of the local system $\Kti_\mini$, so it is a local system.\finl

\bi

\sub{Jordan decomposition\label{sec jordan}}~
In order to get explicit informations on the Jordan decomposition of companion 
matrices, we need to use Corollary \ref{inverse omega} which is a partial 
result in this direction. For this purpose, we need to make the 
following hypothesis (the reader may check that the general case may 
be deduced from this particular one by slight modifications)~:

\bi

\noindent{\bf Hypothesis :} {\it From now on, and until the end of 
\SEC\ref{companion section}, 
we assume that the family $(\sdo_\a,x_\a,x_{-\a})_{\a \in \D}$, 
that $\wdo_0$, and that $\cdo$ are chosen in such a way that Corollary 
\ref{inverse omega} holds. For instance, we have $\cdo=\sdo_{\a_1} \dots \sdo_{\a_r}$.}

\bi

Let $\Xbt_\cdo^\pr$ be the inverse image of $\Xbt_\cdo$ in $\Xbt'$ under 
the \mor $\fba$. It is also equal to the inverse image of $\Xb_\cdo$ under 
the \mor $\piba'$. Moreover, since $\Xb_\cdo \incl \Xb_\mini$, we have $\Xbt_\cdo^\pr 
\incl \Xbt_\mini^\pr$. Finally, let $\fti_\cdo : \Xbt_\cdo^\pr \to \Xbt_\cdo$ 
and $\pi_\cdo^\pr : \Xbt_\cdo^\pr \to \Xb_\cdo$ denote the restrictions 
of the morphisms $\fti_\mini$ and  $\pi_\mini^\pr$ respectively. 
The aim of this subsection is to give an explicit description 
of the variety $\Xbt_\cdo^\pr$. 

\bi

Let us fix here a subtorus $\Sb$ 
of $\Tb$ such that $\Tb=\Sb \times \Zb(\Lb)^\ci$. Let $Z_\Lb=\Zb(\Lb) \cap \Sb$. 
Then $Z_\Lb \simeq \ZC(\Lb)$. 
We set
$$\Xbt_\cdo^1=\{(z,t) \in \Zb(\Lb)^\ci \times \Sb~|~\forall 1 \le k \le r',~
\ch_{i_k}(z)\a_{i_k}(\lexp{w_0}{t})^{-1}=1\}.$$
We will construct here an $A_\Lb(v)$-equivariant \iso $\Xbt_\cdo^1 
\times A_{\Ub_\Lb}(v) \to \Xbt_\cdo^\pr $. Here the action of $A_\Lb(v) 
\simeq Z_\Lb \times A_{\Ub_\Lb}(v)$ is given by right translation 
by $Z_\Lb$ on $\Sb$ and by right translation of $A_{\Ub_\Lb}(v)$ on 
itself.

Let us introduce some further notation. 
By Corollary \ref{inverse omega}, we have, for every $(z,t) \in \Xbt_\cdo^1$,
\equat\label{tgb}
\lexp{t\wdo_0^{-1}}{(\o^{-1}(z))} \in zv\Vb.
\endequat
Now, let 
$$\fonction{\vb}{\Xbt_\cdo^1}{\Vb}{(z,t)}{v^{-1}z^{-1}
\lexp{t\wdo_0^{-1}}{(\o^{-1}(z))}.}$$
Then $\vb$ is a \mor of varieties, and, by construction, we have, for every 
$(z,t) \in \Xbt_\cdo^1$, 
\equat\label{bien defined}
\lexp{t\wdo_0^{-1}\l(z)}{\m(z)}=zv\vb(z,t).
\endequat
Now, let 
$$\fonction{\xbt_\cdo^\pr}{\Xbt_\cdo^1 \times 
A_{\Ub_\Lb}(v)}{\Xbt_\cdo^\pr}{(z,t,a)}{\l(z)^{-1}\wdo_0 t^{-1} *_\Pb 
(aC_\Lb^\ci(v),z,\vb(z,t)).}$$
It is a well-defined \mor of varieties by \ref{bien defined}. Finally, let
$$\gti_\cdo : \Xbt_\cdo^1 \times A_{\Ub_\Lb}(v) \longto \Zb(\Lb)^\ci$$
denote the projection on the first factor of $\Xbt_\cdo^1$.

\bi

\begin{theo}\label{iso eta}
Assume that $\Gb$ is semisimple and simply connected. Then 
$\xbt_\cdo^\pr$ is an $A_\Lb(v)$-equivariant \iso of varieties, and the diagram
$$\diagram
\Xbt_\cdo^1 \times A_{\Ub_\Lb}(v) 
\rrto^{\DS{\xbt_\cdo^\pr}} \ddto_{\DS{\gti_\cdo}} && \Xbt_\cdo^\pr \ddto^{\DS{\fti_\cdo}} \\
&&\\
\Zb(\Lb)^\ci \rrto^{\DS{\xbt_\cdo}} && \Xbt_\cdo
\enddiagram$$
is commutative.
\end{theo}

\bi

\proof The commutativity of the diagram is clear. The equivariance under the action 
of $A_\Lb(v)$ is also clear. Now, let us prove that $\xbt_\cdo^\pr$ is 
an isomorphism of varieties. For this, we construct an explicit inverse 
of $\xbt_\cdo^\pr$. Let 
$$\dot{\Xb}_\cdo=\{(g,l,z,x) \in \Gb \times \Lb \times \Zb(\Lb)^\ci 
\times \Vb~|~\lexp{g}{(lzvl^{-1}x)} \in \Xb_\cdo\}.$$
Let $(g,l,z,x) \in \dot{\Xb}_\cdo$. Then, by Lemma \ref{V conjugue}, 
the semisimple part of $\lexp{g}{(lzvl^{-1}x)}$ is conjugate to $z$, 
so $\lexp{g}{(lzvl^{-1}x)}=\m(z)$. In particular,
$$zv(l^{-1}xl)=\lexp{l^{-1}g^{-1}\l(z)^{-1}\wdo_0}{(\wdo_0^{-1}\o^{-1}(z)\wdo_0)}.$$
But $zv(l^{-1}xl)$ and $\wdo_0^{-1}\o^{-1}(z)\wdo_0$ are regular elements of 
$\Bb$ having the same image under $\pi_\Tb$ (here, $z$). Since $\Gb$ is 
semisimple and simply connected, we get from Lemma \ref{conjugue regulier} (3) 
that $l^{-1}g^{-1}\l(z)^{-1}\wdo_0$. We denote by $t=\sb(g,l,z,x)$ the projection 
of $\pi_\Tb(l^{-1}g^{-1}\l(z)^{-1}\wdo_0)$ on $\Sb$. Then $\sb : \dot{\Xb}_\cdo 
\to \Sb$ is a \mor of varieties. Moreover, by comparing the 
coefficients of $zv(l^{-1}xl)$ and $\wdo_0^{-1}\o^{-1}(z)\wdo_0 \in \Tb\Ub_*$ 
on the simple roots belonging to $\D_\Lb$, we get that $(t,z) \in \Xbt_\cdo^1$. 
In particular, $\pi_\Lb(l^{-1}g^{-1}\l(z)^{-1}\wdo_0 t^{-1})$ belongs to 
$\Zb(\Lb)^\ci C_{\Ub_\Lb}(v)$. We denote by $\ab(g,l,z,x)$ the image of 
$\pi_\Lb(l^{-1}g^{-1}\l(z)^{-1}\wdo_0 \sb(g,l,z,x)^{-1})$ in $A_\Lb(v)$~: in fact, 
it belongs to $A_{\Ub_\Lb}(v)$. We have constructed a morphism
$$\fonction{f_0}{\dot{\Xb}_\cdo}{\Xbt_\cdo^1 \times A_{\Ub_\Lb}(v)}{(g,l,z,x)}{
(z,\sb(g,l,z,x),\ab(g,l,z,x)).}$$
It is straightforward that $f_0$ factorizes through the canonical quotient morphism
$\dot{\Xb}_\cdo \to \Xbt_\cdo^\pr$. We denote by
$$\nablat_\cdo^\pr : \Xbt_\cdo^\pr \longto \Xbt_\cdo^1 \times A_{\Ub_\Lb}(v)$$
the morphism obtain after factorization. An immediate computation shows 
that $\nablat_\cdo^\pr \ci \xbt_\cdo^\pr =\Id_{\Xbt_\cdo^1 \times A_{\Ub_\Lb}(v)}$. 
To conclude the proof of Theorem \ref{iso eta}, we just need to show 
that $\nablat_\cdo'$ is injective.

Let $g*_\Pb (lC_\Lb^\ci(v),z,x)$ and $g'*_\Pb (l'C_\Lb^\ci(v),z',x')$ 
be two elements of $\Xbt_\cdo$ having the same image under $\nablat_\cdo^\pr$. 
By changing of representatives, we may, and we will, assume that 
$l=l'=1$. Also, it is clear that $z=z'$. 
But then, we have 
$$\lexp{g}{(zvx)}=\m(z)=\m(z')=\lexp{g'}{(zvx')}$$
that is
$$zvx=\lexp{g^{-1}g'}{(zvx')}.$$
By Lemma \ref{conjugue regulier} (3), we deduce that $g^{-1}g' \in \Bb$. 
Moreover, the projection of $g^{-1}g'$ on $\Lb$ belongs to $C_\Lb(v)$. 
But, since $\sb(g,1,z,x)=\sb(g',1,z,x')$ and $\ab(g,1,z,x)=\ab(g',1,z,x')$, 
it is easy to deduce that this projection belongs to $C_\Lb^\ci(v)$. 
We write $g^{-1}g'=my$, where $m \in C_\Lb^\ci(v)$ and $y \in \Vb$. 
Then, by the definition of the action of $\Pb$, we get that 
$g*_\Pb (C_\Lb^\ci(v),z,x)=g'*_\Pb (l'C_\Lb^\ci(v),z',x')$.\fin

\bi

\remark{combinaison} 
Combining Theorems \ref{key} and \ref{iso eta}, we get that the diagram
$$\diagram
\Xbt_\cdo^\pr \ddto^{\DS{\fti_\cdo}} 
\rrto^{\DS{\nablat_\cdo^\pr}} && \Xbt_\cdo^1 \times A_{\Ub_\Lb}(v) 
\ddto^{\DS{\gti_\cdo}} \\
&&\\
\Xbt_\cdo \ddto^{\DS{\pi_\cdo}} \rrto^{\DS{\nablat_\cdo}} && 
\Zb(\Lb)^\ci \ddto^{\DS{\nabla}}\\
&& \\
\Xb_\cdo \rrto^{\DS{\na_\cdo}} && \na(\Zb(\Lb)^\ci) 
\enddiagram$$
is commutative. Each horizontal map is an isomorphism.\finl

\bi

%%The aim of this Appendix is to provide a proof for the results given by Table 
%%\ref{tabletable}. By Remark \ref{reduction max}, we only need to make the 
%%computation whenever 
%%$\Gb$ is semisimple, simply connected, and quasisimple, and $\Lb$ is 
%%a maximal quasi-cuspidal \levi of a \para of $\Gb$. Using the classification 
%%of quasi-cuspidal Levi subgroups given by Table \ref{tabletable}, 
%%one sees that we only need to investigate the following 
%%cases~:

\sec{The morphism $\ph_{\Lb,v}^\Gb$ for $v$ regular\label{phi sec}}~

\med

We turn back to the general situation, that is, we no longer assume that 
$\Gb$ is semisimple and simply connected. 
In view of applications to Gel'fand-Graev characters, we need to compute 
explicitly the morphism $\ph_{\Lb,v}^\Gb$ for $v$ regular, 
and $\Lb$ {\it cuspidal} (see \SEC\ref{cuspidal sous} for the 
definition). However, we are only able to do it whenever $p$ is good~: 
the result is then given by Table \ref{tabletable}. The proof 
is done by a case-by-case analysis, after some reduction arguments 
coming from Part I. It involves explicit 
computations in the root system of $\Gb$, which are necessary 
to determine the structure of the variety $\Xbt_\cdo^\pr$ as it is shown 
by Theorem \ref{iso eta}. This proof is done in Sections 
\ref{cas A} and \ref{cas B}. It requires mainly technical computations.

\bigskip

\subsection{Cuspidal groups\label{cuspidal sous}} 
The morphism $h_\Lb^\Gb : \ZC(\Gb) \to \ZC(\Lb)$ 
(recall that $\ZC(\Gb)=\Zb(\Gb)/\Zb(\Gb)^\ci$) is well-known 
to be surjective (see for instance \cite[Lemma 1.4]{DLM1}). If there is no ambiguity, 
we will denote it simply by $h_\Lb$. The group $\Gb$ is said {\it cuspidal} 
if $\Ker h_\Lb \not=\{1\}$ for every \levi $\Lb$ of a proper \para of $\Gb$.

\bi

\remark{comparaison} Note that this definition coincide with the one 
given in \cite[\SEC 1]{Bonnafe} or \cite[\SEC 2]{bonnafe regular}.\finl

\bi

We recall some results on cuspidal groups. For a proof, the reader may refer 
to \cite[\SEC 2 and 3]{bonnafe regular} for instance.

\bigskip

\begin{prop}\label{cuspidal isotypique}
$({\mathrm{a}})$ 
If $\Gb$ is cuspidal and if $\s : \Gbh \to \Gb$ is an isotypic morphism, 
then $\Gbh$ is cuspidal.

\tete{b} A cuspidal group is universally self-opposed (see \SEC\ref{distingue sub} 
for the definition).

\tete{c}  If $\Gb$ is cuspidal then any irreducible component of 
$\Phi$ is of type $A$.

\tete{d} If $\Lb$ is cuspidal, then $A_\Lb(v)=A_\Gb(v)$.
\end{prop}

\bigskip

Proposition \ref{cuspidal isotypique} (d) and Corollary \ref{CGU} show 
that the morphism $\ph_{\Lb,v}^\Gb$ is well-defined whenever $\Lb$ is cuspidal.

\bigskip

\noindent{\bf Hypothesis :} {\it From now on, and until the end of this paper, 
we assume that $\Lb$ is cuspidal and that $p$ is good for $\Gb$.}

\bigskip

We start this section by a compatibility property of the family of morphisms 
$\ph_{\Lb,v}^\Gb$ whenever $\Lb$ runs over the set of cuspidal Levi 
subgroups of \para of $\Gb$. We prove this fact {\it a priori}. 
It could have been checked using the explicit computation of these 
morphisms but we have thought that it is more interesting to 
give a proof not relying on the classification. 
This compatibility is crucial for the application to 
Gel'fand-Graev characters. 

\bi

\subsection{Compatibility\label{compatible sub}} 
Let $\Pb'$ be a \para of $\Gb$ containing $\Pb$ and let $\Lb'$ denote 
the unique \levi of $\Pb'$ containing $\Lb$. Let $\Vb'$ denote the unipotent 
radical of $\Pb'$ and let $v'=\pi_{\Lb'}(u)$. So we are in the 
set-up of \SEC\ref{para sub red}. 
Note that $v'$ is regular and that $v=\pi_\Lb(v')=\pi_\Lb(u)$. 
In \cite[\SEC 4]{bonnafe regular}, 
the author defined a morphism
$$\r_{\Lb,\Lb',u} : N_\Gb(\Lb',v') \to N_\Gb(\Lb,v).$$
By \cite[\SEC 4]{bonnafe regular}, this morphism induces 
a morphism
$$\rhoba_{\Lb,\Lb',u} : W_\Gb(\Lb',v') \to W_\Gb(\Lb,v).$$
and its restriction to $A_{\Lb'}(v') \simeq \ZC(\Lb')$ coincide 
with the map $h_\Lb^{\Lb'} : \ZC(\Lb') \to \ZC(\Lb)$. 

\bigskip

\begin{prop}\label{inclusion stabilisateurs}
If $\Lb$ is cuspidal and $p$ is good for $\Gb$, then 
$$\rhoba_{\Lb,\Lb',u}(H_\Gb(\Lb',v')) \incl H_\Gb(\Lb,v).$$
\end{prop}

\bigskip

\proof By the commutativity of the Diagram 4.7 in \cite{bonnafe regular}, 
we may replace $u$ by any other regular unipotent element. Therefore, we may 
(and we will) 
assume that $u=\beta(\sum_{\a \in \D} e_\a)$, where 
$\beta : \gGB_\nil \to \Gb_\uni$ is a $\Gb$-equivariant morphism 
of varieties inducing an homeomorphism between the underlying topological spaces. 
Let $\r^\vee$, $\r_\Lb^\vee$ and $\r_{\Lb'}^\vee$ denote the sum of the 
positive coroots (relative to $\Tb$) 
of $\Gb$, $\Lb$ an $\Lb'$ respectively. By \cite[End of \SEC 5.B]{bonnafe regular}, 
the canonical morphisms 
$\Bigl(N_\Gb(\Lb,v) \cap C_\Gb(\Im \r_\Lb^\vee)\Bigr)/\Zb(\Lb)^\circ \to W_\Gb(\Lb,v)$ 
and $\Bigl(N_\Gb(\Lb',v') \cap C_\Gb(\Im \r_{\Lb'}^\vee)\Bigr)/\Zb(\Lb')^\circ 
\to W_\Gb(\Lb',v')$ are isomorphisms. Therefore, 
if $w$ belongs to $W_\Gb(\Lb,v)$ or $W_\Gb(\Lb',v')$, we denote by $\wdo$ a 
representative of $w$ lying in $N_\Gb(\Lb,v) \cap C_\Gb(\Im \r_\Lb^\vee)$ 
or $N_\Gb(\Lb',v') \cap C_\Gb(\Im \r_{\Lb'}^\vee)$ respectively. 
In particular, $\wdo$ normalizes 
$C_\Lb(\Im \r_\Lb^\vee)=C_{\Lb'}(\Im \r_{\Lb'}^\vee)=\Tb$. 
On the other hand, we have \cite[5.10]{bonnafe regular} 
$$N_\Gb(\Lb',v') \cap C_\Gb(\Im \r_{\Lb'}^\vee) \incl 
N_\Gb(\Lb,v) \cap C_\Gb(\Im \r_\Lb^\vee)$$
and the morphism $\rhoba_{\Lb,\Lb',u}$ is induced by this inclusion.

The Borel subgroup $\Bb_{\Lb'}=\Bb \cap \Lb'$ of $\Lb'$ contains $C_{\Lb'}^\circ(v')$. 
Therefore, 
$\Lb'/C_{\Lb'}^\circ(v') \simeq \Lb' \times_{\Bb_{\Lb'}} \Bb_{\Lb'}/C_{\Lb'}^\circ(v')$.   
Consequently, 
$$\Sigb_{\Lb'}' \times \Vb' \simeq \Pb' \times_{\Bb} 
(\Bb_{\Lb'}/C_{\Lb'}^\circ(v) \times \Zb(\Lb')^\circ \times \Vb'),\leqno{(\#)}$$
where $\Sigb_{\Lb'}'$ is the analogous of the variety $\Sigb'$ but defined using 
$\Lb'$ instead of $\Lb$, and $\Bb=\Bb_{\Lb'}.\Vb'$ acts on 
$(\Bb_{\Lb'}/C_{\Lb'}^\circ(v') \times \Zb(\Lb')^\circ \times \Vb')$ by restriction 
of the action of $\Pb'$. 
Finally, we get an isomorphism
$$\Xbt_{0,\Lb'}' \simeq \Gb \times_\Bb 
(\Bb_{\Lb'}/C_{\Lb'}^\circ(v') \times \Zb(\Lb')^\circ \times \Vb'),$$
where $\Xbt_{0,\Lb'}'$ is the analogous of the variety $\Xbt_0^\prime$ 
but defined using $\Lb'$ instead of $\Lb$. We can then define a morphism 
$$\fonction{\O}{\Xbt_{0,\Lb'}'}{\Xbt_0^\prime}{g *_\Bb 
(b' C_{\Lb'}^\circ(v'),z,x)}{g *_\Bb (\pi_\Lb(b') C_\Lb^\circ(v),z, 
\pi_\Lb(b'v'b^{\prime -1})^{-1}b'v'b^{\prime -1} x).}$$ 
The reader can check that $\O$ is a well-defined 
morphism of varieties. Let $\uti_{\Lb'}'$ denote the element of $\Xbt_{\Lb'}^\prime$ 
defined analogously to $\uti'$ (see \SEC{O2}). 
%%Moreover, $\O$ is $\Gb$-equivariant (this fact 
%%will not be used in this proof). 
It is clear that $\O(\uti_{\Lb'}')=\uti_\Lb'$. By restriction to $\Xbt_{\min,\Lb'}'$, 
we get a morphism $\O_\mini : \Xbt_{\min,\Lb'}' \to \Xbt_\mini'$. So 
it remains to show that the morphism $\O_\mini$ is $W_\Gb(\Lb',v')$-equivariant.

By density, it is sufficient to prove that its restriction $\O_\reg$ to 
$\Xbt_{0,\Lb',\reg}^\prime$ is $W_\Gb(\Lb',v')$-equivariant. But, 
$\Xbt_{0,\Lb',\reg}^\prime \simeq \Gb/C_{\Lb'}^\circ(v') \times \Zb(\Lb')^\circ_\reg$ 
and
$$\fonction{\O_\reg}{\Xb_{0,\Lb',\reg}^\prime}{\Xbt_0'}{(g C_{\Lb'}^\circ(v'),z)}{
(g *_\Bb (C_\Lb^\circ(v),z,v^{-1}v'))}$$
Let $\Zb(\Lb)_{\reg,+}^\circ=\{z \in \Zb(\Lb)^\circ~|~C_\Gb(z) \incl \Lb'\}$. Then 
$\Zb(\Lb)_{\reg,+}^\circ$ is an open subset of $\Zb(\Lb)^\circ$ containing 
$\Zb(\Lb)_\reg^\circ \cup \Zb(\Lb')_\reg^\circ$. We set 
$$\Xbt_{0,+}^\prime=\Gb \times_\Bb 
(\Bb_\Lb/C_\Lb^\circ(v) \times \Zb(\Lb)_{\reg,+}^\circ \times \Vb).$$
Then the image of $\O_\reg$ is contained in $\Xbt_{0,+}^\prime$. Moreover, 
by Corollary \ref{U conjugue bis}, we have
$$\Xbt_{0,+}^\prime\simeq\Gb \times_{\Bb_{\Lb'}} 
(\Bb_\Lb/C_\Lb^\circ(v) \times \Zb(\Lb)_{\reg,+}^\circ \times (\Vb \cap \Lb')).$$
Now, if $w \in W_\Gb(\Lb',v')$, then $w$ acts on $\Xbt_{0,+}^\prime$ as follows. 
If $g *_{\Bb_{\Lb'}} (bC_\Lb^\circ(v),z,x) \in \Xbt_{0,+}^\prime$ then 
$$(g *_{\Bb_{\Lb'}} (bC_\Lb^\circ(v),z,x)).w=g\wdo *_{\Bb_{\Lb'}} (\wdo^{-1}b\wdo 
C_{\Lb}^\circ(v),\wdo^{-1} z \wdo, \wdo^{-1} x \wdo).$$
Indeed, this follows fromthe fact that $\wdo$ normalizes $\Bb_{\Lb'}$, 
$\Vb \cap \Lb'$ and $C_\Lb(v)$ (see \cite[\SEC 5]{bonnafe regular}
This action coincides with the action of $\rhoba_{\Lb,\Lb',u}(w) \in W_\Gb(\Lb,v)$, 
as it can be checked by restriction to the open subset $\Ybt'$. 
Therefore, $\O_\reg$ is $W_\Gb(\Lb',v')$-equivariant.\fin

\bigskip

\begin{coro}\label{compatible prop}
If $\Lb$ and $\Lb'$ are both cuspidal, 
then $\rhoba_{\Lb,\Lb',u}(W_\Gb^\circ(\Lb',v')) \incl W_\Gb^\circ(\Lb,v)$ and 
%%$\rhoba_{\Lb,\Lb',u}(H_\Gb(\Lb',v')) \incl H_\Gb(\Lb,v)$. Therefore, 
the diagram 
$$\diagram
W_\Gb^\circ(\Lb',v') \rrto^{\DS{\ph_{\Lb',v'}^\Gb}} \ddto_{\DS{\rhoba_{\Lb,\Lb',u}}} && 
\ZC(\Lb') \ddto^{\DS{h_\Lb^{\Lb'}}} \\
&&\\
W_\Gb^\circ(\Lb,v) \rrto^{\DS{\ph_{\Lb,v}^\Gb}} && \ZC(\Lb) \\
\enddiagram$$
is commutative.
\end{coro}

\bigskip

\proof The fact that 
$\rhoba_{\Lb,\Lb',u}(W_\Gb^\circ(\Lb',v')) \incl W_\Gb^\circ(\Lb,v)$ 
has been proved in \cite[Proposition 6.3]{bonnafe regular}. The 
commutativity of the diagram 
then follows from Proposition \ref{inclusion stabilisateurs}.\fin

\bigskip

\sub{The classification} 
The Table \ref{tabletable} gives the values of the morphism $\ph_{\Lb,v}^\Gb$ 
whenever $v$ is regular, $\Gb$ is semisimple, simply connected and 
quasi-simple, $\Lb$ is cuspidal, and $p$ is good for $\Lb$ 
(for the classification of pairs $(\Gb,\Lb)$ where $\Gb$ is semisimple, 
simply connected and quasi-simple and $\Lb$ is cuspidal, the reader may refer 
to \cite[Table 2.17]{bonnafe regular}). 

Table \ref{tabletable} must be read as follows. 
In the first column, the diagram of $\Gb$ is given. The diagram 
of $\Lb$ corresponds to the black nodes (except if in the first row of 
the Table in which case $\Lb$ is of type $A_{d-1} \times \dots A_{d-1}$. 
To each white node 
(which corresponds to the elements $\a$ of $\D-\D_\Lb$) is associated 
an element of order $2$ of $W_\Gb(\Lb)\simeq 
W_\Gb^\ci(\Lb,v)$ (\cf Proposition \ref{propriete distingues}~: it is a reflection 
for its action on $X(\Zb(\Lb)^\ci)$). The number ($1$ or $-1$ or $\varpi_1^\vee$) 
written just behind this white node is equal to the value of 
$\ph_{\Lb,v}^\Gb$ at this element of order $2$. Note that it belongs 
to $A_\Lb(v) \simeq \ZC(\Lb)$ which is given in the 
third column. 
The second column gives the type of the Weyl group 
$W_\Gb^\circ(\Lb,v) \simeq W_\Gb(\Lb)$. 
The last column gives some conditions on $p$ or $r$ 
for which the situation is possible. 

We need to give an explanation for the value $\varpi_1^\ve$ which appears in type 
$D_r$ for $r$ even. If we number the simple roots of $\Gb$ 
as in \cite[Page 256]{bourbaki}, that is
$$
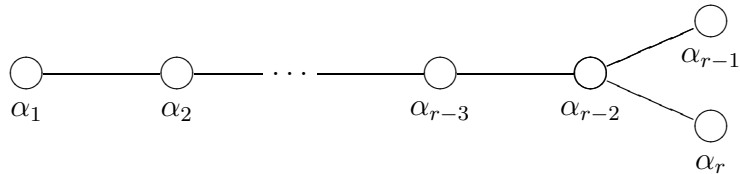
\xy (0,0) *++={\phan} *\frm{o} ; (20,0) *++={\phan} *\frm{o} **@{-};
(35,0) *++={\dots} **@{-}; (55,0) *++={\phan} *\frm{o} **@{-};
(75,0) *++={\phan} *\frm{o} **@{-};(91,7) *++={\phan} *\frm{o} **@{-};
(75,0) *++={\phan} *\frm{o} **@{-};(91,-7) *++={\phan} *\frm{o} **@{-};
(0,-5) *+={\a_1};(20,-5) *+={\a_2}; (55,-5) *+={\a_{r-3}}; 
(75,-5) *+={\a_{r-2}}; (91,2) *+={\a_{r-1}}; (91,-12) *+={\a_r}\endxy$$
then $\varpi_1^\vee$ denotes the minuscule weight associated to $\a_1$, 
and we identify it with its image in $\Zb(\Gb)$, which is isomorphic 
to $\ZC(\Lb)$ in the unique situation where it appears.

\bi

\rem As it is explained in Section \ref{ele sec}, Table \ref{tabletable} 
is sufficient for determining the morphism $\ph_{\Lb,v}^\Gb$ 
for every reductive group $\Gb$ (not necessarily semisimple and 
simply connected) under the following hypothesis~: $v$ is regular, $p$ is good 
for $\Lb$, and $\Lb$ is cuspidal.\finl

\bi

\begin{table}
\begin{tabular}{|c|c|c|c|}
\hline
$\vphantom{\DS{\sum_A^B}}$   $\Gb$, $\Lb$, $\ph_{\Lb,v}^\Gb $ & Type of $W_\Gb(\Lb)$ & 
 $\ZC(\Lb)$ & Remarks \\
\hline
$\vphantom{\xy (0,8) *+={a} ; (0,-2) *+={a};\endxy }$ 
$\xy
(0,0) *+={~A_{d-1}~} *\frm{-} ;(17,0) *+={\phan} *\frm{o} **@{-} ;
(34,0) *+={~A_{d-1}~} *\frm{-}**@{-}; (51,0) *+={\phan} *\frm{o} **@{-} ;
(61,0) *+={~\dots~} **@{-};(71,0) *+={\phan} *\frm{o} **@{-};
(88,0) *+={~A_{d-1}~} *\frm{-}**@{-};
(17,-4) *+={(-1)^{d-1}} ; (51,-4) *+={(-1)^{d-1}} ;
(71,-4) *+={(-1)^{d-1}} ;
\endxy $ & 
$A_{{r+1 \over d}-1}$&$\mub_d$ & $\xy (0,3) *+={p \not| ~d}; (0,-3) *+={d~|~r+1}\endxy$ \\
\hline
$\vphantom{\xy (0,4) *+={a} ; (0,-6) *+={a};\endxy }$
$\xy (0,0) *+={\bullet} *\frm{o} ; (15,0) *+={\phan} *\frm{o} **@{-} ;
(25,0) *+={~\dots~}  **@{-} ; (35,0) *+={\phan} *\frm{o} **@{-} ;
(50,0) *+={\bullet} *\frm{o} **@{-};(65,0) *+={\phan} *\frm{o} **@{-};
(80,0) *+={\bullet} *\frm{o}**@{=}; (72.5,0) *+={>};
(15,-3) *+={-1} ; (35,-3) *+={-1} ;
(65,-3) *+={1}\endxy$ & 
$B_{r-1 \over 2}$ &$\mub_2$ & $\xy (0,3) *+={p \not= 2};(0,-3) *+={r~\text{odd}}\endxy$\\
\hline
$\vphantom{\xy (0,4) *+={a} ; (0,-6) *+={a};\endxy }$
$\xy (0,0) *+={\bullet} *\frm{o} ; (15,0) *+={\phan} *\frm{o} **@{-} ;
(30,0) *+={\bullet} *\frm{o} **@{-} ; (45,0) *+={\phan} *\frm{o} **@{-} ;
(55,0)  *+={~\dots~} **@{-};(65,0) *+={\bullet} *\frm{o} **@{-};
(80,0) *+={\phan} *\frm{o} **@{=}; (72.5,0) *+={>};
(15,-3) *+={-1} ; (45,-3) *+={-1} ;
(80,-3) *+={-1}\endxy$ & 
$B_{r \over 2}$ & $\mub_2$ & $\xy (0,3) *+={p \not= 2}; 
(0,-3) *+={r~\text{even}}\endxy$\\
\hline
$\vphantom{\xy (0,4) *+={a} ; (0,-6) *+={a};\endxy }$
$\xy (0,0) *+={\phan} *\frm{o};  (15,0) *+={~\dots~} **@{-};
(30,0) *+={\phan} *\frm{o} **@{-}; (50,0) *+={\phan} *\frm{o} **@{-};
(70,0)*+={\bullet} *\frm{o} **@{=}; (60,0) *+={<};
(0,-3) *+={1} ; (30,-3) *+={1}; (50,-3) *+={-1}\endxy$ & 
$B_{r-1}$ & $\mub_2$ & $p \not= 2$ \\
\hline
$\vphantom{\xy (0,6) *+={a} ; (0,-6) *+={a};\endxy }$
$\xy (0,0) *+={\bullet} *\frm{o} ; (15,0) *+={\phan} *\frm{o} **@{-} ;
(30,0) *+={\bullet} *\frm{o} **@{-} ; (45,0) *+={\phan} *\frm{o} **@{-} ;
(55,0)  *+={~\dots~} **@{-};(65,0) *+={\bullet} *\frm{o} **@{-};
(80,0) *+={\phan} *\frm{o} **@{-}; (92,4) *+={\bullet} *\frm{o} **@{-};
(80,0) *+={\phan} *\frm{o} **@{-};(92,-4) *+={\bullet} *\frm{o} **@{-};
(79,-3) *+={1}; (15,-3) *+={\varpi_1^\ve}; (45,-3) *+={\varpi_1^\ve};  \endxy$ & 
$B_{r-2 \over 2}$ &
$\mub_2 \times \mub_2$ & $\xy (0,3) *+={p \not= 2}; (0,-3) *+={r~\text{even}}\endxy$ \\
\hline
$\vphantom{\xy (0,6) *+={a} ; (0,-6) *+={a};\endxy }$ 
$\xy (0,0) *+={\bullet} *\frm{o} ; (15,0) *+={\phan} *\frm{o} **@{-} ;
(30,0) *+={\bullet} *\frm{o} **@{-} ; (45,0) *+={\phan} *\frm{o} **@{-} ;
(55,0)  *+={~\dots~} **@{-};(65,0) *+={\bullet} *\frm{o} **@{-};
(80,0) *+={\phan} *\frm{o} **@{-}; (92,4) *+={\phan} *\frm{o} **@{-};
(80,0) *+={\phan} *\frm{o} **@{-};(92,-4) *+={\bullet} *\frm{o} **@{-};
(79,-3) *+={-1}; (15,-3) *+={-1}; (45,-3) *+={-1}; (95,2) *+={1} \endxy$ & 
$B_{r \over 2}$&$\mub_2$ & $\xy (0,3) *+={p \not= 2}; (0,-3) *+={r~\text{even}}\endxy$ \\
\hline
$\vphantom{\xy (0,6) *+={a} ; (0,-6) *+={a};\endxy }$ 
$\xy (0,0) *+={\bullet} *\frm{o} ; (15,0) *+={\phan} *\frm{o} **@{-} ;
(30,0) *+={\bullet} *\frm{o} **@{-} ; (45,0) *+={\phan} *\frm{o} **@{-} ;
(55,0)  *+={~\dots~} **@{-};(65,0) *+={\bullet} *\frm{o} **@{-};
(80,0) *+={\phan} *\frm{o} **@{-}; (92,4) *+={\bullet} *\frm{o} **@{-};
(80,0) *+={\phan} *\frm{o} **@{-};(92,-4) *+={\phan} *\frm{o} **@{-};
(79,-3) *+={-1}; (15,-3) *+={-1}; (45,-3) *+={-1}; (95,-2) *+={1} \endxy$ & 
$B_{r \over 2}$&$\mub_2$ & $\xy (0,3) *+={p \not= 2}; (0,-3) *+={r~\text{even}}\endxy$ \\
\hline
$\vphantom{\xy (0,6) *+={a} ; (0,-6) *+={a};\endxy }$ 
$\xy (0,0) *+={\bullet} *\frm{o} ; (15,0) *+={\phan} *\frm{o} **@{-} ;
(30,0) *+={\bullet} *\frm{o} **@{-} ; (40,0) *+={~\dots~}  **@{-} ;
(50,0)  *+={\bullet} *\frm{o}**@{-};(65,0) *+={\phan} *\frm{o} **@{-};
(80,0) *+={\bullet} *\frm{o} **@{-}; (92,4) *+={\bullet} *\frm{o} **@{-};
(80,0) *+={\bullet} *\frm{o} **@{-};(92,-4) *+={\bullet} *\frm{o} **@{-};
 (15,-3) *+={-1}; (65,-3) *+={-1};  \endxy$ & 
$B_{r-3 \over 2}$&$\mub_4$ & 
$\xy (0,3) *+={p \not= 2}; (0,-3) *+={r~\text{odd}}\endxy$ \\
\hline
$\vphantom{\xy (0,6) *+={a} ; (0,-6) *+={a};\endxy }$ 
$\xy (0,0) *+={\phan} *\frm{o} ; (15,0) *+={\phan} *\frm{o} **@{-} ;
(30,0) *+={~\dots~}  **@{-} ;
(40,0)  *+={\phan} *\frm{o}**@{-};
(55,0) *+={\phan} *\frm{o} **@{-}; (67,4) *+={\bullet} *\frm{o} **@{-};
(55,0) *+={\phan} *\frm{o} **@{-};(67,-4) *+={\bullet} *\frm{o} **@{-};
 (0,-3) *+={1}; (15,-3) *+={1};(40,-3) *+={1}; (55,-3) *+={-1};   \endxy$ & 
$B_{r-2}$ & $\mub_2$ & $p \not= 2$ \\
\hline
$\vphantom{\xy (0,5) *+={a} ; (0,-12) *+={a};\endxy }$ 
$\xy (0,0) *+={\bullet} *\frm{o} ; (15,0) *+={\bullet} *\frm{o} **@{-} ;
(30,0) *+={\phan} *\frm{o} **@{-}; (45,0) *+={\bullet} *\frm{o} **@{-} ; 
(60,0)*+={\bullet} *\frm{o} **@{-} ; 
(30,0) *+={\phan} *\frm{o}  ; (30,-9) *+={\phan} *\frm{o} **@{-}; 
(27,-3) *+={1} ; (27,-9) *+={1};
\endxy$ & $G_2$ & $\mub_3$ & $p \not= 3$ \\
\hline
$\vphantom{\xy (0,5) *+={a} ; (0,-12) *+={a};\endxy }$ 
$\xy (0,0) *+={\phan} *\frm{o} ; (15,0) *+={\phan} *\frm{o} **@{-} ;
(30,0) *+={\phan} *\frm{o} **@{-}; (45,0) *+={\bullet} *\frm{o} **@{-} ; 
(60,0)*+={\phan} *\frm{o} **@{-} ; (75,0) *+={\bullet} *\frm{o} **@{-} ;
(30,0) *+={\phan} *\frm{o}  ; (30,-9) *+={\bullet} *\frm{o} **@{-}; 
(0,-3) *+={1} ; (15,-3) *+={1}; (34,-3) *+={-1}; (60,-3) *+={-1};
\endxy$ & $F_4$ & $\mub_2$ & $p \not= 2$ \\
\hline
\end{tabular}

\bi

\begin{centerline}{\refstepcounter{theo}\label{tabletable}\noindent\bf Table   
\arabic{section}.\arabic{theo}}\end{centerline}
\end{table}

\bi

%\sub{A compatibility property\label{co su}} 
%Let $\Pb'$ be a \para of $\Gb$ containing $\Pb$ 
%and let $\Lb'$ be the unique \levi of $\Pb'$ containing $\Lb$. 
%As in \SEC\ref{para sub red}, we define $v'=\pi_{\Lb'}(u)$. Then $v'$ is 
%a regular unipotent element of $\Lb'$. We assume in this subsection that 
%{\it $\Lb$ and $\Lb'$ are quasi-cuspidal}. Let 
%$$\WC_\Lb=\{w \in W~|~w(\D_\Lb)=\D_\Lb\}\hspace{0.7cm}\text{and}\hspace{0.7cm}
%\WC_{\Lb'}=\{w \in W~|~w(\D_{\Lb'})=\D_{\Lb'}\}.$$
%Then $\WC_\Lb \simeq W_\Gb(\Lb) \simeq W_\Gb^\ci(\Lb,v)$, and, similarly,
%$\WC_{\Lb'}\simeq W_\Gb^\ci(\Lb',v')$. 
%ut, by Proposition \ref{propriete distingues}, we have $\WC_{\Lb'} \incl 
%\WC_\Lb$. This gives rise to an injective morphism $\iota_{\Lb',\Lb} : 
%_\Gb^\ci(\Lb',v') \injto W_\Gb^\ci(\Lb,v)$. %
%
%Then it is straightforward to check from Table \ref{tabletable} 
%that the following result holds~:%
%
%\bi
%
%\begin{prop}\label{compatibilite regulier}
%If $\Lb$ and $\Lb'$ are quasi-cuspidal, and if $p$ is good for $\Lb'$, then 
%the diagram
%$$\diagram
%W_\Gb^\ci(\Lb',v') \rrto^{\DS{\iota_{\Lb',\Lb}}} \ddto_{\DS{\ph_{\Lb',v'}}}&& 
%W_\Gb^\ci(\Lb,v) \ddto^{\DS{\ph_{\Lb,v}}} \\
%&&\\
%A_{\Lb'}(v') \rrto_{\DS{h_\Lb^{\Lb'}}} && A_\Lb(v) 
%\enddiagram$$
%is commutative.
%\end{prop}
%
%\bi

\sub{Reductions} Before coming to the proof of the results given by Table 
\ref{tabletable} (which will be done in Sections \ref{cas A} and \ref{cas B}), 
we reduce the number of cases which remain to be treated by brute force 
using results from Part I. By Proposition 1.13 (b) and Corollary 4.7, 
it is sufficient to deal with the following cases~:

\med

{\sc Case A :} $(\Gb,\Lb)$ is of type $(A_r,A_{d-1} \times A_{d-1})$, 
with $r=2d-1$ odd and $p$ does not divide $d$.

\med

{\sc Case B :} $(\Gb,\Lb)$ is of type $(B_3,A_1 \times A_1)$, and $p \not= 2$.

\med

{\sc Case C :} $(\Gb,\Lb)$ is of type $(C_2,A_1)$, $p \not= 2$, and the root 
of $\Lb$ is a long root for $\Gb$.

\med

{\sc Case D :} $(\Gb,\Lb)$ is of type $(D_5,A_1 \times A_3)$, and $p\not= 2$.

\med

{\sc Case E :} $(\Gb,\Lb)$ is of type $(D_4,A_1 \times A_1 \times A_1)$, and $p\not= 2$.

\bi

In Section \ref{cas A}, we will give a proof for the case $A$, while 
Section \ref{cas B} is devoted to the remaining cases, for which 
we use Theorem \ref{iso eta} to perform the computations  
in the root system of $\Gb$. 

\bi

\sec{Groups of type $A$\label{cas A}}~

\med

{\it We assume in this section, and only in this section, 
that $\Gb=\Sb\Lb_{r+1}(\FM)$, that 
$$\Lb=\{\diag(A_1,\dots,A_k)~|~\Bigl(\forall 1 \le i \le k,~A_i \in \Gb\Lb_d(\FM) \Bigr)
{\mathrm{~and~}} \det(A_1\dots A_k)=1\},$$
where $r+1=dk$, and that $p$ does not divide $d$. 
For simplifying notations, we set $n=r+1$. Note that $r'=(d-1)k=n-k$.}

\bi

\rem As noticed in the introduction of this Appendix, it would have been 
sufficient to handle with the case $k=2$. However, the general 
case may be treated as well.\finl

\bi

\sub{Notation} 
We choose for $\Bb$ the group of lower triangular matrices of $\Gb$, for $\Tb$ the group 
of diagonal matrices of $\Gb$, and for $\cdo$ the following matrix~:
$$\cdo=\matrice{0 & \dots & \dots & 0 & (-1)^{n-1} \\
                1 & 0 & \dots &  \dots & 0 \\
                0 & \ddots & \ddots & & \vdots \\
               \vdots & \ddots & \ddots & \ddots & \vdots \\
               0 & \dots & 0 & 1 & 0}.$$ 
Note that $\det \cdo=1$, so that $\cdo \in \Gb$. 
We denote by $\mub_d$ the group of $d$th roots of unity in $\FM^\times$. 
We have $|\mub_d|=d$ since $d$ is prime to $p$. 
If $z \in \FM^\times$, we denote by $\t(z)$ the matrix 
$$\diag(\unb_d,\dots,\unb_d,z \unb_d),$$
where $\unb_d$ denotes the $d \times d$ identity matrix. If $z \in \mub_d$, 
then $\t(z) \in \Zb(\Lb)$ and the morphism 
$\mub_d \to \Zb(\Lb)$, $z \mapsto \t(z)$ induces an isomorphism $\mub_d \simeq 
\ZC(\Lb)$. Note that, if $z \not\in \mub_d$, then $\t(z) \not\in \Gb$. 

If $m$ is a non-zero natural number, we denote 
by $J_m $ the $m \times m$ matrix
$$J_m=\matrice{1 & 0 & \dots & \dots & 0 \\
               1 & 1 & \ddots &  & \vdots \\
            0 & \ddots & \ddots & \ddots & \vdots \\
            \vdots & \ddots & \ddots & \ddots & 0 \\
              0 & \dots & 0 & 1 & 1}.$$
Finally, we may, and we will, assume that 
$$v=\diag(J_d,\dots,J_d).$$ 

\med

\sub{The varieties $\Tb/W$, $\Ub \cdo \cap \cdo \Ub^-$, and $\na(\Zb(\Lb)^\ci)$} 
Let $\FM_n[X]$ denote the set of monic polynomials $P$ 
in the indeterminate $X$ of degree $n$ with coefficients in $\FM$ and 
such that $P(0)=(-1)^n$. Then $\Tb/W$ is isomorphic to $\FM_n[X]$ and, 
via this isomorphism, the map $\na$ is identified with
$$\fonction{\na}{\Gb}{\FM_n[X]}{g}{\det(X\unb_n - g).}$$

If $P \in \FM_n[X]$ and if $Q \in \FM[X]$, we denote by 
$\overline{Q}^{\SSS{P}}$ the class of $Q$ in $\FM[X]/(P)$. Let $M(P)$ denote the 
matrix of the multiplication by $\overline{X}^{\SSS{P}}$ in $\FM[X]/(P)$, computed 
in the basis $\BC_P=(\overline{1}^{\SSS{P}},\overline{X}^{\SSS{P}},\dots,
(\overline{X}^{\SSS{P}})^{n-1})$. If $P=X^n+a_{n-1} X^{n-1}+\dots+a_1 X + (-1)^n$, then
$$M(P)=\matrice{0 & \dots & \dots & 0 & (-1)^{n-1} \\
                1 & \ddots &  &  \vdots & -a_1\\
                0 & \ddots & \ddots & \vdots & \vdots \\
               \vdots & \ddots & \ddots & 0 & \vdots \\
               0 & \dots & 0 & 1 & -a_{n-1}} \in \Ub\cdo \cap \cdo\Ub^-.$$ 
Hence the map $M : \FM_n[X] \to \Ub\cdo \cap \cdo\Ub^-$ is an \iso of varieties and is 
the inverse of $\na_\cdo$. 

\bi

\begin{prop}\label{technique}
The variety $\na(\Zb(\Lb)^\ci)\simeq \Xb_\cdo$ is smooth~: it 
is isomorphic to the affine space $\Ab^{k-1}(\FM)$. 
\end{prop}

\bi

\proof With these identifications, 
$\na(\Zb(\Lb)^\ci)$ is isomorphic to the image of the morphism of varieties
$$\fonction{\r}{\FM_k[X]}{\FM_n[X]}{P}{P^d.}$$
To prove Proposition \ref{technique}, it is sufficient to prove that $\r$ 
induces an \iso between $\FM_k[X]$ and its image (which is a closed subvariety of 
$\FM_n[X]$ since the composite \mor $\Zb(\Lb)^\ci \to \Tb \to \Tb/W$ 
is finite as the composite of two finite morphisms). 
But, if $P=X^k+a_{k-1}X^{k-1}+\dots+a_1 X + (-1)^k \in \FM_k[X]$ and if 
$1 \le i \le k-1$, then the coefficient of $X^{(d-1)k+i}$ in $P^d$ is of the form
$$d a_i + {\mathrm{polynomial~in}} ~a_{i+1},\dots, a_{k-1}.$$
The proof of Proposition \ref{technique} is completed by the 
fact that $d$ is prime to $p$.\fin

\bi

Now, let
$$\Ob =\{(z_1,\dots,z_k) \in (\FM^\times)^k~|~
z_1\dots z_k=1\hspace{0.3cm}{\mathrm{~and~}} 
\hspace{0.3cm}\forall 1 \le i < j \le k,~ z_i \not= z_j\}$$
$$\fonction{\db}{\Ob}{\Zb(\Lb)^\ci_\reg}{(z_1,\dots,z_k)}{\diag(z_1 \unb_d, 
\dots, z_k \unb_d).}
\leqno{\mathrm{and}}$$
Then $\db$ is an isomorphism of varieties. We will often identify in 
this subsection $\Zb(\Lb)_\reg^\ci$ with $\Ob$ via the \iso $\db$. 

If $z=(z_1,\dots,z_k)$ is an element 
of $\Ob$, we denote by $P_z$ the polynomial
$$P_z=\prod_{i=1}^k (X-z_i)^d \in \FM_n[X].$$
Note that $\m(\db(z)) = M(P_z)$ (indeed, $\na(\db(z))=P_z$), where 
$\m : \Tb \to \Ub\cdo \cap \cdo\Ub^-$ is the morphism defined at the beginning 
of Subsection \ref{compa para subsub} (\cf also Remark \ref{nabla nabla}.

\bi

\begin{prop}\label{connexite +}
There exists a \mor of varieties $\ph : \Ob \to \Gb$ such that 
$$\lexp{\ph(z)}{(\db(z)v)}=M(P_z)$$
for any $z \in \Ob$. 
\end{prop}

\proof Let $z \in \Ob$. 
We denote by $\r_z : \FM[X]/(P_z) \to \FM[X]/(P_z)$, $a \mapsto 
\overline{X}^{\SSS{P_z}} a$. 
If $1 \le i \le k$ and $1 \le j \le d$, we set
$$P_{z,(i-1)d+j}=z_i^{d-j} (X-z_i)^{j-1} 
\prod_{{\SS{s=1}} \atop {\SS{s \not= i}}}^k
(X-z_s)^d.$$
If $1 \le i \le k$, let $E_i$ denote the subspace of $\FM[X]/(P_z)$ spanned by the
$P_{z,(i-1)d+j}$, $1 \le j \le d$. Since $\deg P_{z,(i-1)d+j} < 
\deg P_{z,(i-1)d+j'} < \deg P_z$ whenever $1 \le j < j' \le d$, the family 
$(\overline{P}_{z,(i-1)d+j}^{\SSS{P_z}})_{1 \le j \le d}$ is a basis of $E_i$. 
Moreover, $E_i$ is contained 
in the characteristic subspace associated to the 
eigenvalue $z_i$ of the endomorphism $\r_z$. 
Since $z_i \not= z_{i'}$ if $1 \le i \not= i' \le k$, and for dimension reasons, 
we have
$$\FM[X]/(P_z)= E_1 \oplus \dots \oplus E_k,$$
$$\forall 1 \le i \le d,~\dim E_i =d,$$
and $(\overline{P}^{\SSS{P_z}}_{z,r})_{1 \le r \le n}$ is a basis of $\FM[X]/(P_z)$. 

It is immediately checked that the matrix of $\r_z$ in this basis is equal to 
$\db(z)v$. On the other hand, $M(P_z)$ is the matrix of $\r_z$ 
in the basis $\BC_{P_z}$. Therefore, if 
we denote by $\ph'(z)$ the transition matrix from the basis $\BC_{P_z}$ 
to the basis $(\overline{P_{z,1}}^{\SSS{P_z}},\dots,\overline{P_{z,n}}^{\SSS{P_z}})$, 
then we have 
$$\lexp{\ph'(z)}{(\db(z)v)}=M(P_z).$$
Moreover, $\ph' : \Ob \to \Gb\Lb_n(\FM)$ is a morphism of varieties, 
by definition of the polynomials $P_{z,r}$ ($1 \le r \le n$). 

Now let 
$$\d(z)=\prod_{1 \le i < j \le k} (z_i-z_j)$$
and let $\ph(z)=\ph'(z) \t(\d(z)^{-d})$. Since $\t(\d(z)^{-d})$ commutes with 
$\db(z)v$, we have 
$$\lexp{\ph(z)}{(\db(z)v)}=M(P_z).$$
Also, $\ph : \Zb(\Lb)_\reg^\ci \to \Gb\Lb_n(\FM)$ is a morphism of varieties. 
Hence, in order to prove Proposition \ref{connexite +}, we only need to prove that 
$\ph(z) \in \Gb$ for every $z \in \Ob$. This is equivalent to the following equality~:

\equat\label{determinant}
\det \ph'(z) = \d(z)^{d^2}.
\endequat

\med

\noindent{\sc Proof of \ref{determinant} - } If $1 \le i \le k$ and $1 \le j \le d$, 
we put
$$Q_{z,(i-1)d+j}=(\overline{X}^{\SSS{P_z}}-z_i)^{j-1} 
\prod_{\SS{s=1} \atop \SS{s \not= i}}^k (\overline{X}^{\SSS{P_z}}-z_s)^d=
z_i^{i-d} P_{z,(i-1)d+j}.$$
We denote by $G(z)$ the transition matrix 
from the basis $\BC_{P_z}$
to the basis $(Q_{z,1},\dots,Q_{z,n})$. 
Since $z_1 \dots z_k=1$, we have $\det G(z)=\det \ph'(z)$. 

Now, if $1 \le i \le k$ and $1 \le j \le d$, we put
$$R_{z,(i-1)d+j}=(\overline{X}^{\SSS{P_z}}-z_i)^{j-1} \prod_{s=1}^{i-1}
(\overline{X}^{\SSS{P_z}}-z_s)^d.$$
We denote by $H(z)$ the transition matrix 
from the basis $\BC_{P_z}$
to the basis $(R_{z,1},\dots,R_{z,n})$. 
It is obvious that $\det H(z)=1$ (indeed, $H(z)$ is unipotent upper triangular). 
Let $I(z)$ denote the transition matrix from the basis 
$(Q_{z,1},\dots,Q_{z,n})$ to the basis $(R_{z,1},\dots,R_{z,n})$. 
It is equal to $G(z)^{-1} H(z)$. Hence, the equality \ref{determinant} is equivalent to 
$$\det I(z) = \d(z)^{-d^2}.\leqno{\boldsymbol{(*)}}$$

Now, let us prove $(*)$. If $m$ is a natural number, we denote by $D_m : 
\FM[X] \to \FM[X]$ the \endo of the $\FM$-vector space $\FM[X]$ such that, 
for any $P \in \FM[X]$ and any $a \in \FM$, we have
$$P=\sum_{m \ge 0} D_m(P)(a) (X-a)^m.$$
(If $p=0$, then $D_m(P)=\DS{\frac{1}{m!}} P^{(m)}$ where $P^{(m)}$ denotes the $m$th 
derivative of $P$.) 
The operators $D_m$ $(m \in \NM$) satisfy the Leibniz rule~:
$$D_m(PQ)=\sum_{s=0}^m D_s(P)D_{m-s}(Q)$$
for all $m \in \NM$, and $P, Q \in \FM[X]$. 

If $1 \le i \le k$ and $1 \le j \le d$, we put
$$\fonction{C_{(i-1)d+j}}{\FM[X]}{\FM}{P}{D_{j-1}(P)(z_i).}$$
It is a linear form on $\FM[X]$ and it factorizes through $\FM[X]/(P_z)$. 
We still denote $C_{(i-1)d+j} : \FM[X]/(P_z) \to \FM$ the factorization of 
$C_{(i-1)d+j}$. Let $1 \le r \not= s \le n$. We have
$$C_r(Q_{z,s})=0.$$
Consequently, the basis $(C_r)_{1 \le r \le n}$ is, up to multiplication by scalars, 
a dual basis of $(Q_{z,r})_{1 \le r \le n}$. These constants are determined 
by the following identity~:
\equat\label{dual C}
C_{(i-1)d+j}(Q_{(i-1)d+j})=\prod_{\SS{s=1} \atop \SS{s \not= i}}^k (z_i-z_s)^d.
\endequat
Let $(C_r^\pr)_{1 \le r \le n+1}$ denote the dual basis to $(Q_{z,r})_{1 \le r 
\le n}$. We have $C_r^\pr(R_{z,s})=0$ if $1 \le s < r$. This proves that the matrix 
$I(z)$ is lower triangular. Hence the determinant of $I(z)$ is determined when 
one knows the diagonal coefficients. Using the identity \ref{dual C}, 
we obtain that, for all $1 \le i \le k$ and $1 \le j \le d$, we have
$$C_{(i-1)d+j}^\pr(R_{(i-1)d+j}) = {\DS{\prod_{s=1}^{i-1} (z_i-z_s)^d } 
\over \DS{\prod_{\SS{s=1} \atop \SS{s \not= i}}^k (z_i-z_s)^d}}
= {1 \over \DS{\prod_{s=i+1}^k (z_i-z_s)^d}}.$$
This gives 
$$\det I(z)^{-1} = \prod_{i=1}^k \prod_{j=1}^d \prod_{s=i+1}^k (z_i-z_s)^d = 
\d(z)^{d^2}$$
as desired. This completes the proofs of \ref{determinant}, 
and of Proposition \ref{connexite +}.\fin

\bi

If we mimic the proof of Theorem \ref{iso eta}, we get~:

\bi

\begin{coro}\label{pi-1}
The map
$$\fonction{\psi}{\Ob \times 
\mub_d}{\Ybt_\cdo^\pr}{(z,a)}{(z,\ph(z)\t(a) C_\Lb^\ci(v))}$$
is an isomorphism of varieties.
\end{coro}

\bi

\sub{Action of $\SG_k$\label{num w}} The group 
$W_\Gb(\Lb)=N_\Gb(\Lb)$ is canonically isomorphic to 
the symmetric group on $k$ letters $\SG_k$. If $w \in \SG_k$, we denote by $\wba$ 
the permutation matrix by blocks associated to $w$~: $\wba \in \Gb\Lb_n(\FM)$ 
and each block is either the zero $d \times d$ matrix either the $d \times d$ 
matrix $\unb_d$. Let
$$\wdo=\wba \t(\e(w))$$
where $\e : \SG_k \to \FM^\times$ is the signature of $\SG_k$. 
Then $\det \wdo = 1$, $\wdo$ is a representative of $w$ in $N_\Gb(\Lb)$, 
and, more precisely,
\equat\label{cgoo}
\wdo \in N_\Gb(\Lb) \cap C_\Gb^\ci(v).
\endequat
The group $\SG_k$ acts on $\Ob$ by permutation of the components~: this 
action corresponds, via the above identifications, to the natural action 
of $W_\Gb(\Lb)$ on $\Zb(\Lb)_\reg^\ci$ by conjugation. Moreover, if $z 
\in \Ob$ and if $w \in \SG_k$, then $\m(\db(z))=\m(\lexp{w}{\db(z)})$. It follows 
from Proposition \ref{connexite +} that 
$$\ph(z)^{-1}\ph(\lexp{w}{z}) \wdo \in C_\Gb(zv)=C_\Lb(v).$$
(Indeed, $C_\Gb(zv)=C_\Gb(z) \cap C_\Gb(v)$ and $C_\Gb(z)=C_\Gb^\ci(z)=\Lb$ 
since $z \in \Ob$ and $\Gb$ is semisimple and simply connected).

\bi

\begin{prop}\label{class}
If $z \in \Ob$ and $w \in \SG_k$, then the class 
of $\ph(z)^{-1}\ph(\lexp{w}{z}) \wdo$ in $A_\Lb(v)\simeq \ZC(\Lb)$ 
depends neither on $z$, nor on $\ph$, nor on the choice 
of a representative $\wdo$ of $w$ in $N_\Gb(\Lb) \cap C_\Gb^\ci(v)$. 
It depends only on $w$ 
and is equal to the class of $\t(\e(w)^{d-1})$.
\end{prop}

\bigskip

\proof The class of $\ph(z)^{-1}\ph(\lexp{w}{z}) \wdo$ in $A_\Lb(v)$
does not depend on $z$ because $\Ob$ is connected. 
If $\phi : \Ob \to \Gb$ is another morphism of 3 
such that $\lexp{\phi(z)}{(zv)}=M(P_z)$ for any $z 
\in \Ob$, then $\ph(z)^{-1}\phi(z) \in C_\Gb(zv)=C_\Lb(v)$. 
Again, the connexity of $\Ob$ implies that the class 
of $\ph(z)^{-1}\phi(z)$ in $A_\Lb(v)$ does not depend on $z$. 
Moreover, if $\wdo'$ is another representative of $w$ in $N_\Gb(\Lb) \cap 
C_\Gb^\ci(v)$, then $\wdo^{-1}\wdo' \in \Lb \cap C_\Gb^\ci(v)
=C_\Lb^\ci(v)$. This proves that the class 
of $\ph(z)^{-1}\ph(\lexp{w}{z}) \wdo$ in $A_\Lb(v) \simeq \ZC(\Lb)$ 
does not depend on the choice of $\wdo$. Let 
us denote it by $a$. If $\g(w)$ (\resp $\g'(w)$) denotes 
the class of $\ph(z)^{-1}\ph(\lexp{w}{z}) \wdo$ 
(\resp $\phi(z)^{-1}\phi(\lexp{w}{z}) \wdo$ in $A_\Lb(v)$), then 
$\g'(w)=a^{-1}\g(w)\lexp{w^{-1}}{a}=\g(w)$ because $\SG_k$ acts 
trivially on $A_\Lb(v) \simeq \ZC(\Lb)$. 

Now let us prove that $\g(w)$ is equal to the class of $\t(\e(w)^{d-1})$. 
If $z \in \Ob$ and if $w \in \SG_k$, then $P_{\lexp{w}{z}}=P_z$. 
The matrix $\ph(z)^{-1}\ph(\lexp{w}{z})$ is the transition matrix from the basis 
$$(P_{z,1},\dots,P_{z,n-d},\d(z)^{-d}P_{z,n+1-d},\dots,\d(z)^{-d}P_{z,n})$$ 
to the basis
$$(P_{\lexp{w}{z},1},\dots,P_{\lexp{w}{z},n-d},
\d(\lexp{w}{z})^{-d}P_{\lexp{w}{z},n+1-d},\dots,\d(\lexp{w}{z})^{-d}
P_{\lexp{w}{z},n}).$$ 
Hence it is equal to
$$\t(\d(z)^d) \wba^{-1} \t(\d(\lexp{w}{z})^{-d}).$$
But $\d(\lexp{w}{z})=\e(w)\d(z)$. In particular,
$$\ph(z)^{-1}\ph(\lexp{w}{z})\wdo = 
\t(\d(z)^d) \lexp{\wba^{-1}}{\t(\d(z)^{-d}\e(w)^{-d})} \t(\e(w)).$$
Since $\SG_k$ acts trivially on $\ZC(\Lb)$, 
the class of $\ph(z)^{-1}\ph(\lexp{w}{z})\wdo$ 
in $A_\Lb(v)$ is equal to the class of 
$$\t(\d(z)^d) \t(\d(z)^{-d}\e(w)^{-d}) \t(\e(w)) = \t(\e(w)^{1-d}).$$
But $\e(w)^{1-d}=\e(w)^{d-1}$ because $\e(w) \in \{1,-1\}$. 
The Proposition \ref{class} follows.\fin

\bi

\begin{coro}\label{enfin sln}
We have, for every $w \in \SG_k$, $\ph_{\Lb,v}^\Gb(w)=\t(\e^{d-1}(w))$, as expected.
\end{coro}

\bi

\proof Let $\WC=\{(w,a) \in W_\Gb^\ci(\Lb,v) \times A_\Lb(v)~|~a=\t(\e^{d-1}(w))\}$. 
Then, by Proposition \ref{class}, $\WC$ is the stabilizer of any irreducible 
component of $\Ybt_\cdo^\pr$. Since $\Xbt_\cdo^\pr$ is smooth (\cf Theorem 
\ref{iso eta}), then any connected component of $\Xbt_\cdo^\pr$ is 
irreducible, and is in fact the closure of an irreducible component 
of $\Ybt_\cdo^\pr$. Therefore, the stabilizer of any element of 
$\Xbt_\cdo^\pr$ must be contained in $\WC$. In particular, 
$H_\Gb(\Lb,v,1) \incl \WC$. But $|H_\Gb(\Lb,v,1)|=|\WC|$ by Proposition \ref{agu} (2), 
so the corollary follows.\fin

\bi

\sec{Small classical groups\label{cas B}}~

\med

\sub{Reductions\label{sub pas mal}} 
We assume in this subsection that $\Gb$ is semisimple, 
simply connected, quasi simple, that $\Lb$ is a \levi of a maximal 
\para of $\Gb$, and that $\Lb$ is cuspidal. We also assume that 
$p \not= 2$, that $p$ is good for $\Gb$, and that the family 
$(\sdo_\a,x_\a,x_{-\a})_{\a \in \D}$ is chosen in such a way that 
Corollary \ref{inverse omega} holds.
Here is a list of consequences of this hypothesis~:

\bi

\begin{lem}\label{liste} 
Under the above hypothesis, we have~:

\tete{1} $|W_\Gb(\Lb)|=2$.

\tete{2} $\dim \Zb(\Lb)^\ci=1$.

\tete{3} If $s$ is the non-trivial element of $W_\Gb(\Lb)$ and if 
$z \in \Zb(\Lb)^\ci$, then $\lexp{s}{z}=z^{-1}$. 

\tete{4} The closed subset $\Zb(\Lb)^\ci-\Zb(\Lb)^\ci_\reg$ of $\Zb(\Lb)^\ci$ 
consists of elements of order $1$ or $2$. 
Consequently, $|\Zb(\Lb)^\ci-\Zb(\Lb)^\ci_\reg|=2$.

\tete{5} Let $z \in \Zb(\Lb)^\ci-\Zb(\Lb)^\ci_\reg$ and $x \in \Gb$ be such 
that $x^{-1}zx \in \Zb(\Lb)^\ci$. Then $x \in C_\Gb(z)=C_\Gb^\ci(z)$.
\end{lem}

\proof (1) follows from Proposition \ref{propriete distingues} (b) and Proposition 
\ref{cuspidal isotypique} (b). (2) follows 
from the semisimplicity of $\Gb$ and the maximality of $\Lb$. 
Since $C_\Gb(\Zb(\Lb)^\ci)=\Lb$, the element $s \in W_\Gb(\Lb)$ 
induces a non-trivial automorphism of $\Zb(\Lb)^\ci$. 
Hence (3) follows from (2) and from the fact that a one-dimensional torus 
has a unique non-trivial automorphism, namely the one sending an element 
to its inverse. 

Now, let $z \in \Zb(\Lb)^\ci$ be an element of order $1$ or $2$. By (3), 
$\sdo \in C_\Gb(z)$, where $\sdo$ is a representative of $s$ in 
$N_\Gb(\Lb)$. But, since $\Gb$ is semisimple and simply-connected, 
$C_\Gb(z)$ is connected \cite[Theorem 8.2]{steinendo}. So $C_\Gb^\ci(z) 
=C_\Gb(z) \not= \Lb$. Therefore, $z \not\in \Zb(\Lb)_\reg^\ci$. 

Conversely, assume that $z \in \Zb(\Lb)^\ci-\Zb(\Lb)_\reg^\ci$. 
Then $\Gb'=C_\Gb^\ci(z)=C_\Gb(z) \not= \Lb$. So, since $\Lb$ is cuspidal, 
$|W_{\Gb'}(\Lb)| = 2$. Hence, it follows 
from (1) that $\sdo \in \Gb'$, that is that $z=z^{-1}$. In other words, 
$z$ is of order $1$ or $2$. The proof of (4) is complete. 

Let $z \in \Zb(\Lb)^\ci-\Zb(\Lb)^\ci_\reg$ and $x \in \Gb$ be such 
that $x^{-1}zx \in \Zb(\Lb)^\ci$. There is nothing to prove if 
$z$ is of order $1$ and, if $z$ is of order $2$, then, by (2), $z$ is the unique 
element of order $2$ of $\Zb(\Lb)^\ci$. So $x \in C_\Gb(z)=C_\Gb^\ci(z)$.\fin

\bi

From now on, we denote by $z$ the unique element of $\Zb(\Lb)^\ci$ of order $2$, 
and we set $\Gb'=C_\Gb(z)$. Then, by Lemma \ref{liste} (5), 
$$W_{\Gb'}(\Lb,v)=W_\Gb(\Lb,v).$$
We denote by $s$ the unique non-trivial element of 
$W_\Gb^\ci(\Lb,v)=W_{\Gb'}^\circ(\Lb,v)$. We 
assume that $\ph_{\Lb,v}^{\Gb'}(s)=1$, which means that $s$ stabilizes 
the element $\uti_z^\pr$. 

\bi

\rem In cases $B$, $C$, $D$ and $E$, we are under these hypothesis. 
Indeed, this may be checked case-by-case using \ref{produit H}.\finl

\bi

Let $\Xb^\ci$ be the irreducible component of $\Xbt_\cdo^\pr$ containing 
$\uti_z^\pr$ (or, equivalently, $\uti$). Since $\ph_{\Lb,v}^{\Gb'}(s)=1$, 
$s$ stabilizes $\Xb^\ci$. 
Let $A$ be the stabilizer, in $A_\Lb(v)$, of $\Xb^\ci$. Then the stabilizer 
of $\uti'$ is contained in $W_\Gb^\ci(\Lb,v) \times A$. 

Moreover, by Theorem \ref{iso eta}, $\Xb^\ci$ is isomorphic to $\FM^\times$, 
and $\Xbt_\cdo$ is also isomorphic to $\FM^\times$. 
We will denote by $\FM^\times_*$ for the variety $\FM^\times$ identified 
with $\Xb^\ci$, and simply $\FM^\times$ for the same variety identified 
with $\Xbt_\cdo$. With these identifications, 
the restriction of $\fti_\cdo$ 
to $\Xb^\ci \to \Xbt_\cdo$ may be  identified with the map 
$f^\ci : \FM^\times_* \to \FM^\times$, $x \mapsto x^n$, where $n=|A|$, and $A$ 
may be identified with $\mub_n$, acting on $\FM^\times_*$ by translation. 
The fact that $f^\ci$ is $W_\Gb^\ci(\Lb,v)$-equivariant 
shows that the action of $s$ on $\FM_*^\times$ is given by 
$$\forall x \in \FM^\times_*,~\lexp{s}{x}=\kappa x^{-1},$$
where $\kappa$ is some element of $\mub_n$. Since the action of 
$A$ commutes with the action of $s$, we deduce that $n=1$ or $2$, so 
that $\kappa=\pm 1$. 

Since $\uti'$ may be identified with $1 \in \FM_*^\times$, we have 
\equat\label{kap}
\kappa=\ph_{\Lb,v}^\Gb(s), 
\endequat
where $\mub_n$ is identified with $A$ which is a subgroup of $A_\Lb(v)$. 

\bi

\begin{lem}\label{simplification}
With the above notation, $\ph_{\Lb,v}^\Gb(s)=(-1)^{n-1}$, with $n$ equal to 
$|A_\Lb(v)|/i$ where $i$ is the number of irreducible components of $\Xbt_\cdo^\pr$.
\end{lem}

\bi

\proof Indeed, if $n=1$, then $\kappa \in\mub_1$, so $\kappa=1$. 
On the other hand, if $n=2$, then $\uti_z^\pr$ may be identified 
with a preimage of $-1$ in $\FM^\times_*$, that is to a fourth 
root of unity, say $i$. But, since $\lexp{s}{i}=i$, we get that 
$\kappa=-1$. The Lemma then follows from Formula \ref{kap}.\fin

\bi

Lemma \ref{simplification} gives 
an easy criterion for computing $\ph_{\Lb,v}^\Gb(s)$. Indeed, 
$n$ is equal to the order of $A_\Lb(v)\simeq \ZC(\Lb)$ 
divided by the number of irreducible components of $\Xbt_\cdo^\pr$. 
This number can be determined using Theorem \ref{iso eta} and easy 
computations in the root system. 
This fact will be used in the next subsections for completing the proof 
of Table \ref{tabletable}.

\bigskip

\sub{Case B\label{sub 1}} {\it Until the end of \SEC\ref{sub 1}, we assume that $\Gb 
\simeq \Sb\pb_4(\FM)$, that $\Lb \simeq \Sb\Lb_2(\FM) 
\times \FM^\times$ and that $p \not=2$.}  So 
$\Gb$ is semisimple and simply connected of type $C_2$. 

We identify 
$X(\Tb)$ with $\ZM^2$ and we denote by $(\e_1,\e_2)$ the canonical 
basis of $X(\Tb)=\ZM^2$. 
Let $\a_1=\e_1-\e_2$ and $\a_2=2\e_2$. We identify $\D$ with $\{\a_1,\a_2\}$, 
so that $\D_\Lb=\{\a_2\}$. Note also that $w_0=-\Id$. Let $(\e_1^*,\e_2^*)$ 
denote the basis of $Y(\Tb)$ dual to $(\e_1,\e_2)$. Then $\a_1^\ve=\e_1^*-\e_2^*$ 
and $\a_2^\ve=\e_2^*$. Hence $(\varpi_{\a_1},\varpi_{\a_2}) = 
(\e_1,\e_1+\e_2)$. Then $Y(\Zb(\Lb)^\ci)$ is spanned by $\e_1^*$. 

We denote by $\Sb$ the subtorus of $\Tb$ such that $Y(\Sb)=<\e_2^*>$. 
We have $\Tb=\Zb(\Lb)^\ci \times \Sb$. By \ref{defi chi}, we have
$$\ch_1=\varpi_{\a_1}-\langle \a_1,\a_2^\ve \rangle \varpi_{\a_2} = 2\e_1+\e_2,$$
$$\ch_2=\varpi_{\a_2}=\e_1+\e_2.\leqno{\mathrm{and}}$$

By Theorem \ref{iso eta}, we have
\eqna
\Xbt_\cdo^\pr &\simeq& \{(z,t) \in \Zb(\Lb)^\ci \times \Sb~|~
~\ch_2(z)\a_2(t) = 1\} \\
&\simeq& \{(a,b) \in \FM^\times \times \FM^\times~|~ab^2=1 \} \\
&\simeq& \FM^\times.
\endeqna
Therefore, $\Xbt_\cdo^\pr$ is connected. Since $|A_\Lb(v)|=2$, we get 
from Lemma \ref{simplification} 
that $\ph_{\Lb,v}^\Gb(s)=-1$, as it is expected from Table \ref{tabletable}.

\bi

\remark{cgzcgz} By the above computation, we have 
$$\ph_{\Lb,v}^{\Gb'} \not= 
\Res_{W_{\Gb'}^\ci(\Lb,v)}^{W_\Gb^\ci(\Lb,v)} \ph_{\Lb,v}^\Gb$$
(see Remark \ref{non restriction}).\finl

\bi

\rem The value of $\ph_{\Lb,v}^\Gb(s)$ in this case  
has been obtained using a different method by J.L. Waldspurger 
\cite[Lemma VIII.9]{waldspurger}. See also Example \ref{wa}.\finl

\bi

\sub{Case C\label{sub 2}} {\it Until the end of \SEC\ref{sub 2}, we assume that 
$\Gb \simeq \Sb\pb\ib\nb_7(\FM)$, that $\Lb$ is of type 
$A_1 \times A_1$, and that $p \not=2$.} So 
$\Gb$ is semisimple and simply connected of type $B_3$. 

We assume that the basis $\D=\{\a_1,\a_2,\a_3\}$ of $\Phi$ is 
numbered in such a way that its Cartan matrix is given by
$$(\langle \a_i , \a_j^\ve\rangle)_{1 \le i,j \le 3} = 
\matrice{2 & -1 & 0 \\ -1 & 2 & -2 \\ 0 & -1 & 2}.$$
In particular, $\D_\Lb=\{\a_1,\a_3\}$. 
Note that $w_0=-\Id$. 
In this case, we get 
\eqna
\varpi_{\a_1}&=&\a_1+\a_2+\a_3,\\ 
\varpi_{\a_2}&=& \a_1 + 2 \a_2 +2\a_3,\\
\varpi_{\a_3} &=& \DS{1 \over 2}(\a_1 + 2\a_2 + 3 \a_3).
\endeqna 
As a consequence, it follows from \ref{defi chi} that 
\eqna
\ch_1&=&2\a_1+3\a_2+3\a_3, \\
\ch_2&=& 2 \a_1 + 4\a_2 + 5 \a_3, \\
\ch_3&=&\DS{1 \over 2}(\a_1 +2\a_2+3\a_3).
\endeqna
Now, let $y=\a_1^\ve+2\a_2^\ve + \a_3^\ve$. Then $Y(\Zb(\Lb)^\ci)=<y>$. 
Let $\Sb$ be the two-dimensional subtorus of $\Tb$ such that 
$Y(\Sb)=<\a_1^\ve,\a_2^\ve >$. Then $\Tb=\Zb(\Lb)^\ci \times \Sb$. 
By Proposition \ref{iso eta}, we have 
\eqna
\Xbt_\cdo^\pr &\simeq& 
\{(z,t) \in \Zb(\Lb)^\ci \times \Sb~|~\ch_1(z)\a_1(t) = 1 {\mathrm{~and~}} 
\ch_3(z)\a_3(t)=1\} \\
&\simeq& \{(a,b,c) \in \FM^\times \times \FM^\times 
\times \FM^\times~|~\ch_1(y(a))\a_1(\a_1^\ve(b)\a_2^\ve(c))=1\\
&& \hphantom{\{(a,b,c) \in \FM^\times \times \FM^\times 
\times \FM^\times~|}
{\mathrm{~and~}} \ch_3(y(a))\a_3(\a_1^\ve(b)\a_2^\ve(c))=1\} \\
&\simeq& \{(a,b,c) \in \FM^\times \times \FM^\times 
\times \FM^\times ~|~a^3b^2c^{-1}=1
{\mathrm{~and~}} ac^{-1}=1\} \\
&\simeq& \{(a,b) \in \FM^\times \times \FM^\times~|~a^2b^2=1\} \\
&\simeq& \FM^\times \times \mub_2.
\endeqna
Then, by \ref{kap} and Lemma \ref{simplification}, we get that $\ph_{\Lb,v}^\Gb(s)=1$. 

\bi

\sub{Case D\label{sub 3}} {\it Until the end of \SEC\ref{sub 3}, we assume that $\Gb 
\simeq \Sb\pb\ib\nb_{10}(\FM)$, that $\Lb$ is of type $A_1 \times 
A_3$, and that $p\not= 2$.} So 
$\Gb$ is semisimple and simply connected of type $B_3$. 

We assume that the basis $\D=\{\a_1,\a_2,\a_3,\a_4,\a_5\}$ of $\Phi$ is 
numbered in such a way that its Cartan matrix is given by
$$(\langle \a_i , \a_j^\ve\rangle)_{1 \le i,j \le 5} = 
\matrice{2 & -1 & 0 & 0 & 0 \\ -1 & 2 & -1 & 0 & 0 \\
0 & -1 & 2 & -1 & -1 \\ 0 & 0 & -1 & 2 & 0 \\ 0 & 0 & -1 & 0 & 2}.$$
In particular, $\D_\Lb=\{\a_1,\a_3,\a_4,\a_5\}$. 
Note that $-w_0$ permutes the simple roots $\a_4$ and $\a_5$ and stabilizes 
the other simple roots. So $(i_1,i_2,i_3,i_4)=(1,3,4,5)$. 
We get
\eqna
\ch_1&=&2\a_1+3\a_2+3\a_3+\DS{3 \over 2}(\a_4+\a_5), \\
\ch_2&=&2\a_1+4\a_2+5\a_3+\DS{5 \over 2}(\a_4+\a_5), \\
\ch_3&=&2\a_1+4\a_2+6\a_3+\DS{7 \over 2}(\a_4+\a_5), \\
\ch_4&=&\DS{1 \over 2}(-\a_1-2\a_2-3\a_3-\DS{1 \over 2}\a_4-\DS{7 \over 2}\a_5), \\
\ch_5&=&\DS{1 \over 2}(\a_1+2\a_2+3\a_3+\DS{3 \over 2} \a_4 + \DS{5 \over 2}\a_5). \\
\endeqna
Now let $y=\a_1^\ve + 2 \a_2^\ve + 2 \a_3^\ve +\a_4^\ve+\a_5^\ve$. Then
$Y(\Zb(\Lb)^\ci)=<y>$. Let $\Sb$ be the four-dimensional subtorus 
of $\Tb$ such that $Y(\Sb)=<\a_2^\ve,\a_3^\ve,\a_4^\ve,\a_5^\ve>$. 
By Theorem \ref{iso eta}, we have
$$\Xbt_\cdo^\pr =\{(z,t) \in \Zb(\Lb)^\ci \times \Sb~|~\left\{\begin{array}{lcl}
\ch_1(z)\a_1(t) &=& 1 \\
\ch_3(z)\a_3(t) &=& 1 \\
\ch_4(z)\a_5(t) &=& 1 \\
\ch_5(z)\a_4(t) &=& 1 \\
\end{array}\right.\}.$$
If we identify $\Zb(\Lb)^\ci$ to $\FM^\times$ via $y$, and 
$\Sb$ to $(\FM^\times)^4$ via the isomorphism 
$$\fonctio{(\FM^\times)^4}{\Sb}{(b,c,d,e)}{\a_2^\ve(b)
\a_3^\ve(c)\a_4^\ve(d)\a_5^\ve(e),}$$
we get that
\begin{eqnarray*}
\Xbt_\cdo^\pr &=&\{(a,b,c,d,e) \in 
(\FM^\times)^5~|~\left\{\begin{array}{lcl}
a^3b^{-1}&=& 1 \\
a^4b^{-1}c^2d^{-1}e^{-1} &=& 1 \\
a^{-1}c^{-1}e^2 &=& 1 \\
ac^{-1}d^2 &=& 1 
\end{array}\right.\}\\
&=&\{(a,c,d,e) \in 
(\FM^\times)^4~|~\left\{\begin{array}{lcl}
ac^2d^{-1}e^{-1} &=& 1 \\
a^{-1}c^{-1}e^2 &=& 1 \\
ac^{-1}d^2 &=& 1 
\end{array}\right.\} \\
&=&\{(c,d,e) \in 
(\FM^\times)^3~|~\left\{\begin{array}{lcl}
cd^{-1}e &=& 1 \\
c^{-3}d^3e &=& 1 
\end{array}\right.\}\\
&=&\{(d,e) \in 
(\FM^\times)^2~|~
e^4= 1 \}.
\end{eqnarray*}
Again, using Lemma \ref{simplification}, we get that 
$\ph_{\Lb,v}^\Gb(s)=1$, as expected.

\bi

\sub{Case E} This case is the easiest. Indeed, the invariance of the problem 
under the order $3$ automorphism of $\Gb$ yields immediately the answer, 
because $1$ is the unique invariant element in $\ZC(\Lb)$. 

\bi

\rem The courageous reader can try to use the previous methods 
to find another proof of this last result. This is of course possible, 
as the author has checked by himself.\finl

\bi

The proof of the results stated in Table \ref{tabletable} 
is now complete.

\bigskip

\section{Lusztig restriction of characteristic functions of regular unipotent 
classes\label{section fini}}~

\med

From now on, and until the end of this paper, we assume that we are 
in the situation of Section \ref{part finite}, that is, we assume 
that $\Gb$ is defined over a finite field. Recall that we also assume that 
$v$ is a regular unipotent element. Finally, 
we assume that $\Cb$ supports a cuspidal local system, and that 
$p$ is good for $\Gb$. Consequently we can use all the results proved 
in the first part and in the previous sections of this second part.

\medskip

The aim of this section is to compute the Lusztig restriction $\lexp{*}{R}_\Mb^\Gb$ 
of the characteristic function of the $\Gb^F$-conjugacy class 
of $u$ (which is a regular unipotent element of $\Gb$). Here, $\Mb$ 
is an $F$-stable \levi of a \para of $\Gb$.  
It is expected to be the characteristic function of the conjugacy class of 
some regular unipotent element $u' \in \Mb^F$~: this result has been 
proved whenever $p$ is good and $q$ is large enough for $\Gb$.
%$(\Gb,F)$ is friendly \cite{DLM2} (\cf Remark 
%\ref{large enough} for the definition of friendly). 
%In this 
In this section, we determine explicitly the $\Gb^F$-conjugacy class 
of $u'$~: this precision was not in \cite{DLM2}. However, as for Digne-Lehrer and 
Michel's theorem, 
our result is proved only whenever $p$ is good and $q$ is large enough. 

This section is organized as follows. In \SEC\ref{res map}, we recall 
the definition of  
a map between the set of $\Gb^F$-conjugacy classes of regular unipotent 
elements in $\Gb^F$ to the corresponding set in $\Mb^F$. This map was 
constructed in \cite[page 279]{Bonnafe}. We also provide a proof for 
some properties of this map which were announced in 
\cite[Proposition 2.2]{Bonnafe}. In \SEC\ref{gelfand}, we 
prove the main result, that is the computation of the Lusztig 
restriction of the characteristic function of the $\Gb^F$-conjugacy 
class of $u$.

\bi

\noindent{\bf Hypothesis and notation :} {\it We assume all along this part 
that we are in the situation of Section \ref{part finite}. We also assume 
that $p$ is good for $\Gb$. Finally, we fix an $F$-stable \levi $\Mb$ 
of a \para $\Qb$ of $\Gb$.}

\bi

As a consequence of these hypothesis, we may (and we will) 
identify $A_\Gb(u)$ with $\ZC(\Gb)$ via the canonical isomorphism. 

\bi

\subsection{Restriction map for regular unipotent classes\label{res map}}~
Let $\Reg_\uni(\Gb^F)$ denote the set of regular 
unipotent classes of $\Gb^F$. We recall in this subsection the construction 
of a map $\res_\Mb^\Gb : \Reg_\uni(\Gb^F) \to \Reg_\uni(\Mb^F)$ which was 
made in \cite[page 279]{Bonnafe}. For a proof of all the results 
used here, the reader may refer to \cite{Bonnafe}. 

Let us first consider the easiest case~: if $\Qb$ is $F$-stable, we set 
$$\fonction{\r_\Mb^\Gb}{\Reg_\uni(\Gb^F)}{\Reg_\uni(\Mb^F)}{[g]_{\Gb^F}}{
\pi_\Mb([g]_{\Gb^F} \cap \Qb^F).}$$
It is well-known that $\r_\Mb^\Gb$ is well-defined and does not 
depend on the choice of $\Qb$.

Now, let us come back to the general case. 
We fix an $F$-stable \borel $\Bb_\Mb$ of $\Mb$, and an $F$-stable 
\tor $\Tb_\Mb$ of $\Bb_\Mb$. We denote by $\Lb_\Mb$ the unique minimal standard 
\levi of a \para $\Pb_\Mb$ of $\Mb$ (standard with respect to the pair 
$(\Tb_\Mb,\Bb_\Mb)$) such that the canonical map 
$\ZC(\Mb) \to \ZC(\Lb_\Mb)$ is an isomorphism. 
By unicity, $\Lb_\Mb$ and $\Pb_\Mb$ are $F$-stable. 
Moreover, the map $\r_{\Lb_\Mb}^\Mb$ is bijective.
 
Now, by minimality, $\Lb_\Mb$ is quasi-cuspidal, so it is 
universally distinguished. Therefore, 
there exists a unique standard \levi $\Lb_\Mb^\pr$ of a \para 
$\Pb_\Mb^\pr$ of $\Gb$ (standard with 
respect to $(\Tb,\Bb)$) which is $\Gb$-conjugate to $\Lb_\Mb$. 
Since $\Lb_\Mb$ is $F$-stable, $\Lb_\Mb^\pr$ and $\Pb_\Mb^\pr$ are $F$-stable. 
Moreover, by Proposition \ref{conjugue M} and Proposition 
\ref{cuspidal isotypique} (d), the map 
$$\fonction{c}{\Reg_\uni(\Lb_\Mb^{\pr F})}{\Reg_\uni(\Lb_\Mb^F)}{
[l]_{\Lb_\Mb^{\pr F}}}{[l]_{\Gb^F} \cap \Lb_\Mb^F}$$
is bijective. 

Using the above notation and properties, we define the map $\res_\Mb^\Gb$ 
to be the composition
$$\diagram
\Reg_\uni(\Gb^F) \rrto^{\DS{\r_{\Lb_\Mb^\pr}^\Gb}}&& \Reg_\uni(\Lb_\Mb^{\pr F}) 
\rrto^{\DS{c}} && \Reg_\uni(\Lb_\Mb^F) \rrto^{\DS{(\r_{\Lb_\Mb}^\Mb)^{-1}}} 
&& \Reg_\uni(\Mb^F).\enddiagram$$

\bi

We recall some properties of these restriction maps (proofs of these 
fact may be found in \cite[\SEC 7]{bonnafe regular}). First, 
\equat\label{11}
\res_\Mb^\Gb~{\mathit{is~independent~of~the~choice~of~}}\Qb.
\endequat 
Moreover, if 
$\Mb$ is an $F$-stable \levi of an $F$-stable \para of $\Gb$, then
\equat\label{22}
\res_\Mb^\Gb = \r_\Mb^\Gb.
\endequat

To state the next property, we need to introduce some further notation. 
First, we denote by $\hb_\Mb^\Gb : H^1(F,\ZC(\Gb)) \to H^1(F,\ZC(\Mb))$ 
the map induced by $h_\Mb^\Gb$ (it is surjective). Now, 
let $\UC \in \Reg_\uni(\Gb^F)$ and let $z \in H^1(F,\ZC(\Gb))$. We denote 
by $g_z$ an element of $\Gb$ such that $g_z^{-1}F(g_z) \in \Zb(\Gb)$ 
and represents $z$. We then define $\UC_z=\lexp{g_z}{\UC}$. 
Then $\UC_z \in \Reg_\uni(\Gb^F)$, and $\UC_z$ depends only on $z$, and 
not on the choice of the element $g_z$. Moreover, since 
$A_\Gb(u)\simeq \ZC(\Gb)$, the map
$$\fonctio{H^1(F,\ZC(\Gb))}{\Reg_\uni(\Gb^F)}{z}{\UC_z}$$ 
is bijective. Also, since $g_z$ may be chosen in $\Mb$, it is immediate 
that
\equat\label{33}
\res_\Mb^\Gb \UC_z=(\res_\Mb^\Gb \UC)_{\hb_\Mb^\Gb(z)}.
\endequat
The last result is a transitivity property \cite[Proposition 7.2 (c)]{bonnafe regular}~:
if $\Mb'$ is an $F$-stable \levi of a \para of $\Gb$ containing 
$\Mb$, then 
\equat\label{44}
\res_\Mb^\Gb=\res_\Mb^{\Mb'} \ci \res_{\Mb'}^\Gb.
\endequat 

\bi

\subsection{Lusztig restriction to $\Lb_w$\label{gelfand}}~
If $w \in W_\Gb^\ci(\Lb,v)$ and if $z \in H^1(F,\ZC(\Lb))$, then 
\equat\label{155}
|C_{\Lb_w}^\ci(v_{w,z})^F|=|\Zb(\Lb)^{\ci wF}| ~q^{r'}.
\endequat

\bigskip

\noindent{\sc Proof of \ref{155} - } Indeed, 
\eqna
|C_{\Lb_w}^\ci(v_{w,z})^F|&=&|C_\Lb(v_z)^{\wdo F}| \\
&=& |\Zb(\Lb)^{\ci \zdo wF}|.|C_{\Ub_\Lb}(v)^{\zdo\wdo F}|.
\endeqna
But $C_{\Ub_\Lb}(v)$ is a connected unipotent group of dimension $r'$, so for every 
Frobenius morphism $F'$ defining an $\fq$-structure 
on $C_{\Ub_\Lb}(v)$, we have $|C_{\Ub_\Lb}(v)^{F'}|=q^{r'}$. The result follows from 
the fact that $\zdo$ acts trivially on $\Zb(\Lb)^\ci$.\fin 

\bigskip

Moreover, since 
$W_\Gb(\Lb)$ acts trivially on $\ZC(\Lb)$, we have
\equat
|C_{\Lb_w}(v_{w,z})^F|=|\Zb(\Lb)^{\ci wF}| ~q^{r'}. |\ZC(\Lb)^F|.
\endequat

\bi

If $\x$ is a linear character of $\ZC(\Gb)$, we set~:
$$\gamh_{u,\x}^\Gb={1 \over |H^1(F,\ZC(\Gb))|} \sum_{z \in H^1(F,\ZC(\Gb))} 
\x(z) \g_{u_z}^\Gb.$$
By \ref{je sais pas}, the characteristic function $\ch_{\EC_w,\ph_w}$ is a multiple 
of $\gamh_{v_w,\z}^{\Lb_w}$. The isomorphism $\ph$ is determined up to multiplication 
by a scalar, we may, and we will, choose $\ph$ such that
\equat
\ch_{\EC,\ph}={1 \over |\Zb(\Lb)^{\ci F}|}~\gamh_{v,\z}^\Lb.
\endequat
In particular, it follows from Corollaries \ref{coro zeta phi} 
and \ref{description chi} that
\equat\label{tiens tiens}
\ch_{\EC_w,\ph_w}={\z(\ph_{\Lb,v}^\Gb(w)) 
\over |\Zb(\Lb_w)^{\ci F}|}\gamh_{v_w,\z}^{\Lb_w}.
\endequat

Our aim in this subsection 
is to determine the Lusztig restriction 
to $\Lb_w^F$ of the function $\gamh_{u,\zett}^\Gb$, where $\zett=\z \ci \hb_\Lb^\Gb$. 
For this purpose, we first compute its Harish-Chandra restriction 
to $\Lb$ using elementary methods. Then we use all the work done during 
this paper to deduce the general case. However, our result is known 
to be valid only if $q$ is large enough. 
By \cite[Proposition 5.3]{DLM1}, we have, for any $z \in H^1(F,\ZC(\Gb))$,
$$\lexp{*}{R}_{\Lb \incl \Pb}^\Gb(\g_{u_z}^\Gb)=\g_{v_{\hb_\Lb^\Gb(z)}}^\Lb.$$
Hence, we get
\equat\label{res gamma}
\lexp{*}{R}_{\Lb \incl \Pb}^\Gb(\gamh_{u,\zett}^\Gb)=\gamh_{v,\z}^\Lb.
\endequat

\begin{theo}\label{lusztig restriction}
Let $w \in W_\Gb^\ci(\Lb,v)$. Assume that $p$ is good for $\Gb$ and that 
$q$ is large enough. Then
$$\lexp{*}{R}_{\Lb_w}^\Gb(\gamh_{u,\zett}^\Gb)=\z(\ph_{\Lb,v}^\Gb(w))
\gamh_{v_w,\z}^{\Lb_w}.$$
\end{theo}

\bigskip

\proof First, we determine the following scalar product
$$\kappa = \langle \XC_{K_1,\ph_1}, \gamh_{u,\zett}^\Gb \rangle_{\Gb^F}.$$
If $\ch$ is a non-trivial irreducible character of $W_\Gb^\ci(\Lb,v)$, 
then $\langle \XC_{K_\ch,\ph_\ch}, \gamh_{u,\zett}^\Gb \rangle_{\Gb^F}=0$ 
since the perverse sheaf $K_\ch$ does not contain $\Cb^\Gb=(u)_\Gb$ in its support. 
Hence, for each $w \in W_\Gb^\ci(\Lb,v)$, we have 
$$\kappa = \langle \sum_{\ch \in (\Irr W_\Gb^\ci(\Lb,v))^F} \chit(wF) 
\XC_{K_\ch,\ph_\ch},\gamh_{u,\zett}^\Gb \rangle_{\Gb^F}.$$
But, by Theorem \ref{formule lusztig}, we have
$$\sum_{\ch \in (\Irr W_\Gb^\ci(\Lb,v))^F} \chit(wF) \XC_{K_\ch,\ph_\ch}
= R_{\Lb_w}^\Gb(\XC_{\EC_w,\ph_w}).$$
From this and from \ref{tiens tiens}, it follows that 
$$|\Zb(\Lb)^{\ci wF}|
\kappa=\z(\ph_{\Lb,v}^\Gb(w))\langle \gamh_{v_w,\z}^{\Lb_w},\lexp{*}{R}_{\Lb_w}^\Gb
(\gamh_{u,\zett}^\Gb) \rangle_{\Gb^F}\leqno{(*)}$$
for every $w \in W_\Gb^\ci(\Lb,v)$. 

Let us now use $(*)$ in the case where $w=1$. By \ref{res gamma}, we obtain 
the following equality~:
$$|\Zb(\Lb)^{\ci F}|
\kappa=\langle \gamh_{v,\z}^\Lb , \gamh_{v,\z}^\Lb \rangle_{\Gb^F}=|C_\Lb^\ci(v)^F|.$$
So $\kappa=q^{r'}$. 
Therefore, we deduce from $(*)$ and this equality that 
$$\langle \gamh_{v_w,\z}^{\Lb_w},\lexp{*}{R}_{\Lb_w}^\Gb
(\gamh_{u,\zett}^\Gb) \rangle_{\Gb^F}=
\z(\ph_{\Lb,v}^\Gb(w))|\Zb(\Lb_w)^{\ci F}|~q^{r'}.$$
In other words,
\equat\label{youhou}
\langle \gamh_{v_w,\z}^{\Lb_w},\lexp{*}{R}_{\Lb_w}^\Gb
(\gamh_{u,\zett}^\Gb) \rangle_{\Gb^F}=
\z(\ph_{\Lb,v}^\Gb(w)) |C_{\Lb_w}^\ci(v_w)^F|.
\endequat
Since $p$ is good for $\Lb$ and since the regular unipotent class of $\Lb$ 
supports a cuspidal local system, then all the irreducible components 
of the root system $\Phi_\Lb$ are of type $A$ \cite[Corollary 1.1.4]{bonnafe torsion}. 
But, by \cite[Lemma 1.8.4 (b) and Corollary 6.2.2 (a)]{bonnafe A} and 
\cite{bonnafe A'}, 
$\lexp{*}{R}_{\Lb_w}^\Gb(\gamh_{u,\zett}^\Gb)$ is a 
multiple of $\gamh_{v_w,\z}^{\Lb_w}$. This completes the proof of the Theorem 
\ref{lusztig restriction}.\fin

%\begin{coro}\label{restriction u}
%Let $w \in W_\Gb^\ci(\Lb,v)$, let $a \in H^1(F,A_\Gb(u))$ 
%and assume that $q \ge q_0$. Then 
%$$\lexp{*}{R}_{\Lb_w}^\Gb(\g_{u_a}^\Gb)=
%\g_{v_{w,\ph_\Lb^\Gb(w)\th(a)}}^{\Lb_w}.$$
%\end{coro}
%
%\proof By \cite[??]{DLM1} and \cite[??]{DLM2}, there exists an element 
%$b \in H^1(F,A_\Lb(v))$ such that, for every $a \in H^1(F,A_\Gb(u))$, 
%we have
%$$\lexp{*}{R}_{\Lb_w}^\Gb(\g_{u_a}^\Gb)=
%\g_{v_{w,b\th(a)}}^{\Lb_w}.$$
%It is then clear from Theorem \ref{lusztig restriction} that $b=\ph_\Lb^\Gb(w)$.\fin

\bi

\subsection{General Lusztig restriction} 
Let $v_\Mb^\pr=\pi_{\Lb_\Mb^\pr}(u)$. Since $\Lb_\Mb$ and $\Lb_\Mb^\pr$ are 
$\Gb$-conjugate, there exists an element $w_\Mb \in W_\Gb^\ci(\Lb_\Mb^\pr,v_\Mb^\pr)$ 
and an element $g_\Mb \in \Gb$ such that $\Lb_\Mb=\lexp{g_\Mb}{\Lb_\Mb^\pr}$ and 
$g_\Mb^{-1}F(g_\Mb)=\wdo_\Mb$. 
Now, let $z_\Mb$ denote the class of $\ph_{\Lb_\Mb^\pr,v}^\Gb(w_\Mb)$ in 
$H^1(F,\ZC(\Mb))$ (note that the element $\ph_{\Lb,v}^\Gb(w_\Mb)$ belongs to 
$\ZC(\Lb_\Mb^\pr)$, and we identify 
$\ZC(\Lb_\Mb^\pr)$ with $\ZC(\Mb)$). 
Finally, let $u_{\Mb,z}$ denote a representative of the regular unipotent 
class of $\Mb^F$ associated to $z \in H^1(F,\ZC(\Mb))$ via 
the Lang Theorem.

\bi

\begin{theo}\label{but dans la vie}
Assume that $p$ is good for $\Gb$ and that $q$ is large enough. Then 
$$\lexp{*}{R}_\Mb^\Gb(\g_u^\Gb)=\g_{u_{\Mb,z_\Mb}}^\Mb.$$
\end{theo}

\bi

%\remark{but remarque} We recall that the definition of friendly has been 
%given in Remark \ref{large enough}, and that $(\Gb,F)$ is friendly 
%whenever $p$ is almost good and $q$ is large enough.\finl
%
%\bi

\proof Since $p$ is good and $q$ is large enough, the Mackey formula holds in $\Gb$ by 
\cite[Theorem 6.1.1]{q grand}. So the functions $R_{\Mb'}^\Mb(f)$, 
where $\Mb'$ is an $F$-stable \levi of a \para of $\Gb$ and $f$ is 
a unipotently supported absolutely cuspidal functions on $\Mb^{\pr F}$, 
span the space of unipotently supported functions on $\Mb^F$. So, it 
is sufficient to prove that 
$$\langle \lexp{*}{R}_\Mb^\Gb(\g_u^\Gb),R_{\Mb'}^\Mb(f) \rangle_{\Mb^F} =
\langle \g_{u_{\Mb,z_\Mb}}^\Mb , R_{\Mb'}^\Mb(f) \rangle_{\Mb^F} \leqno{\mathrm{(*)}}$$
for all such pairs $(\Mb',f)$. By adjonction, $(*)$ is equivalent to
$$\langle \lexp{*}{R}_{\Mb'}^\Gb(\g_u^\Gb),f \rangle_{\Mb^{\pr F}}=
\langle \lexp{*}{R}_{\Mb'}^\Mb \g_{u_{\Mb,z_\Mb}}^\Mb ,f \rangle_{\Mb^{\pr F}}.
\leqno{\mathrm{(**)}}$$
But, by \cite[Theorem 1.14 (b)]{lugf}, a basis of the space 
of unipotently supported absolutely cuspidal functions on $\Mb^{\pr F}$ is given 
by the characteristic functions of $F$-stable cuspidal local systems on $F$-stable 
unipotent classes of $\Mb'$. So we may assume that $f$ is such a characteristic 
function. If the unipotent class which is the support of $f$ is not regular, 
then both sides of $(**)$ are equal to $0$. Therefore, we may assume that 
$\Mb'=\Lb_{w'}$ for some $w' \in W_\Gb^\ci(\Lb,v)$, and that 
$f=\gamh_{v_{w'},\z}^{\Lb_{w'}}$. 

First,
$$\g_u^\Gb=\sum_{\x \in H^1(F,\ZC(\Gb))^\we} \gamh_{u,\x}^\Gb$$
and, for instance by \cite[1.8.3 and Lemma 1.8.4 (b)]{bonnafe A}, 
$$\langle \lexp{*}{R}_{\Lb_{w'}}^\Gb(\gamh_{u,\x}^\Gb),\gamh_{v_{w'},\z}^{\Lb_{w'}} 
\rangle_{\Lb_{w'}^F}=0$$
if $\x \not= \zett$. Therefore, by Theorem \ref{lusztig restriction}, 
$$\langle \lexp{*}{R}_{\Lb_{w'}}^\Gb(\g_u^\Gb),\gamh_{v_{w'},\z}^{\Lb_{w'}} 
\rangle_{\Lb_{w'}^F}= \z(\ph_{\Lb,v}^\Gb(w')) 
\langle \g_{v_{w'},\z}^{\Lb_{w'}},\g_{v_{w'},\z}^{\Lb_{w'}} \rangle_{\Lb_{w'}^F}.
\leqno{(A)}$$

We consider now the Lusztig restriction from $\Mb$ to $\Lb_w$. First, note 
that 
$$\g_{u_{\Mb,z_\Mb}}^\Gb =\sum_{\x \in H^1(F,\ZC(\Mb))^\we} 
\x(z_\Mb) \gamh_{u_\Mb,\x}^\Mb.$$
Let 
$\Mb_0$ be a \levi of a \para of $\Gb$ conjugate to $\Mb$.  
By Proposition \ref{propriete distingues} (b), $W_{\Mb_0}^\ci(\Lb,v)$ is 
a standard parabolic subgroup of $W_\Gb^\ci(\Lb,v)$. Therefore, we may 
write $w'=w_1w_2$ where $w_1 \in W_{\Mb_0}^\ci(\Lb,v)$ and $w_2$ is the image 
of $w_\Mb$ under the morphism $\rhoba_{\Lb,\Lb',u}$ defined in 
\SEC\ref{compatible sub}. Recall that 
$w_\Mb \in W_\Gb^\ci(\Lb,v)$ has been defined in the introduction of this 
subsection.

By \ref{44}, by Theorem 
\ref{lusztig restriction} applied to $\Mb$, and by 
\cite[1.8.3 and Lemma 1.8.4 (b)]{bonnafe A}, we have
$$\lexp{*}{R}_{\Lb_w}^\Mb(\gamh_{u_\Mb,\x}^\Mb)=\left\{\begin{array}{ll}
\z(\ph_{\Lb,v}^{\Mb_0}(w_1))\gamh_{v_{w'},\z}^{\Lb_{w'}} & {\mathrm{if}}~\x=\z 
\ci h_\Lb^{\Mb_0} \\
0 & {\mathrm{otherwise}}.\end{array}\right.$$
So, by Proposition \ref{compatible prop}, 
$$\langle \lexp{*}{R}_{\Lb_w}^\Mb(\g_{u_{\Mb,z_\Mb}}^\Mb) , \gamh_{v_{w'},\z}^{\Lb_{w'}} 
\rangle_{\Lb_{w'}^F} = \z(\ph_{\Lb,v}^{\Mb_0}(w_1))\z(\ph_{\Lb,v}^\Gb(w_2)) 
\langle \gamh_{v_{w'},\z}^{\Lb_{w'}}, 
\gamh_{v_{w'},\z}^{\Lb_{w'}} \rangle_{\Lb_{w'}^F}.$$
Therefore, by Corollary \ref{coro compa para}, 
$$\langle \lexp{*}{R}_{\Lb_w}^\Mb(\g_{u_{\Mb,z_\Mb}}^\Mb) , \gamh_{v_{w'},\z}^{\Lb_{w'}} 
\rangle_{\Lb_{w'}^F} = \z(\ph_{\Lb,v}^\Gb(w')) 
\langle \gamh_{v_{w'},\z}^{\Lb_{w'}}, 
\gamh_{v_{w'},\z}^{\Lb_{w'}} \rangle_{\Lb_{w'}^F}.\leqno{(B)}$$
By comparing $(A)$ and $(B)$, we get $(**)$, so Theorem \ref{but dans la vie} 
is proved.\fin

\newpage

%%%%%%%%%%%%%%%%%%
%%%%%%%%%%%%%%%%%%
%%%%%%%%%%%%%%%%%%
%%%%%%%%%%%%%%%%%%
%%%%%%%%%%%%%%%%%%
%%%%%%%%%%%%%%%%%%
%%%%%%%%%%%%%%%%%%
%%%%%%%%%%%%%%%%%%
%%%%%%%%%%%%%%%%%%
%%%%%%%%%%%%%%%%%%
%%%%%%%%%%%%%%%%%%
%%%%%%%%%%%%%%%%%%
%%%%%%%%%%%%%%%%%%
%%%%%%%%%%%%%%%%%%
%%%%%%%%%%%%%%%%%%
%%%%%%%%%%%%%%%%%%
%%%%%%%%%%%%%%%%%%
%%%%%%%%%%%%%%%%%%
%%%%%%%%%%%%%%%%%%
%%%%%%%%%%%%%%%%%%
%%%%%%%%%%%%%%%%%%
%%%%%%%%%%%%%%%%%%

\bigskip

\bigskip

{\sc C. Bonnaf\'e}

CNRS - UMR 6623, 

Laboratoire de Math\'ematiques de Besan\c{c}on, 

16 Route de Gray, 

25030 BESAN\c{C}ON Cedex, 

FRANCE, 

{\tt bonnafe\arobas math.univ-fcomte.fr}

\end{document}